\documentclass{m2an}
\usepackage{amssymb}
\usepackage{amsmath, amsthm}
\usepackage{hyperref, cleveref}
\usepackage{stmaryrd}
\usepackage{siunitx}
\usepackage{commath}
\usepackage{enumitem}
\usepackage[caption=false]{subfig}
\usepackage{algorithm, algorithmicx}
\usepackage{algpseudocode}
\usepackage{euscript,amscd,verbatim,amsbsy}

 \usepackage{array,multirow,graphicx}
 \usepackage{float}
\usepackage{xspace,color}

\newcommand{\rchi}{\raisebox{2pt}{$\chi$}}

\makeatletter
\newcommand{\tnorm}{\@ifstar\@tnorms\@tnorm}
 \newcommand{\@tnorms}[1]{%
 }
\newcommand{\@tnorm}[2][]{%
  \mathopen{#1|\mkern-1.5mu#1|\mkern-1.5mu#1|}
  #2
  \mathclose{#1|\mkern-1.5mu#1|\mkern-1.5mu#1|}
}
\makeatother

\newcommand{\jump}[1]{\llbracket #1 \rrbracket}

\begin{document}
\title{A hybridizable discontinuous Galerkin method for the fully coupled time-dependent Stokes/Darcy--transport problem}
%
\author{Aycil Cesmelioglu}\address{Oakland University, Department of Mathematics and Statistics, Rochester Hills, MI, 48307, USA, \email{cesmelio@oakland.edu\ \&\ ddpham@oakland.edu}}
\author{Dinh Dong Pham}\sameaddress{1}
\author{Sander Rhebergen}\address{University of Waterloo, Department of Applied Mathematics, Waterloo, ON, Canada, \email{srheberg@uwaterloo.ca}}
%
%
\begin{abstract} 
    We present a high-order hybridized discontinuous Galerkin (HDG) method for the fully coupled time-dependent Stokes–Darcy-transport problem where the fluid viscosity and source/sink terms depend on the concentration and the dispersion/diffusion tensor depends on the fluid velocity. This HDG method is such that the discrete flow equations are compatible with the discrete transport equation. Furthermore, the HDG method guarantees strong mass conservation in the $H^{\rm div}$ sense and naturally treats the interface conditions between the Stokes and Darcy regions via facet variables. We employ a linearizing decoupling strategy where the Stokes/Darcy and the transport equations are solved sequentially by time-lagging the concentration. We prove well-posedness and optimal a priori error estimates for the velocity and the concentration in the energy norm. We present numerical examples that respect compatibility of the flow and transport discretizations and demonstrate that the discrete solution is robust with respect to the problem parameters.
\end{abstract}
%
%
\maketitle
\section{Introduction}
\label{sec:introduction}
Coupled free fluid and porous media flow is encountered in many engineering applications \cite{Discacciati:thesis-2004,Hanspal:2006} and can be modeled by the Stokes/Darcy equations. Adding a transport equation to this coupled system brings forth a model that can be used to simulate the spread of contaminants towards groundwater resources \cite{Bear:book} or biochemical transport in hemodynamics \cite{Dangelo:2011}.

The accuracy and stability of numerical discretizations of the stationary Stokes/Darcy equations \cite{Discacciati:2002,Layton:2002,Burman:2005,Cao:2010,Camano:2015,Badia:2009,Gatica:2009,Marquez:2015,Riviere:2005} and advection-diffusion type transport equations \cite{Nguyen:2009,Brooks:1982,Cockburn:1998,Wells:2011} are well studied. However, accuracy and stability are not automatically guaranteed when these discretizations are coupled.
In particular,
\emph{compatible} discretizations, as defined by \cite{Dawson:2004}, are desired to avoid loss of accuracy and loss of conservation properties of the numerical methods used for the transport equation.

The first numerical study on the coupling of the stationary Stokes/Darcy equations with a transport equation was given in \cite{Vassilev:2009} where a mixed finite element method (MFEM) is used for the flow problem and the local discontinuous Galerkin method is used for the transport problem. They considered one-way coupling; the concentration is affected by the flow velocity, but the velocity is not affected by a change in concentration.
The same problem was studied in \cite{Riviere:2014} by using discontinuous Galerkin (DG) methods for both flow and transport equations. In \cite{Ervin:2019}, Ervin et. al considered a fully time-dependent version of the one-way coupled problem where they developed partitioned time-stepping methods by imposing the interface conditions weakly using penalties. One-way coupling was considered also in \cite{Cesmelioglu:2021} in which the flow problem was discretized by a strongly mass conservative Embedded-Hybridized DG (EDG-HDG) method while the transport equation was discretized by an EDG method.

Less studied is the fully-coupled problem in which, apart from the transport equation depending on the flow velocity, the flow solution is time-dependent and the fluid viscosity and source/sink terms depend on the concentration. 
To the best of our knowledge, there are only two papers that focus on this fully coupled problem. First, \cite{Cesmelioglu:2012} presented an analysis of a weak solution for the case where free flow is governed by the Navier--Stokes equations. The analysis in \cite{Cesmelioglu:2012}, however, also holds when free flow is governed by the Stokes equations. The only numerical paper on this topic, \cite{Rui:2017}, introduced a stabilized mixed finite element method using nonconforming piece-wise linear Crouzeix–Raviart finite elements for the velocity, a piece-wise constant approximation for the pressure, and a conforming, piece-wise linear, finite element method for the transport equation based on a skew-symmetric formulation.

In this paper, we extend the work in \cite{Cesmelioglu:2021} to the fully coupled case. To deal with the non-linearity, we consider a linearizing decoupling strategy, where the Stokes/Darcy and the transport equations are solved sequentially by time-lagging the concentration. We use HDG methods \cite{Cockburn:2009a} for both the Stokes/Darcy and transport sub-problems at each time step and prove well-posedness and a priori error estimates. These results can easily be extended to the EDG-HDG discretization used for the Stokes/Darcy problem in \cite{Cesmelioglu:2021}. Our HDG method for the flow problem provides the transport sub-problem at each time step with an exactly mass conserving and $H({\rm div})$-conforming velocity field. This renders our scheme robust with respect to the problem parameters. By choosing the polynomial degree in a specific way our flow/transport scheme is also compatible.

Here is an outline for the remainder of this article. 
In \Cref{sec:stokesdarcytransport}, we present the fully coupled Stokes/Darcy-transport model and specify the assumptions on the problem parameters. \Cref{sec:hdgstokesdarcytransport} sets notation, describes in detail the semi-discrete HDG scheme, and lists the attractive properties of the numerical discretization.  Next, \Cref{sec:properties} summarizes standard inequalities and shows continuity, coercivity, and the inf-sup condition for the discretization of the Stokes/Darcy sub-problem. A full discretization of the problem based on a sequential decoupling strategy is introduced in \Cref{sec:sequential_strategy} while the main results, i.e., a priori error estimates for the velocity, pressure, and concentration, are presented in \Cref{sec:main_results}. Finally, we present some numerical experiments in \Cref{sec:numerical_examples} followed by conclusions in \Cref{sec:conclusions}.

\section{The Stokes/Darcy--transport system}
\label{sec:stokesdarcytransport}
Let $\Omega \subset \mathbb{R}^{\dim}$, ${\dim}=2,3$, be a bounded polygonal
domain.
\begin{figure}[!ht]
  \centering
  \includegraphics[width=0.25\textwidth]{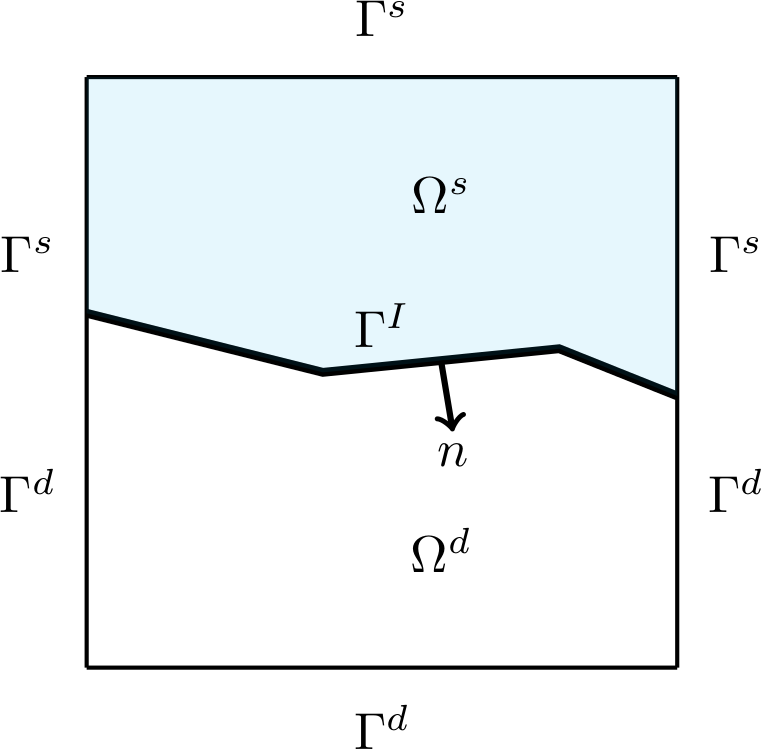}
  \caption{A depiction of a domain $\Omega\subset \mathbb{R}^2$ with its two sub-domains $\Omega^s$ and $\Omega^d$.}
  \label{fig:domaindescription}
\end{figure}
We denote its boundary by $\partial\Omega$ and the outward
unit normal to $\partial \Omega$ by $n$. Domain $\Omega$ consists of two
non-overlapping polygonal regions, a free flow region $\Omega^s$ and a Darcy flow region $\Omega^d$, such that $\Omega = \Omega^s \cup \Omega^d$. The
polygonal interface between $\Omega^s$ and $\Omega^d$ is denoted by $\Gamma^I$
and the external boundary of $\Omega^j$ is denoted by
$\Gamma^j := \partial \Omega \cap \partial \Omega^j$, $j=s,d$. See
\Cref{fig:domaindescription} for a depiction of a domain when ${\dim}=2$. 

We denote the time interval of interest by $J = [0,T]$.
The fully coupled Stokes/Darcy--transport system for the velocity field $u :\Omega \times J \to
\mathbb{R}^{\dim}$, fluid pressure $p : \Omega \times J \to \mathbb{R}$
and concentration
$c:\Omega\times J \rightarrow \mathbb{R}$ is given
by
\begin{subequations}
  \label{eq:system}
  \begin{align}
    \label{eq:momentum}
    \partial_t u-\nabla\cdot (2\mu(c)\varepsilon(u)) + \nabla p &= f^s(c) & &
    \text{in}\ \Omega^s \times J,
    \\
    \label{eq:d_velocity}
    \mathbb{K}^{-1}(c)u + \nabla p &=\mathbb{K}^{-1}(c)f^d(c) & & \text{in}\ \Omega^d \times J,
    \\
    \label{eq:mass}
    - \nabla\cdot u &= \rchi_d (g_{p} - g_{i}) & & \text{in}\ \Omega\times J,
    \\
    \label{eq:concentration}
    \phi\partial_t c + \nabla \cdot(cu - \widetilde{D}(u)\nabla c) &=
    \rchi_d (c_{I} g_{i} - c g_p) & & \text{in}\ \Omega\times J,
    \\
    \label{eq:bc_s}
    u &= 0 & & \text{on}\ \Gamma^s\times J,
    \\
    \label{eq:bc_d}
    u\cdot n &= 0 & & \text{on}\ \Gamma^d\times J,
    \\
    \label{eq:bc_c}
    \widetilde{D}(u)\nabla c\cdot n &= 0 & & \text{on}\ \partial \Omega\times J,
  \end{align}
\end{subequations}
where $\varepsilon(u) := \del{\nabla u + (\nabla u)^T}/2$ is the strain
rate tensor and $\rchi_d$ is the characteristic function that takes the
value 1 in $\Omega^d$ and 0 in $\Omega^s$.
Here the fluid viscosity $\mu$, the matrix $\mathbb{K}=\frac{\kappa}{\mu}$, where $\kappa$ is the permeability matrix of the porous medium, and the body force terms $f^s$ and $f^d$ are concentration dependent functions. The porosity $\phi$ of the medium in $\Omega^d$ is a spatially varying function. In $\Omega^s$ we set $\phi = 1$. The functions $g_{i}$ and $g_{p}$ denote the source and sink terms related to injection and production wells and $c_{I}$ is the injected concentration.
Furthermore, in the Stokes region $\widetilde{D}(u)=d\mathbb{I}$, where $\mathbb{I}$ is the $\dim\times\dim$ identity matrix and $d$ is the
diffusion coefficient. In the Darcy region $\widetilde{D}(u)=D(u)$, where $D(u)$ denotes the diffusion dispersion tensor in $\Omega^d$.

We will denote the restriction of the velocity $u$, pressure $p$, and concentration $c$ to sub-domain $\Omega^j$, $j=s,d$ by, respectively, $u^j$, $p^j$, and $c^j$. Then, on the interface $\Gamma^I$, choosing the unit normal vector $n$ to be pointing from $\Omega^s$ to $\Omega^d$, we prescribe the following interface conditions that hold for $t \in J$:
\begin{subequations}
  \label{eq:interface}
  \begin{align}
    \label{eq:bc_I_u}
    u^s\cdot n &= u^d\cdot n,
    \\
    \label{eq:bc_I_p}
    p^s - (2\mu(c^s)\varepsilon(u^s)n)\cdot n &= p^d,
    \\
    \label{eq:bc_I_slip}
    -2(\varepsilon(u^s)n)\cdot \tau^{\ell}&=\gamma^{\ell} u^s\cdot \tau^{\ell}, 
   \quad \ell = 1, \hdots, {\dim}-1,\\
    \label{eq:bc_I_c}
    c^s &= c^d,
    \\
    \label{eq:bc_I_c_flux}
    d\nabla c^s\cdot n &= D(u^d)\nabla c^d\cdot n,
  \end{align}
\end{subequations}
where $\tau^{\ell}$, $\ell=1, \hdots, {\dim}-1$ denote the unit tangent vectors on $\Gamma^I$.
These conditions enforce the normal continuity of the velocity \eqref{eq:bc_I_u},
the normal continuity of the normal component of the stress \eqref{eq:bc_I_p}, continuity of the concentration \eqref{eq:bc_I_c}, and normal continuity of the concentration flux \eqref{eq:bc_I_c_flux}. \Cref{eq:bc_I_slip}, where  $\gamma^{\ell}=\alpha/\sqrt{\tau^{\ell}\cdot \kappa \tau^{\ell}}$ with $\alpha > 0$ a constant, is the Beavers--Joseph--Saffman law which enforces a condition on the tangential component of the normal stress~\cite{Beavers:1967, Saffman:1971}.

To close the model, we assume the following initial conditions:
\begin{subequations}
  \label{eq:initial}
  \begin{align}
    \label{eq:u_init}
    u^s(x,0)&=u^s_0(x) & & \text{in}\ \Omega^s,
    \\
    \label{eq:c_init}
    c(x,0)&=c_0(x) & & \text{in}\ \Omega.
  \end{align}
\end{subequations}

We end this section by discussing some assumptions we make on the various functions used in the Stokes/Darcy--transport model.
The dispersion-diffusion tensor $D(u)$ in $\Omega^d$ satisfies for $u, v \in \mathbb{R}^{\dim}$:
\begin{subequations}
  \begin{align}
  \label{eq:D_min}
    D_{\min} |\xi|^2 &\le \xi^T D(u) \xi \quad \forall \xi \in \mathbb{R}^{\dim},
    \\
    \label{eq:upperboundDuabs}
    |D(u)| &\le C(1 + \envert{u}),
    \\
    \label{eq:DLipschitz}
    \envert{ D(u) - D(v) } &\le C \envert{u - v},
  \end{align}
\end{subequations}
where $D_{\min}$ and $C$ are positive constants and $|\cdot|$  denotes the Euclidean norm.
We assume that $\mu$ is Lipschitz continuous in $c$ with Lipschitz constant $\mu_L$ 
and that there exist constants $\phi_*, \phi^*, \mu_*, \mu^* > 0$ such that
\begin{subequations}
  \begin{align}
  \label{eq:assumption_phi}
  \phi_*\leq \phi(x)&\leq \phi^* && \forall x\in \Omega^d,
  \\
  \label{eq:assumption_mu}
  \mu_*\leq \mu(c)&\leq \mu^* && \forall c\in \mathbb{R}.
  \end{align}
\end{subequations}
%
The permeability matrix $\kappa$ is symmetric, uniformly bounded, and elliptic, that is, there exist positive constants $\kappa_{*}< \kappa^{*}$ such that
\begin{equation}
  \label{eq:assumption_kappa}
  \kappa_* |\xi|^2\leq \xi^T\kappa(x) \xi\leq \kappa^* |\xi|^2\quad \forall  \xi\in \mathbb{R}^{\dim}, \quad \forall x\in \bar{\Omega}^d.
\end{equation}
From \eqref{eq:assumption_kappa} and \eqref{eq:assumption_mu}, we deduce that
\begin{equation}\label{eq:assumption_K}
  K_{*}|\xi|^2\leq \xi^T \mathbb{K}(c, x) \xi\leq K^{*} |\xi|^2\quad \forall  \xi\in \mathbb{R}^{\dim},\quad  \forall (c,x)\in \mathbb{R}\times  \bar{\Omega}^d,
\end{equation}
where $K_{*}=\kappa_{*}/\mu^{*}$ and $K^{*}=\kappa^{*}/\mu_{*}$.

The body force functions $f^s$ and $f^d$ are assumed to be Lipschitz continuous in $c$ with Lipschitz constants $L_f^s$ and $L_f^d$. 
Note that $f^s$ and $f^d$ depend on $x$ and $c$, but they do not depend explicitly on $t$.
We will further assume that $0\leq c_I \leq 1$ a.e. in $\Omega^d$ and that $g_{i}, g_{p}\geq 0$, $g_{i}, g_{p}\in L^{\infty}(J; L^2(\Omega^d))$ are such that
\begin{equation*}
  \int_{\Omega^d}(g_{i}(x, t)-g_{p}(x, t)) \dif x= 0 \quad \forall t\in J.  
\end{equation*}
A weak formulation of the problem defined by \cref{eq:system,eq:interface,eq:initial} was presented in \cite{Rui:2017}. The analysis for a weak solution of a more general version of this problem, in which the free fluid flow is governed by the Navier--Stokes equations, can be found in \cite{Cesmelioglu:2012}.

\section{The hybridized discontinuous Galerkin method}
\label{sec:hdgstokesdarcytransport}

\subsection{Preliminaries}
We use the same notation that we used previously in \cite{Cesmelioglu:2020,Cesmelioglu:2021}. 
Let $\mathcal{T}^j := \cbr{K}$ be a shape-regular triangulation of
$\Omega^j$, $j = s, d$, into non-overlapping elements (we only consider simplices) such that $\mathcal{T}^s$ and $\mathcal{T}^d$ match at the interface~$\Gamma^I$. We define the triangulation of the entire domain $\Omega$ as $\mathcal{T} := \mathcal{T}^s \cup \mathcal{T}^d$. The maximum diameter
over all elements is $h= \max_{K\in \mathcal{T}} h_K$, where $h_K$ stands for the diameter of an element $K$. The boundary of an element $K$ and its outward unit normal are denoted by $\partial K$ and $n$, respectively. 
A facet of an element boundary is an interior facet if it is shared by two neighboring elements and it is a boundary facet if it is a part of $\partial \Omega$. 
The set of all interior facets and all boundary facets in $\bar{\Omega}^j$ are denoted by $\mathcal{F}_i^j$ and $\mathcal{F}_b^j$, $j=s,d$, respectively. We also collect the facets that lie on the interface $\Gamma^I$ in the set $\mathcal{F}^I$. 
The set of all facets that lie in $\bar{\Omega}$ and in $\bar{\Omega}^j$ are denoted by $\mathcal{F}$ and $\mathcal{F}^j$, respectively. We point out that $\mathcal{F}^j=\mathcal{F}_i^j\cup\mathcal{F}_b^j\cup\mathcal{F}^I$, $j=s,d$. Furthermore, we define $\Gamma_0:=\cup_{F\in \mathcal{F}}F$ and $\Gamma_0^j:=\cup_{F\in \mathcal{F}^j}F$, $j = s, d$.  

The finite element function spaces on~$\Omega$ for the velocity and pressure are given by
\begin{equation}
  \label{eq:DGcellspaces}
  \begin{split}
    V_h &:= \cbr[1]{v_h \in \sbr[0]{L^2(\Omega)}^{\dim} : \ v_h \in
      \sbr[0]{P_{k_f}(K)}^{\dim}, \ \forall\ K \in \mathcal{T}},
    \\
    Q_h &:= \cbr[1]{q_h \in L^2(\Omega) : \ q_h \in P_{k_f-1}(K) ,\
      \forall \ K \in \mathcal{T}}\cap L^2_0(\Omega),
  \end{split}
\end{equation}
where $P_k(K)$ denotes the space of polynomials of degree at most $k$ defined on the element $K$. The finite element spaces for the velocity and pressure traces are given by
\begin{equation}
  \label{eq:DGfacetspaces}
  \begin{split}
    \bar{V}_h &:= \cbr[1]{\bar{v}_h \in \sbr[0]{L^2(\Gamma_0^s)}^{\dim}:\
      \bar{v}_h \in \sbr[0]{P_{k_f}(F)}^{\dim}\ \forall\ F \in \mathcal{F}^s,\
      \bar{v}_h = 0 \ \mbox{on}\ \Gamma^s},\\
    \bar{Q}_h^j &:= \cbr[1]{\bar{q}_h^j \in L^2(\Gamma_0^j) : \ \bar{q}_h^j
      \in P_{k_f}(F) \ \forall\ F \in \mathcal{F}^j},\quad j=s,d.
  \end{split}
\end{equation}
Here $P_k(F)$ denotes the space of polynomials of degree at most $k$ defined on the facet $F$. Note that functions in $\bar{V}_h$ are not defined on $\Gamma_0^d\backslash\Gamma_I$. The finite element function spaces for the concentration and its trace are defined as
\begin{equation}
  \label{eq:spaces_transport}
  \begin{split}
    C_h = \cbr[0]{c_h\in L^2(\Omega) : \ c_h \in P_{k_c}(K) ,\
      \forall \ K \in \mathcal{T}},
\quad
    \bar{C}_h = \cbr[0]{\bar{c}_h \in L^2(\Gamma_0) : \ \bar{c}_h \in
      P_{k_c}(F) \ \forall\ F \in \mathcal{F}}.
  \end{split}
\end{equation}
The semi-discrete and fully-discrete HDG methods for the flow and transport equations considered in this article are compatible when $k_c = k_f - 1$ \cite{Cesmelioglu:2021}. For this reason we set $k_c=k_{f}-1$.

To reduce the notational burden, we define $\boldsymbol{V}_h := V_h \times \bar{V}_h$,
$\boldsymbol{Q}_h := Q_h \times \bar{Q}_h^s \times \bar{Q}_h^d$, and
$\boldsymbol{Q}_h^{j} := Q_h^j \times \bar{Q}_h^j$, $j=s,d$. We denote elements in these product spaces by $\boldsymbol{v}_h := (v_h, \bar{v}_h) \in \boldsymbol{V}_h$, $\boldsymbol{q}_h := (q_h, \bar{q}_h^s, \bar{q}_h^d) \in
\boldsymbol{Q}_h$, and $\boldsymbol{q}_h^{j} := (q_h^j, \bar{q}_h^j) \in
\boldsymbol{Q}_h^{j}$, $j=s,d$. In addition, we set $\boldsymbol{X}_h := \boldsymbol{V}_h \times \boldsymbol{Q}_h$.
Similarly, we introduce $\boldsymbol{C}_h = C_h \times \bar{C}_h$ and
denote the corresponding elements by $\boldsymbol{c}_h := (c_h,\bar{c}_h) \in \boldsymbol{C}_h$.

Next, let us define the function spaces
\begin{align*}
    V &:= \big\{v\in [L^2(\Omega)]^{\dim}\ :\ v^s\in \sbr[0]{H^2(\Omega^s)}^{\dim},\
      v^d \in \sbr[0]{H^1(\Omega^d)}^{\dim}, \\
      &\hspace{3.5cm} v=0\ \text{on}\ \Gamma^s,\ v \cdot n=0 \ \text{on}\ \Gamma^d,
      \ v^s\cdot n = v^d\cdot n\ \text{on}\ \Gamma^I\big\}, \nonumber
    \\
    Q &:= \big\{q \in L_0^2(\Omega)\ :\ q^s \in H^1(\Omega^s),\ q^d \in H^2(\Omega^d)\big\},\\
    C &:= H^2(\Omega), 
\end{align*}
and set $X := V \times Q$. As before, we use a superscript $^j$ to specify the restriction of these spaces to $\Omega^j$, $j=s,d$. The trace spaces of $V$
restricted to $\Gamma_0^s$, $Q^j$ restricted to $\Gamma_0^j$, and $C$ restricted to $\Gamma_0$ are denoted by, respectively, $\bar{V}$, $\bar{Q}^j$, and $\bar{C}$. The trace operator
$\gamma_V : V^s \to \bar{V}$ restricts functions in $V^s$ to
$\Gamma_0^s$, and similarly the trace operators
$\gamma_{Q^j} : Q^j \to \bar{Q}^j$ restrict functions in $Q^j$ to
$\Gamma_0^j$, $j=s,d$. However, when it is clear from the context, we omit the
subscript in the trace operator. Analogous to the discrete case, we
 introduce $\boldsymbol{V} := V\times \bar{V}$,
$\boldsymbol{Q} := Q\times \bar{Q}^s \times \bar{Q}^d$, and $\boldsymbol{C} := C\times \bar{C}$. We then define extended function spaces as
\begin{equation*}
  \boldsymbol{V}(h) := \boldsymbol{V}_h + \boldsymbol{V}, \quad
  \boldsymbol{Q}(h) := \boldsymbol{Q}_h + \boldsymbol{Q}, \quad
  \boldsymbol{C}(h) := \boldsymbol{C}_h + \boldsymbol{C}, 
\end{equation*}
and set $\boldsymbol{X}(h) := \boldsymbol{V}(h) \times \boldsymbol{Q}(h)$.

We close this section by listing various norms on the spaces described above. We refer the reader to \cite{Adams:book} for the definitions of the standard Sobolev spaces $W^{m,p}(D)$ and their corresponding norms $\|\cdot\|_{W^{m,p}(D)}$. For ease of notation, we write $\|\cdot\|_{m,p,D}$ instead of $\|\cdot\|_{W^{m,p}(D)}$ with the following simplifications. When $m=0$, $W^{0,p}(D)$ coincides with $L^p(D)$ and when $p=2$, $H^m(D)=W^{m,p}(D)$. For $p=2$, we write $\|\cdot\|_{m,D}$ to denote $\|\cdot\|_{W^{m,2}(D)}$ and for $m=0$, $p=2$, we write $\|\cdot\|_{D}$ instead of $\|\cdot\|_{0,D}$.

On $\boldsymbol{V}^s(h)$ we define the standard HDG-norm and its strengthened version as follows:
\begin{equation*}
  \tnorm{\boldsymbol{v}}_{v,s}^2 := \sum_{K\in \mathcal{T}^s}\del[1]{\norm[0]{\nabla v}_K^2
  + h_K^{-1}\norm{v-\bar{v}}^2_{\partial K}},
\quad  \tnorm{\boldsymbol{v}}_{v',s}^2 := \tnorm{\boldsymbol{v}}_{v,s}^2
  + \sum_{K\in \mathcal{T}^s} h_K^2 \envert{v}_{H^2(K)}^2.
\end{equation*}
On $\boldsymbol{V}(h)$ we then introduce the norms
\begin{equation*}
    \begin{split}
    \tnorm{\boldsymbol{v}}_{v}^2
    &:=
    \tnorm{\boldsymbol{v}}_{v,s}^2
    + \norm{v}_{\Omega^d}^2
    +  \sum_{j=1}^{\dim}\gamma^j\| \bar{v}\cdot \tau_j \|^2_{\Gamma^I}, 
    \\
  \tnorm{\boldsymbol{v}}_{v'}^2
  &:= \tnorm{\boldsymbol{v}}_{v}^2
  + \sum_{K\in \mathcal{T}^s} h_K^2 \envert{v}_{H^2(K)}^2
  =\tnorm{\boldsymbol{v}}_{v',s}^2  + \norm{v}_{\Omega^d}^2
  + \sum_{j=1}^{\dim}\gamma^j\| \bar{v}\cdot \tau_j \|^2_{\Gamma^I},
  \end{split}
\end{equation*}
and note that $\tnorm{\boldsymbol{v}}_{v}$ and $\tnorm{\boldsymbol{v}}_{v'}$ are equivalent on $\boldsymbol{V}_h$ due to the fact that $\tnorm{\cdot}_{v,s}$ and $\tnorm{\cdot}_{v',s}$ are equivalent on $\boldsymbol{V}_h^s$ (see, for example, \cite[eq.~(5.5)]{Wells:2011}).

On the pressure spaces $\boldsymbol{Q}^{j}(h)$, $j = s,d$ and $\boldsymbol{Q}(h)$, we define, respectively,
\begin{equation*}
  \tnorm{\boldsymbol{q}^j}_{p,j}^2 := \norm{q}_{\Omega^j}^2
  + \sum_{K\in \mathcal{T}^j} h_K \norm[0]{\bar{q}^j}_{\partial K}^2,
  \qquad   \tnorm{\boldsymbol{q}}_{p}^2
  := \sum_{j=s,d}\tnorm{\boldsymbol{q}^j}^2_{p,j}.
\end{equation*}

Finally, for $\boldsymbol{w}_h\in \boldsymbol{C}(h)$, we define the following semi-norm:
\begin{equation}
    \tnorm{\boldsymbol{w}_h}_c^2=\sum\limits_{K\in\mathcal{T}}(\|\nabla w_h\|_K^2+h_K^{-1}\|w_h-\bar{w}_h\|_{\partial K}^2).
\end{equation}

\subsection{Semi-discrete HDG scheme}
\label{sec:semi-discrete}
The semi-discrete method we propose for the Stokes/Darcy--transport system in
\cref{eq:system,eq:interface} is as follows: for $t>0$, find
$\big(\boldsymbol{u}_h(t), \boldsymbol{p}_h(t)\big)\in \boldsymbol{X}_h$ and $\boldsymbol{c}_h(t) \in \boldsymbol{C}_h$ such that
\begin{subequations}
\label{eq:semidiscreteHDG}
\begin{multline}
  \label{eq:hdgsd}
   \sum_{K\in \mathcal{T}^s}\int_{K}\partial_t u_h \cdot v_h \dif x  +  B_h^{sd}(\boldsymbol{c}_h; (\boldsymbol{u}_h, \boldsymbol{p}_h), (\boldsymbol{v}_h, \boldsymbol{q}_h) )
   = \sum_{K\in \mathcal{T}^s}\int_{K} f^s(c_h)\cdot v_h \dif x 
   \\
   +\sum_{K\in \mathcal{T}^d}\int_{K} \mathbb{K}^{-1}(c_h)f^d(c_h)\cdot v_h\dif x  
   + \sum_{K\in \mathcal{T}^d}\int_K (g_p-g_i)\,q_h \dif x
\end{multline}
and
\begin{equation}
  \label{eq:hdgtr}
  \sum_{K\in \mathcal{T}}\int_K  \phi\, \partial_tc_h w_h \dif x 
  + B_h^{tr}(u_h; \boldsymbol{c}_h, \boldsymbol{w}_h) +  \sum_{K\in \mathcal{T}^d} \int_K c_h\,g_{p}\, w_h \dif x
   =\sum_{K\in \mathcal{T}^d} \int_K c_{I}\,g_{i}\, w_h \dif x,    
\end{equation}
\end{subequations}
for all $\big(\boldsymbol{v}_h, \boldsymbol{q}_h\big)\in \boldsymbol{X}_h$ and $\boldsymbol{w}_h \in \boldsymbol{C}_h$.
 
The form $B_h^{sd}$ in \cref{eq:hdgsd} collects the discretization terms for the Stokes/Darcy momentum and mass conservation equations as follows:
\begin{multline}
  \label{eq:Bhsd}
  B_h^{sd}(\boldsymbol{c};  (\boldsymbol{u}, \boldsymbol{p}), (\boldsymbol{v}, \boldsymbol{q}) )
    :=
    a_h(\boldsymbol{c};  \boldsymbol{u}, \boldsymbol{v})
    + \sum_{j=s,d} \del{b_h^j(\boldsymbol{p}^j, v) + b_h^{I,j}(\bar{p}^j, \bar{v})}
    + \sum_{j=s,d} \del{b_h^j(\boldsymbol{q}^j, u) + b_h^{I,j}(\bar{q}^j, \bar{u})}.
\end{multline}
Here $a_h(\cdot, \cdot)$ is defined as
\begin{equation}
  \label{eq:def_ah}
  a_h(\boldsymbol{c};  \boldsymbol{u}, \boldsymbol{v})
  := a_h^s(c;  \boldsymbol{u}, \boldsymbol{v})
  + a_h^d(c;  u, v) + a_h^I(\bar{c};  \bar{u}, \bar{v}),
\end{equation}
where
\begin{align*}
    a_h^s(c;\boldsymbol{u}, \boldsymbol{v}):=&\sum_{K\in\mathcal{T}^s} \int_K 2\mu(c) \varepsilon(u) : \varepsilon(v) \dif x
    +\sum_{K\in\mathcal{T}^s} \int_{\partial K} \frac{2\beta_s\mu(c)}{h_K}(u-\bar{u}) \cdot (v-\bar{v}) \dif s
    \\
    &-\sum_{K\in\mathcal{T}^s} \int_{\partial K} 2\mu(c)\varepsilon(u)n^s \cdot (v-\bar{v})\dif s
    -\sum_{K\in\mathcal{T}^s} \int_{\partial K} 2\mu(c)\varepsilon(v)n^s \cdot (u-\bar{u})\dif s,\\\
    a_h^d(c;u, v)
    :=& \int_{\Omega^d} \mathbb{K}^{-1}(c) u\cdot v \dif x ,
    \\
    a_h^I(\bar{c};\bar{u}, \bar{v}) :=&  \sum_{\ell=1}^{\rm{\dim}-1}\int_{\Gamma^I}  \gamma^\ell \mu(\bar{c})(\bar{u}\cdot \tau^\ell)(\bar{v}\cdot \tau^\ell)\dif s,
  \end{align*}
and $\beta_s > 0$ is a penalty parameter. The bilinear forms
$b_h^j(\cdot, \cdot)$ and $b_h^{I,j}(\cdot, \cdot)$ in \cref{eq:Bhsd}, $j=s,d$ are defined as
\begin{subequations}
\label{eq:bhterms}
  \begin{align}
    \label{eq:bh_j}
    b_h^{j}(\boldsymbol{p}^{j},v )
      &:= -\sum_{K\in\mathcal{T}^j} \int_K p \nabla\cdot v \dif x
      + \sum_{K\in\mathcal{T}^j} \int_{\partial K} \bar{p}^j v\cdot n^j \dif s,
    \\
    \label{eq:bh_Ij}
    b_h^{I,j}(\bar{p}^j, \bar{v} )
      &:= -\int_{\Gamma^I}\bar{p}^j\bar{v}\cdot n^j \dif s.
  \end{align}
\end{subequations}
Before defining the  terms related to the transport equation, we point out that \cref{eq:Bhsd,eq:def_ah,eq:bhterms} are the same as in \cite{Cesmelioglu:2020} when the viscosity and $\kappa$ are both constants.

The form $B_h^{tr}(u; \boldsymbol{c}(t), \boldsymbol{w})$ in \cref{eq:hdgtr} discretizes the advective and diffusive parts of the transport equation:
\begin{equation}
  \label{eq:bilinearform_Bh}
  B_h^{tr}(u; \boldsymbol{c}, \boldsymbol{w}) =
  B_h^a(u; \boldsymbol{c}, \boldsymbol{w})
  + B_h^d(u; \boldsymbol{c}, \boldsymbol{w}).
\end{equation}
The advective part is defined as
\begin{multline}
  \label{eq:bilinearform_Bha}
  B_h^a(u; \boldsymbol{c}, \boldsymbol{w})
  := - \sum_{K\in \mathcal{T}} \int_K c \, u \cdot \nabla w \dif x
  + \sum_{K\in \mathcal{T}} \int_{\partial K} c\, u \cdot n \,(w-\bar{w}) \dif s
  - \sum_{K\in \mathcal{T}} \int_{\partial K^{\rm in}}  u\cdot n\, (c-\bar{c})\,(w-\bar{w}) \dif s,
\end{multline}
where $\partial K^{\rm in}$ denotes the inflow portion of the boundary on which
$u_h \cdot n < 0$, and the diffusive part is defined as
\begin{multline}
  \label{eq:bilinearform_Bhd}
    B_h^d(u; \boldsymbol{c}, \boldsymbol{w}) := 
    \sum_{K\in \mathcal{T}}\int_K \widetilde{D}(u)\nabla c \cdot \nabla w\dif x + \sum_{K\in \mathcal{T}} \tfrac{\beta_{tr} }{h_K} \int_{\partial K} [\widetilde{D}(u) n](c-\bar{c})\cdot(w-\bar{w}) n\dif s\\
    - \sum_{K\in \mathcal{T}}\int_{\partial K} [\widetilde{D}(u)\nabla c] \cdot n\,(w-\bar{w})\dif s- \sum_{K\in \mathcal{T}} \int_{\partial K}  ([\widetilde{D}(u)\nabla w] \cdot n) \,(c-\bar{c})\dif s,
\end{multline}
where $\beta_{tr} > 0$ is a penalty parameter. Here we pause again to mention that \cref{eq:bilinearform_Bh,eq:bilinearform_Bha,eq:bilinearform_Bhd} are the same as in \cite{Cesmelioglu:2021}, and a standard extension of the discretization analyzed in \cite{Wells:2011}.

To complete the discretization, we project the initial conditions $u_0$ and $c_0$ \cref{eq:initial} into $\boldsymbol{V}_h$ and $\boldsymbol{C}_h$, respectively.

\subsection{Properties of the numerical scheme}
The semi-discrete HDG scheme presented in \Cref{sec:semi-discrete} has various attractive features.
Besides local momentum conservation, a property of all HDG methods, this particular HDG method also conserves mass strongly, according to the definition defined in \cite{Kanschat:2010}. 

To be specific, the discrete velocity enjoys the following properties:
\begin{subequations}
\label{eq:properties_discvel}
\begin{align}
  \label{eq:masscons-1}
  -\nabla \cdot u_h &= \chi^d \Pi_Q (g_p-g_i) &&\forall x\in K, \ \forall K \in \mathcal{T},\\
  \label{eq:masscons-2}
  \jump{u_h\cdot n}&=0 &&\forall x\in F,\ \forall F\in \mathcal{F} \backslash \mathcal{F}^I,\\
  \label{eq:masscons-3}
  u_h^j\cdot n&=\bar{u}_h\cdot n &&\forall x\in F,\ \forall F\in \mathcal{F}^I,\, j=s,d,
\end{align}
\end{subequations}
where $\jump{\cdot}$ is the usual jump operator and $n$ is the unit
normal vector on~$F$. Note that \cref{eq:masscons-2,eq:masscons-3} imply that $u_h$ is $H(\text{div})$-conforming on the whole domain.
More details on \cref{eq:properties_discvel} can be found in \cite[Section 3.3]{Cesmelioglu:2020}.
Additionally, the scheme is consistent, that is, the solution to \cref{eq:system,eq:interface,eq:initial} satisfies \cref{eq:semidiscreteHDG}, as we discuss next.
\begin{lmm}[Consistency]
  \label{lem:consistency}
  Suppose that the solution $(u,p,c)$ to the Stokes/Darcy--transport system
  \cref{eq:system,eq:interface,eq:initial} satisfies $(u, p) \in L^2(0,T;X)$,  $\partial_t u\in L^2(0;T;L^2(\Omega^s))$, $c\in L^2(0,T;C)$, and $\partial_tc\in L^2(0,T;L^2(\Omega))$. Then $(\boldsymbol{u}(t),\boldsymbol{p}(t),\boldsymbol{c}(t))$, where  $\boldsymbol{u}:=(u, \gamma(u))$, $\boldsymbol{p}:= (p, \gamma(p^s), \gamma(p^d))$, $\boldsymbol{c}:=(c,\gamma(c))$,  satisfy the semi-discrete HDG scheme \cref{eq:semidiscreteHDG} for all $t>0$.
\end{lmm}
\begin{proof}
  The proof is as in \cite[Lemma 1]{Cesmelioglu:2020} and \cite[Lemma 6]{Cesmelioglu:2021} with minor modifications. We do not repeat the proof here, but mention that it is based on integration by parts in $a_h^s(c;\boldsymbol{u},\boldsymbol{v})$,  $\sum_{j=s,d} b_h^j(\boldsymbol{p}^{j}, v)$, $B_h^a(u,\boldsymbol{c},\boldsymbol{w})$, and $B_h^d(u,\boldsymbol{c},\boldsymbol{w})$, using  $\gamma(u)=u$ on $\Gamma_0^s$, $\gamma(p^j)=p^j$ on $\Gamma_0^j$, $j=s,d$, $\gamma(c)=c$ on $\Gamma_0$, the smoothness of the solution, the continuity of $\mu$ and $D$, and \cref{eq:system,eq:interface}.
\end{proof}
%
\section{Continuity, coercivity, and an inf-sup condition}
\label{sec:properties}
Let us recollect various known inequalities. Throughout this article we denote by $C>0$ a generic constant that is independent of the mesh size and the time step. From \cite[Lemma~1.46, Remark~1.47]{DiPietro:book}, for any $K\in\mathcal{T}$, we have
\begin{equation}
  \label{ineq:trace_disc}
  \norm{v}_{\partial K} \le C h_K^{-1/2} \norm{v}_K \quad \forall v\in P_k(K).
\end{equation}
We will also use the following versions of the continuous trace
inequality \cite[Theorem 1.6.6,(10.3.8)]{Brenner:book}:
\begin{align}
  \label{ineq:trace_cont}
  \norm{v}_{\partial K}^2 &\le C \del[1]{h_K^{-1} \norm{v}_K^2
  + h_K \norm{v}_{1,K}^2} && \forall v\in H^1(K),\\
  \label{ineq:trace_cont_inf}
  \norm{v}_{0,\infty,\partial K} &\le C && \forall v\in W^{1,\infty}(K),
\end{align}
where $C$ in \cref{ineq:trace_cont_inf} depends on $\|v\|_{1,\infty,K}.$
Regarding the trace on the interface, we have \cite[(1.24)]{Girault:2009}, \cite[Theorem 1.6.6]{Brenner:book}:
\begin{align}
  \label{eq:continuoustrace}
  \norm{v}_{\Gamma^I} &\le C \norm[0]{\nabla v}_{\Omega^s} && \forall
  v \in \cbr[1]{ v \in H^1(\Omega^s)\ :\ v=0\ \text{on}\ \Gamma^s},\\
  \label{eq:continuoustrace_s}
  \norm{v}_{\Gamma^I} &\le C\norm[0]{v}_{1,\Omega^s} && \forall
  v \in H^1(\Omega^s).
\end{align}
Furthermore, by 
\cite[Theorem 4.4]{Girault:2009}, for any $\boldsymbol{v}_h\in \boldsymbol{V}_h$, for $k_f\geq 1$, 
\begin{equation}\label{ineq:trace_interface_tnorm_vs}
\|v_h^s\|_{\Gamma^I}\leq C \tnorm{\boldsymbol{v}_h}_{v,s} \leq C \tnorm{\boldsymbol{v}_h}_{v}.    
\end{equation}
Similarly, 
for any $w_h\in C_h$, for $k_c\geq 1$, 
\begin{equation}\label{ineq:trace_inteface_tnorm_c}
\|w_h^s\|_{\Gamma^I}\leq C  \tnorm{\boldsymbol{w}_h}_{c}.  
\end{equation}
Next, we recall some inverse inequalities from \cite[Lemma 1.44, Lemma 1.50]{DiPietro:book}:
\begin{align}
      \label{ineq:inverse}
      \|\nabla v\|_K& \leq Ch_K^{-1}\|v\|_K  && \forall v\in P_k(K),\\
      \label{ineq:inverse-0-infty}
      \|v\|_{0,\infty,K}& \leq Ch_K^{-{\dim}/2}\|v\|_K && \forall v\in P_k(K).
\end{align}
The following Poincar\'{e}-type inequality follows from \cite[Proposition 4.5]{Girault:2009}, \cite[Remark 1.1]{Brenner:2003}:
\begin{equation}\label{ineq:bound_L2_by_tnorm_vs}
\|v\|_{\Omega^s}\le C\tnorm{\boldsymbol{v}}_{v,s} \quad\forall \boldsymbol{v}:=(v,\mu)\in H^1(\mathcal{T}_h^s)\times \bar{V}_h.
\end{equation}
The following version of Korn's first inequality is a consequence of \cite[(1.19)]{Brenner:2004}, \cite[Proposition 4.7]{Girault:2009}, \cite[p.110]{Riviere:book}:
\begin{equation}\label{ineq:Korn}
    \tnorm{\boldsymbol{v}_h}_{v,s}\leq C\sum_{K\in \mathcal{T}^s}(\|\varepsilon(v_h)\|_K^2+h_K^{-1}\|v_h-\bar{v}_h\|_{\partial K}^2).
\end{equation}
Continuity and coercivity $a_h(\cdot, \cdot)$, follow from \cite[Lemma 2, Lemma 3]{Cesmelioglu:2020} keeping in mind that $\mu$ satisfies \cref{eq:assumption_mu} and $\mathbb{K}$ satisfies \cref{eq:assumption_K}. They can be stated as follows:
\begin{lmm}[Continuity and coercivity of $a_h$]
\label{lem:a_h_coercivity}
  There exists a constant $C > 0$, independent of $h$, such that for all $\boldsymbol{u}, \boldsymbol{v}\in \boldsymbol{V}(h)$ and $\boldsymbol{c}\in \boldsymbol{C}(h)$,
  \begin{align}
    a_h(\boldsymbol{c};\boldsymbol{u},\boldsymbol{v})&\leq C\tnorm{\boldsymbol{u}}_{v'}\tnorm{\boldsymbol{v}}_{v'}.\label{thm:SDT_a_h_continuityvp}
  \end{align}
In addition, there exists a constant $C_a > 0$, independent of $h$ but dependent on $\kappa^*, \mu_*$, and $\mu^*$, and a constant $\beta^s_0 > 0$ such that if $\beta_s>\beta^0_s$, then
  \begin{equation}
  \label{eq:a_h_coercivity}
    a_h(c_h;\boldsymbol{v}_h,\boldsymbol{v}_h)\geq C_a\tnorm{\boldsymbol{v}_h}_v^2 \quad \forall \boldsymbol{v}_h \in \boldsymbol{V}_h, \quad\forall c_h\in C_h.
  \end{equation}
\end{lmm}
The inf-sup condition on the discrete spaces $\boldsymbol{V}_h$ and $\boldsymbol{Q}_h$ was proved in \cite{Cesmelioglu:2020} in the case of a continuous discrete velocity trace space $\bar{V}_h \cap \sbr[0]{C^0(\Gamma_0^s)}^{\dim}$. It is straightforward to show that the inf-sup condition also holds when the discrete velocity trace space is the larger discontinuous $\bar{V}_h$ space. 
\begin{thrm}
  \label{thm:infsup}
  There exists a constant $c_{{\rm inf}}^{\star} > 0$, independent of $h$, such
  that for any $\boldsymbol{q}_h \in \boldsymbol{Q}_h$,
  \begin{equation}
    c_{{\rm inf}}^{\star} \tnorm{\boldsymbol{q}_h}_{p} \le
    \sup_{\substack{\boldsymbol{v}_h \in \boldsymbol{V}_h \\ \boldsymbol{v}_h \ne \boldsymbol{0}}}
    \frac{\sum_{j=s,d}\del[1]{b_h^{j}(\boldsymbol{q}_h^j, v_h) +
    b_h^{I,j}(\bar{q}_h^j,\bar{v}_h)}}{\tnorm{\boldsymbol{v}_h}_{v}}.
  \end{equation}
\end{thrm}
Note that the proof of this theorem as well as the error analysis requires appropriate interpolation operators onto $V_h$ and $\bar{V}_h$. For $V_h$ we consider the BDM interpolation operator $\Pi_{V} : \sbr[0]{H^1(\Omega)}^{\dim} \rightarrow V_h$ which is such that if $u \in \sbr[0]{H^{k_f+1}(K)}^{\dim}$, $K\in \mathcal{T}$, then (see, for example, \cite[Lemma 7]{Hansbo:2002} and \cite[Section III.3]{Brezzi:book}):
\begin{subequations}
\label{lem:BDM}
\begin{align}
    \label{lem:BDM-i}
    \int_K q(\nabla \cdot u-\nabla \cdot \Pi_V u) \dif x &=0, && \forall q\in P_{k_f-1}(K),
    \\
    \label{lem:BDM-ii}
    \int_F \bar{q} (u-\Pi_Vu) \cdot n \dif s &=0, && \forall \bar{q}\in P_{k_f}(F),\ F \subset \partial K,
    \\
    \label{lem:BDM-iii}
    \Pi_Vu & \in H({\rm div};\Omega),
    \\
    \label{lem:BDM-iv}
    \norm{u-\Pi_V u}_{m,K} &\le C h_K^{\ell-m} \norm{u}_{\ell, K}, && m=0,1,2\ m \le \ell \le k_f+1,
\end{align}
\end{subequations}
where we remark that $F$ is an edge if $\dim=2$ and a face if $\dim = 3$.
Furthermore, for $u\in W^{1,\infty}(K)$, \cite[(2.33)]{Guzman:2016},
%
\begin{equation}\label{ineq:inf-PiV}
    \|u-\Pi_V u\|_{\infty,K} + h_K \|\nabla (u-\Pi_V u)\|_{\infty,K} \leq C h_K\|u\|_{1,\infty,K}.    
\end{equation}
The interpolant onto the trace space $\bar{V}_h$ is defined by $\bar{\Pi}_V:\sbr[0]{H^1(\Omega^s)}^{\dim} \to \bar{V}_h$ such that
\begin{equation*}
  \bar{\Pi}_V u=
  \begin{cases}
    (P_{\bar{V}}u)|_{F} & \text{ if } F\in \mathcal{F}^s\backslash\mathcal{F}^I,
    \\
    (\Pi_V u)^s|_F  & \text{ if } F\in \mathcal{F}^I,
  \end{cases}
\end{equation*}
where $P_{\bar{V}}$ is the $L^2$-projection onto $\bar{V}_h$. It is straightforward to deduce the following estimates using \cref{ineq:trace_cont} and the fact that $h= \max_{K\in \mathcal{T}}h_K$:
For $v \in \sbr[0]{H^{\ell}(\Omega^s)}^{\dim}$, $1 \leq \ell \leq k_f+1$,
 \begin{align}
    \label{eq:l2projectionVbar}
    \norm[0]{v-\bar{\Pi}_V v}_{\partial K} &\le C h_K^{\ell-1/2} \norm[0]{v}_{\ell,K},\\
        \label{eq:l2projectionVbar-2}
    \norm[0]{\Pi_V v-\bar{\Pi}_V v}_{\partial K} &\le C h_K^{\ell-1/2} \norm[0]{v}_{\ell,K}.
\end{align}
We finish this section by noting that by \cref{lem:BDM-i,lem:BDM-ii,lem:BDM-iii} the solution $u$ of \cref{eq:system,eq:interface} under the assumption that $u\in [H^{k+1}(K)]^{\dim}$, for all $K\in \mathcal{T}$, satisfies
\begin{equation}
    \label{eq:bhPiV}
    \sum_{j=s,d} \big( b_h^j(\boldsymbol{q}_h^{j}, u-\Pi_Vu) + b_h^{I,j}(q_h^j, \gamma(u^s)-\bar{\Pi}_Vu) \big)
    = 0 \quad \forall \boldsymbol{q}_h\in \boldsymbol{Q}_h.
\end{equation}
%
\section{Fully discrete numerical scheme}
\label{sec:sequential_strategy}
Let us now describe the fully discrete HDG method and decoupling strategy used to solve the Stokes/Darcy and transport problems sequentially. 
%
For the time discretization, we partition the time interval $J$ as: $0=t^0 < t^1 < \hdots <
t^N=T$. For simplicity, we assume a uniform partition with $t^{n+1} - t^n = \Delta t$ for $0 \le n \le N-1$. We denote a function $h(t)$ at time level $t^n$ by $h^n := h(t^n)$ and for a sequence $\cbr{u^n}_{n\geq 1}$ we denote by
$d_tu^{n} = (u^{n} - u^{n-1})/\Delta t$ a first order difference operator. 

In the first step of our sequential algorithm, given an initial velocity $u_h^0$ in the Stokes domain and an initial concentration $\boldsymbol{c}_h^0$, we solve the Stokes/Darcy problem and obtain a velocity in the entire region. This velocity, with properties given by \cref{eq:properties_discvel}, is then substituted into the concentration problem. This approach is repeated for all time steps with the initial velocity and concentration being replaced by the last computed velocity and concentration solutions. We summarize the fully discrete problem in \Cref{alg:uhfirst}.
\begin{algorithm}
  \caption{Sequential algorithm}
  \label{alg:uhfirst}
  \begin{algorithmic}
    \State Set $u_h^0 = \Pi_{V} u_0$, $\boldsymbol{c}_h^0 = (\Pi_C c_0, \bar{\Pi}_C c_0)$. \;
    \For {$n=1, \hdots, N$}
      \State 1.
        Find $(\boldsymbol{u}_h^n, \boldsymbol{p}_h^n)\in \boldsymbol{X}_h$ such that for all $(\boldsymbol{v}_h, \boldsymbol{q}_h) \in \boldsymbol{X}_h$
        \begin{multline}\label{eq:hdgsd_sequence}
         \sum_{K\in \mathcal{T}^s}\int_{\partial K}d_t u_h^n \cdot v_h \dif x +  B_h^{sd}(\boldsymbol{c}_h^{n-1}; (\boldsymbol{u}_h^n, \boldsymbol{p}_h^n), (\boldsymbol{v}_h, \boldsymbol{q}_h) )= \sum_{K\in \mathcal{T}^s}\int_K f^s(c_h^{n-1})\cdot v_h \dif x
          \\
          + \sum_{K\in \mathcal{T}^d}\int_K \mathbb{K}^{-1}(c_h^{n-1})f^d(c_h^{n-1})\cdot v_h\dif x + \sum_{K\in \mathcal{T}^d}\int_K (g_p^n-g_i^n)\,q_h \dif x.
        \end{multline}
      \State 2.
        Find $\boldsymbol{c}_h^n\in \boldsymbol{C}_h$ such that for all $\boldsymbol{w}_h\in \boldsymbol{C}_h$
        \begin{equation}\label{eq:hdgtr_sequence}
          \sum_{K\in \mathcal{T}}\int_K \phi\, d_t c_h^n w_h \dif x
          + B_h^{tr}(u_h^{n}; \boldsymbol{c}_h^{n}, \boldsymbol{w}_h)+ \sum_{K\in \mathcal{T}^d}\int_Kg_p^nc_h^nw_h\dif x= \sum_{K\in \mathcal{T}^d} \int_K c_{I} g_{i}^n  w_h \dif x. 
        \end{equation}
    \EndFor
  \end{algorithmic}
\end{algorithm}
\begin{rmrk}
    In \Cref{alg:uhfirst}, $\Pi_C$ denotes the $L^2$-projection onto $C_h$ and 
    $\Pi_V u_0$ is understood as $\Pi_V$ applied to the extension of $u_0$ to $\Omega$ by zero assuming $u_0\in H^1_0(\Omega^s)$. 
    We note that this choice of $u_h^0$ satisfies normal continuity across the interfaces in $\Gamma_0^s$ and has zero divergence in $\Omega^s$ under the additional assumption that $\nabla \cdot u_0=0$ in $\Omega^s$. We further remark that the properties in \cref{eq:properties_discvel} hold for $\boldsymbol{u}_h^n$ for each time step $n=0,\hdots, N$.
\end{rmrk}
We conclude this section by stating some preliminary results obtained by Taylor's theorem \cite[Lemma 3.2]{Cesmelioglu:2020a}.
For a function $z$ defined on $D\times[0,T]$, assuming enough regularity, we have the following results:
  \begin{subequations}
  \label{ineq:taylor}
  \begin{align} 
  \label{ineq:taylor1}
    \sum\limits_{m=1}^n\|\partial_tz^m-d_tz^m\|_{D}^2 &\leq C\Delta t\|\partial_{tt}z\|^2_{L^2(0,T;L^2(D))},
    \\
    \label{ineq:taylor3}
    \Delta t\sum\limits_{m=1}^n\|d_tz^m\|_{\ell,D}^2&\leq\|\partial_t z\|^2_{L^2(0,T;H^{\ell}(D))}, \quad \ell=0,1,
  \end{align}
  \end{subequations}
 where $\|f\|_{L^2(a,b;X)}:=\big(\int_{a}^b\|f(t)\|_X^2\dif t\big)^{1/2}$. Note that the inequalities in \cref{ineq:taylor} for $\ell=0$ were presented in \cite[Lemma 3.2]{Cesmelioglu:2020a} and that it is straightforward to extend \cref{ineq:taylor3} to $\ell=1$.
%
\section{Main results}
\label{sec:main_results}
In this section we present our main results. For the error estimates we will make use of the following definition of the discrete in time norm:
\begin{equation*}
    \|f\|_{\ell^2(0,T;X)}=\big(\Delta t\sum_{m=1}^n\|f^m\|_X^2\big)^{1/2}.     
\end{equation*}
Before proving a priori error estimates for the discrete velocity, pressure, and concentration, we first state the well-posedness of the discrete Stokes/Darcy problem \cref{eq:hdgsd_sequence}. Well-posedness of the discrete transport problem \cref{eq:hdgtr_sequence} is proven in \Cref{ss:existuniq_ch} as it depends on results obtained in \Cref{sec:error-estimates}.
\begin{thrm}\label{thm:existence-sd}
  Let $\beta_s>\beta_s^0$ be as in Lemma~\ref{lem:a_h_coercivity} and $n\geq 1$. Then given $u_h^{n-1}\in V_h$ and $\boldsymbol{c}_h^{n-1}\in \boldsymbol{C}_h$, there exists a unique solution $(\boldsymbol{u}_h^n,
  \boldsymbol{p}_h^n)\in \boldsymbol{X}_h$ to \cref{eq:hdgsd_sequence} that satisfies
  \begin{equation} \label{ineq:uhnphnbd}
      \tnorm{\boldsymbol{u}_h^n}_{v}+\tnorm{\boldsymbol{p}_h^n}_{p}
      \leq C\Big(\dfrac1{\Delta t}\|u_h^{n-1}\|_{\Omega^s}+\|f^s(c_h^{n-1})\|_{\Omega^s}+\dfrac1{K_*}\|f^d(c_h^{n-1})\|_{\Omega^d}+\|g_p^n-g_i^n\|_{\Omega^d}\Big).
  \end{equation}
  %
\end{thrm}
\begin{proof}
    The result follows by applying the abstract theory for saddle point problems \cite[Theorem 2.34]{Ern:book} together with Theorem~\ref{thm:infsup} and \cref{eq:a_h_coercivity}.
\end{proof}
\subsection{Error estimates for the discrete velocity}
\label{sec:error-estimates}
In this section we derive estimates for the error $u_h^n-u^n$, for each $n\geq 0$, given error estimates for the discrete concentration in previous time steps.
To do so, we define the following:
\begin{align*}
  \xi_u^n & : = u^n - \Pi_V u^n,&\, \zeta_u^n & : = u_h^n-\Pi_V u^n,  & \xi_p^n & : = p^n-\Pi_Q p^n, &  \, \zeta_p^n & : = p_h^n-\Pi_Q p^n,\\
  \bar{\xi}_u^n & : = \gamma(u^{sn})-\bar{\Pi}_{V} u^{n},& \,\bar{\zeta}_u^n & : =\bar{u}_h^n - \bar{\Pi}_V u^{n}, &\bar{\xi}_p^{jn} & : = \gamma(p^{jn})-\bar{\Pi}_Q^j p^{jn}, &  \, \bar{\zeta}_p^{jn} & : = \bar{p}_h^{j}-\bar{\Pi}_Q^j p^{jn},
\end{align*}
where $\Pi_Q$ is the $L^2$-projection onto $Q_h$ and $\bar{\Pi}_Q^j$ is the $L^2$-projection onto $\bar{Q}_h^j$, $j=s,d$. For the case $n=0$, $\Pi_V u^0$ is understood as $\Pi_V$ applied to the extension of $u_0$ to $\Omega$ by zero assuming $u_0\in H^1_0(\Omega^s)$. Therefore, $\zeta_u^0=0$.
Furthermore, note that the following identities hold:
\begin{align}\label{eq:split_error}
  u^n-u_h^n&=\xi_u^{n}-\zeta_u^n, & \gamma(u^n)-\bar{u}_h^n&=\bar{\xi}_u^n-\bar{\zeta}_u^n,  \\
  p^n-p_h^n&=\xi_p^n-\zeta_p^n, & \gamma(p^{jn})-\bar{p}_h^{jn}&=\bar{\xi}_p^{jn}-\bar{\zeta}_p^{jn}, \quad j=s, d.
\end{align}
To be consistent with the notation used in previous sections, we set
$\boldsymbol{\ell}_u^n=(\ell_u^n,\bar{\ell}_u^n)$, $\boldsymbol{\ell}_p^n :=
(\ell_p^n,\bar{\ell}_p^{sn},\bar{\ell}_p^{dn})$, and $\boldsymbol{\ell}_p^{jn}:=
(\ell_p^n,\bar{\ell}_p^{jn})$, for $\ell=\xi, \zeta$ and $j=s, d$.

Here we recall the following results on the interpolation errors \cite[Lemma 7, Lemma 8]{Cesmelioglu:2020}. Suppose that $u$ is such that $u^{s}\in [H^{\ell}(\Omega^s)]^{\dim}$ and $u^{d}\in [H^{\ell-1}(\Omega^d)]^{\dim}$ for $2\le \ell\leq k_f+1$, and that $p^{j}\in H^{r}(\Omega^j)$ for $0\leq r \leq k_f$ and $j=s, d$. Then
\begin{subequations}
\begin{align}
    \label{eq:xi_u-1}
    \tnorm{\boldsymbol{\xi}_u}_{v',s} &\le Ch^{\ell-1}\|u\|_{H^{\ell}(\Omega^s)},
    \\
    \label{eq:xi_u-2}
    \tnorm{\boldsymbol{\xi}_u}_{v'}  &\le Ch^{\ell-1}(\|u\|_{H^{\ell}(\Omega^s)} + \|u\|_{H^{\ell-1}(\Omega^d)}),
    \\
    \label{eq:xi_p}
    \tnorm{\boldsymbol{\xi}_p^{j}}_{p,j} &\leq Ch^r \|p\|_{H^{r}(\Omega^j)}.
\end{align}
\end{subequations}

\begin{lmm} \label{lem:vanishingbh}
Let $(u,p)$ be the velocity solution of \cref{eq:system,eq:interface,eq:initial}, $\bar{u}=\gamma(u^s)$.
Then for any $n\geq 1$,
  \begin{align}
  \label{eq:vanishingbh2}
    \sum_{j=s,d}(b_h^j(\boldsymbol{\xi}_p^{jn}, v_h) + b_h^{I,j}(\bar{\xi}_p^{jn}, \bar{v}_h)\big)&=0, \quad \forall \boldsymbol{v}_h\in \boldsymbol{V}_h,\\
  \label{eq:vanishingbh3}
    \sum_{j=s,d} (b_h^j(\boldsymbol{q}_h^j, \xi_u^n) + b_h^{I,j}(\bar{q}_h^j, \bar{\xi}_u^n)\big)&=0, \quad \forall \boldsymbol{q}_h\in \boldsymbol{Q}_h.
  \end{align}
\end{lmm}
\begin{proof}
  The proof is based on the properties of the numerical scheme  \cref{eq:properties_discvel}, the properties of the BDM projection $\Pi_V$ in \cref{lem:BDM}, and the properties of the $L^2$-projections $\Pi_Q$ and $\bar{\Pi}_Q^j$, $j=s,d$.
Indeed, \cref{eq:vanishingbh2} follows after noting that $\nabla\cdot v_h\in P_{k_f-1}(K)$, $v_h\cdot n^j, \bar{v}_h\cdot n^j\in P_{k_f}(F)$, and using the definitions of the $L^2$-projections $\Pi_Q$ and $\bar{\Pi}_Q^j$,  $j = s, d$,
 \begin{multline*}
\sum_{j=s,d}\big(b_h^j(\boldsymbol{\xi}_p^{jn}, v_h) + b_h^{I,j}(\bar{\xi}_p^{jn}, \bar{v}_h)\big) =\sum_{j=s,d}\Big( -\sum_{K\in\mathcal{T}^j} \int_K (p^n-\Pi_Q p^n) \nabla\cdot v_h \dif x\\
    + \sum_{K\in\mathcal{T}^j} \int_{\partial K} (\gamma(p^{jn})-\bar{\Pi}_Q^j p^{jn})v_h \cdot n^j \dif s- \int_{\Gamma^I}(\gamma(p^{jn})-\bar{\Pi}_Q^j p^{jn})\bar{v}_h\cdot n^j \dif s\Big) =0,
\end{multline*}
while \cref{eq:vanishingbh3} is exactly the same as \cref{eq:bhPiV}, evaluated at $t=t_n$, and rewritten by using the definitions of $\xi_u^n$ and $\bar{\xi}_u^n$.
\end{proof}
\begin{thrm}[Error equation for \cref{eq:hdgsd_sequence}] \label{thm:error_equation}
  There holds
  \begin{multline}\label{eq:error_equation_0}
    \sum\limits_{K\in\mathcal{T}^s}\int_K (d_tu_h^n-\partial_t u^n)\cdot v_h\,{\rm d}x  +a_h(\boldsymbol{c}_h^{n-1}; \boldsymbol{\zeta}^n_u, \boldsymbol{v}_h)\\
    +\sum_{j=s,d} \big(b_h^j(\boldsymbol{\zeta}_p^{jn}, v_h)+ b_h^{I,j}(\bar{\zeta}_p^{jn}, \bar{v}_h)\big)
    +\sum_{j=s,d} \big(b_h^j(\boldsymbol{q}_h^{j}, \zeta_u^n)+ b_h^{I,j}(\bar{q}_h^{j}, \bar{\zeta}_u^n)\big)
    \\
    =a_h(\boldsymbol{c}_h^{n-1}; \boldsymbol{\xi}^n_u, \boldsymbol{v}_h)+ a_h(\boldsymbol{c}^n; \boldsymbol{u}^n,\boldsymbol{v}_h)-a_h(\boldsymbol{c}_h^{n-1}; \boldsymbol{u}^n, \boldsymbol{v}_h)\\
    + \sum\limits_{K\in\mathcal{T}^s}\int_K [f^s(c_h^{n-1})-f^s(c^{n})]\cdot v_h \dif x  
  +  \sum\limits_{K\in\mathcal{T}^d}\int_K  [\mathbb{K}^{-1}(c_h^{n-1})f^d(c_h^{n-1})-\mathbb{K}^{-1}(c^{n})f^d(c^{n})]\cdot v_h\dif x.
  \end{multline}
\end{thrm}
\begin{proof}
  By Lemma~\ref{lem:consistency} at time  $t=t^n$, we have that for all $(\boldsymbol{v}_h, \boldsymbol{q}_h) \in \boldsymbol{X}_h$,
  \begin{multline*}
    \sum_{K\in \mathcal{T}^s}\int_{ K}\partial_t u^n \cdot v_h\, \mathrm{d}x + B_h^{sd}(\boldsymbol{c}^n; (\boldsymbol{u}^n, \boldsymbol{p}^n), (\boldsymbol{v}_h, \boldsymbol{q}_h) )
    =    \sum\limits_{K\in\mathcal{T}^s}\int_K  f^s(c^n) \cdot v_h \,{\rm d} x\\
    +  \sum\limits_{K\in\mathcal{T}^d}\int_K  \mathbb{K}^{-1}(c^n)f^d(c^n)\cdot v_h\,{\rm d} x +\sum_{K\in \mathcal{T}^d}\int_K(g_p^n-g_i^n)q_h\,{\rm d}x. \nonumber
  \end{multline*}
  Subtracting this equation from \cref{eq:hdgsd_sequence}, we obtain
  \begin{multline}\label{eq:error_equation_1}
    \sum\limits_{K\in\mathcal{T}^s}\int_K (d_tu_h^n-\partial_t\,u^n)\cdot v_h\,{\rm d}x
    + B_h^{sd}(\boldsymbol{c}^{n-1}_h; (\boldsymbol{u}^n_h, \boldsymbol{p}^n_h), (\boldsymbol{v}_h, \boldsymbol{q}_h) ) -  B_h^{sd}(\boldsymbol{c}^n; (\boldsymbol{u}^n, \boldsymbol{p}^n), (\boldsymbol{v}_h, \boldsymbol{q}_h) )\\
   = \sum\limits_{K\in\mathcal{T}^s}\int_K  [f^s(c_h^{n-1})-f^s(c^{n})]\cdot v_h \dif x 
    +  \sum\limits_{K\in\mathcal{T}^d}\int_K [\mathbb{K}^{-1}(c_h^{n-1})f^d(c_h^{n-1})-\mathbb{K}^{-1}(c^{n})f^d(c^{n})]\cdot v_h\dif x,
  \end{multline}
  for all $(\boldsymbol{v}_h, \boldsymbol{q}_h) \in \boldsymbol{X}_h$.
  Then, by \cref{eq:split_error}, the $B_h^{sd}$ terms can be rewritten as
  \begin{equation}
  \label{eq:B_sd_unh_B_sd_un}
      \begin{split}
    B_h^{sd}(\boldsymbol{c}^{n-1}_h; (\boldsymbol{u}^n_h, \boldsymbol{p}^n_h), (\boldsymbol{v}_h, \boldsymbol{q}_h) ) &- B_h^{sd}(\boldsymbol{c}^n; (\boldsymbol{u}^n, \boldsymbol{p}^n), (\boldsymbol{v}_h, \boldsymbol{q}_h) ) 
    \\
    &=B_h^{sd}(\boldsymbol{c}_h^{n-1}; (\boldsymbol{u}^n_h-\boldsymbol{u}^n, \boldsymbol{p}^n_h-\boldsymbol{p}^n), (\boldsymbol{v}_h, \boldsymbol{q}_h) ) 
    \\
    &\quad +B_h^{sd}(\boldsymbol{c}_h^{n-1}; (\boldsymbol{u}^n, \boldsymbol{p}^n), (\boldsymbol{v}_h, \boldsymbol{q}_h) )-B_h^{sd}(\boldsymbol{c}^n; (\boldsymbol{u}^n, \boldsymbol{p}^n), (\boldsymbol{v}_h, \boldsymbol{q}_h) ) 
    \\
    &=a_h(\boldsymbol{c}_h^{n-1}; \boldsymbol{\zeta}^n_u,\boldsymbol{v}_h) - a_h(\boldsymbol{c}_h^{n-1}; \boldsymbol{\xi}^n_u,\boldsymbol{v}_h) +a_h(\boldsymbol{c}_h^{n-1}; \boldsymbol{u}^n, \boldsymbol{v}_h )- a_h(\boldsymbol{c}^n; \boldsymbol{u}^n,\boldsymbol{v}_h) 
    \\
    &\quad+\sum_{j=s,d} \big(b_h^j(\boldsymbol{\zeta}_p^{jn}, v_h) + b_h^{I,j}(\bar{\zeta}_p^{jn}, \bar{v}_h)\big)+\sum_{j=s,d} \big(b_h^j(\boldsymbol{q}_h^j, \zeta_u^n) + b_h^{I,j}(\bar{q}_h^j, \bar{\zeta}_u^n)\big),          
      \end{split}
  \end{equation}
    where we applied \cref{eq:vanishingbh2,eq:vanishingbh3}. Combining this with \cref{eq:error_equation_1} yields the result.  
\end{proof}
The velocity error at each time step depends on the error in concentration from the previous time step. Therefore, for the velocity error estimates, we will need some auxiliary results related to the concentration error. To estimate the error of the concentration, we use the continuous
interpolant $\mathcal{I}c\in C_h\cap \mathcal{C}^0(\bar{\Omega})$ of $c$ \cite{Brenner:book},
and we set
$\bar{\mathcal{I}} c(t) = \mathcal{I} c|_{\Gamma^0}(t)\in \bar{C}_h$. Denoting the restriction of $c$ to $\Gamma^0$ by $\bar{c}$, we define 
\begin{align}
    \xi_c^n &= c^n - \mathcal{I} c^n, & \zeta_c^n &= c_h^n- \mathcal{I} c^n, &  \boldsymbol{\xi}_c^n &=(\xi_c^n, \bar{\xi}_c^n),\label{eq:xiczetac}\\
    \bar{\xi}_c^n &= \bar{c}^n - \bar{\mathcal{I}} c^n, &  \bar{\zeta}_c^n   &= \bar{c}_h^n - \bar{\mathcal{I}} c^n,
    & \boldsymbol{\zeta}_c^n &=(\zeta_c^n, \bar{\zeta}_c^n).\nonumber
\end{align}
Note that:
\begin{equation}\label{eq:err-split-c}
  c^n - c_h^n = \xi_c^n -  \zeta_c^n, \qquad
  \bar{c}^n - \bar{c}_h^n = \bar{\xi}_c^n - \bar{\zeta}_c^n.
\end{equation}
Furthermore, we have the following interpolation estimate \cite[Section 4.4]{Brenner:book} for $r=0, 1$ and $1 \le \ell \le k_c$:
\begin{align}
  \label{eq:est_interpolant_brenner}
  \|\xi_c\|_{r,K} &\leq C h_K^{\ell+1-r} \|c\|_{\ell+1,K},
\end{align}
\begin{thrm}
Let $c_0^s\in H^{k_c+1}(\Omega^s)$, $c^s \in L^2(0, T; H^{k_c+1}(\Omega^s))$, $c^d\in H^1(0,T;L^2(\Omega^d))$ such that $\partial_t c^s\in L^2(0,T, H^1(\Omega^s))$, and $\bar{c}=\gamma(c)$ on $\Gamma_0$. Then we have the following estimates:
\begin{subequations}
\label{eq:boundcmestimates}
\begin{equation}\label{eq:bound_cm_chm_minus_1_1}
\Delta t\sum\limits_{m=1}^n \sum\limits_{K\in \mathcal{T}^s} h_K^2 \|c^m-c^{m-1}_h\|_{1,K}^2
\leq C\Big(h^2 (\Delta t)^2 \|\partial_t c\|^2_{L^2(0,T;H^1(\Omega^s))}+ \Delta t\sum\limits_{m=1}^n\sum_{K\in \mathcal{T}^s}h_K^2 \|c^{m-1}-c_h^{m-1}\|_{1,K}^2\Big), 
\end{equation}
\begin{multline}\label{eq:bound_cm_chm_minus_chm_interface}
    \Delta t\sum\limits_{m=1}^n\|\bar{c}^m-\bar{c}^{m-1}_h\|_{\Gamma^I}^2 \leq C\Big(\Delta t\sum\limits_{m=1}^n\tnorm{\boldsymbol{\zeta}_c^{m-1}}_c^2\\
    +(\Delta t)^2\|\partial_t c\|_{L^2(0,T;H^1(\Omega^s))}^2+ h^{2k_c+1} \big(\Delta t \|c_0\|_{k_c+1,\Omega^s}^2 +\|c\|_{\ell^2(0,T;H^{k_c+1}(\Omega^s))}^2\big)\Big),
\end{multline}
and 
\begin{equation}\label{eq:bound_cm_chm_minus_1_omega_d}
\Delta t\sum\limits_{m=1}^n\Big(\sum\limits_{K\in \mathcal{T}} (\|c^m-c^{m-1}_h\|_K^2\Big)\leq C\Big( (\Delta t)^2\|\partial_t c\|_{L^2(0,T;L^2(\Omega))}^2 
+ \Delta t\sum\limits_{m=1}^n\big(\sum\limits_{K\in \mathcal{T}} \|c^{m-1}-c_h^{m-1}\|_{K}^2\big)\Big).
\end{equation}
\end{subequations}
\end{thrm}
\begin{proof}
We first prove \cref{eq:bound_cm_chm_minus_1_1}. 
By the triangle inequality and the definition of $d_t$,
\begin{align*}
\sum\limits_{K\in \mathcal{T}^s}h_K^2 \|c^m-c^{m-1}_h\|_{1,K}^2 
&\leq \sum\limits_{K\in \mathcal{T}^s}2h_K^2 \Big((\Delta t)^2\|d_t c^m\|_{1,K}^2+ \|c^{m-1}-c^{m-1}_h\|_{1,K}^2\Big).
\end{align*}
Multiplying this by $\Delta t$, summing from $m=1$ to $n$, and using \cref{ineq:taylor3}, we obtain \cref{eq:bound_cm_chm_minus_1_1}.
We now prove \cref{eq:bound_cm_chm_minus_chm_interface}. We have
\begin{align}\label{ineq:cbarm-cbarhm-1}
\|\bar{c}^m-\bar{c}^{m-1}_h\|_{\Gamma^I}^2&\leq 2 (\|\bar{c}^m-\bar{c}^{m-1}\|_{\Gamma^I}^2+\|\bar{c}^{m-1}-\bar{c}^{m-1}_h\|_{\Gamma^I}^2).
\end{align}
Since $\bar{c}=c^s|_{\Gamma_I}$ on $\Gamma^I$, and $(c^s)^m-(c^s)^{m-1} \in H^1(\Omega^s)$, by \cref{eq:continuoustrace_s}, the first term on the right side of \cref{ineq:cbarm-cbarhm-1} is bounded as follows:
\begin{align}\label{ineq:cbarm-cbarm-1}
\|\bar{c}^m-\bar{c}^{m-1}\|_{\Gamma^I}^2&\leq C\|c^m-c^{m-1}\|_{1,\Omega^s}^2=C(\Delta t)^2\|d_t c^m\|_{1,\Omega^s}^2.
\end{align}
Splitting the second term on the right side of \cref{ineq:cbarm-cbarhm-1} by using $\mathcal{I}c=\bar{\mathcal{I}}c$ for any  $F\in \mathcal{F}^I$ gives:
\begin{align}\label{ineq:cbarm-1-cbarm-1}
\|\bar{c}^{m-1}-\bar{c}^{m-1}_h\|_{\Gamma^I}^2&\leq 2(\|c^{m-1}-\mathcal{I} c^{m-1}\|_{\Gamma^I}^2+\|\bar{c}^{m-1}_h-\bar{\mathcal{I}}c^{m-1}\|_{\Gamma^I}^2 )= 2(\|\xi_c^{m-1}\|_{\Gamma^I}^2+\|\bar{\zeta}_c^{m-1}\|_{\Gamma^I}^2).
\end{align}
The first term on the right hand side of \cref{ineq:cbarm-1-cbarm-1} is bounded by \cref{ineq:trace_cont} and \cref{eq:est_interpolant_brenner} as follows:
\begin{equation}\label{ineq:xicm-1}
    \|\xi_c^{m-1}\|_{\Gamma^I}^2
    \leq \sum\limits_{K\in \mathcal{T}^s}C (h_K^{-1} \|\xi_c^{m-1}\|_K^2   + h_K \|\xi_c^{m-1}\|_{1,K}^2) 
     \leq Ch^{2k_c+1} \|c^{m-1}\|_{k_c+1,\Omega^s}^2.
\end{equation}
Using \cref{ineq:trace_inteface_tnorm_c} and the definition of $\tnorm{\cdot}_c$,
\begin{multline}\label{ineq:zetacm-1}
    \|\bar{\zeta}_c^{m-1}\|_{\Gamma^I}^2 \leq 2\big(\|\bar{\zeta}_c^{m-1}-(\zeta_c^s)^{m-1}\|_{\Gamma^I}^2 +\|(\zeta_c^s)^{m-1}\|_{\Gamma^I}^2\big)\\
   \leq 2\Big(h\sum\limits_{K\in \mathcal{T}^s}h_K^{-1}\|\bar{\zeta}_c^{m-1}-\zeta_c^{m-1}\|_{\partial K}^2 +\|(\zeta_c^s)^{m-1}\|_{\Gamma^I}^2\Big)\leq C\tnorm{\boldsymbol{\zeta}_c^{m-1}}_c^2.
\end{multline}
Collecting \crefrange{ineq:cbarm-cbarhm-1}{ineq:zetacm-1}, we obtain
\begin{equation*}
    \|\bar{c}^m-\bar{c}^{m-1}_h\|_{\Gamma^I}^2
    \leq C\Big((\Delta t)^2\|d_t c^m\|^2_{1,\Omega^s}
   +h^{2k_c+1} \|c^{m-1}\|_{k_c+1,\Omega^s}^2+\tnorm{\boldsymbol{\zeta}_c^{m-1}}_c^2\Big).
\end{equation*}
\Cref{eq:bound_cm_chm_minus_chm_interface} now follows after multiplying the above inequality by $\Delta t$, summing from $1$ to $n$, and using \cref{ineq:taylor3}.
We next prove \cref{eq:bound_cm_chm_minus_1_omega_d}. By the triangle inequality, 
\begin{align}\label{ineq:cm-chm-1}
\sum\limits_{K\in \mathcal{T}}\|c^m-c^{m-1}_h\|_K^2
    &\leq \sum\limits_{K\in \mathcal{T}}2\big((\Delta t)^2\|d_t c^m\|_{K}^2+\|c^{m-1}-c^{m-1}_h\|_K^2\big).
\end{align}
The results follows by multiplying \cref{ineq:cm-chm-1} by $\Delta t$, summing from $m=1$ to $n$, and using \cref{ineq:taylor3} as before.
\end{proof}
Now that the auxiliary result is established, we proceed with proving error estimates the velocity.
\begin{thrm} \label{thm:vel-error}
Let $c$ and $u$ be the solutions of \cref{eq:system,eq:interface,eq:initial} such that
\begin{align*}
&u^s \in L^2(0,T;[H^{k_f+1}(\Omega^s)]^{\dim})\cap L^{\infty}(0,T;[W^{1,\infty}(\Omega^s)]^{\dim}), \\
&\partial_{t}u^s \in L^2(0,T;[H^{k_f}(\Omega^s)]^{\dim}), \partial_{tt}u^s \in L^2(0,T;[L^2(\Omega^s)]^{\dim}),\\
&u^d\in L^2(0,T;[H^{k_f}(\Omega^d)]^{\dim})\cap L^{\infty}(0,T;[L^{\infty}(\Omega^d)]^{\dim}),\\
&p\in L^2(0,T;L^2(\Omega)), \nabla p^d\in L^{\infty}(0,T;[L^{\infty}(\Omega^d)]^{\dim}),\\
&c^s\in L^2(0,T;H^{k_c+1}(\Omega^s)), \partial_t c^s \in L^2(0,T;H^1(\Omega^s)), \\
&c^d\in L^2(0,T;L^2(\Omega^d)), \partial_t c^d \in L^{\infty}(0,T;L^2(\Omega^d)),\\
&u_0\in [H^1_0(\Omega^s)]^{\rm \dim}, \nabla \cdot u^0=0, c_0^s\in H^{k_c+1}(\Omega^s), 
\end{align*}
and let $k_f, k_c\geq {\dim}-1$, and $n\geq 1$. Suppose that $u_h^0, \hdots, u_h^{n-1}\in V_h$ and $\boldsymbol{c}_h^0, \hdots, \boldsymbol{c}_h^{n-1}\in \boldsymbol{C}_h$, the solutions of \cref{eq:hdgsd_sequence} and \cref{eq:hdgtr_sequence}, respectively, are known and satisfy for $1\leq i\leq n$,
\begin{align}\label{hyp:induction}
    \|c_h^{i-1}-c^{i-1}\|^2_{\Omega}+\Delta t \sum_{m=1}^{i} \tnorm{\boldsymbol{c}_h^{m-1}-\boldsymbol{c}^{m-1}}_c^2&\leq C \big((\Delta t)^2+ h^{2k_f}+h^{2k_c}\big).
\end{align}
Then $\boldsymbol{\zeta}_u^n$ satisfies:
\begin{subequations}
\label{eq:estimates_for_zeta}
\begin{align}
    \label{eq:estimate_zeta_u_n_s}
    \|\zeta_u^n\|^2_{\Omega^s}&+(\Delta t)^2\sum_{m=1}^n\|d_t\zeta_u^m\|^2_{\Omega^s}+\Delta t\sum_{m=1}^n\tnorm{\boldsymbol{\zeta}_u^m}_v^2 
    \leq C\big((\Delta t)^2+h^{2k_f}+h^{2k_c}\big),
    \\
    \label{eq:estimate_zeta_u_n_d}
    \|{\zeta}^n_u\|_{\Omega^d}^2 &\leq C''\big((\Delta t)^2+h^{2k_f}+h^{2k_c}\big),
\end{align}
\end{subequations}
where the constants depend on $\mu^*, \kappa_*, \mu_L, L_f^s, L_f^d, \sum_{j=1}^{{\dim}-1}\gamma^j$, and the regularity of $u_0,u,c_0$, and $c$ but are independent of the mesh size.
\end{thrm}
\begin{proof}
\textit{Proof of \cref{eq:estimate_zeta_u_n_s}:}\\
  Setting $(\boldsymbol{v}_h, \boldsymbol{q}_h) = (\boldsymbol{\zeta}_u^n,
  -\boldsymbol{\zeta}_p^n)$ in Theorem~\ref{thm:error_equation}, using $a(a-b)= \frac12(a^2 -b^2+(a-b)^2)$, and the coercivity of $a_h$ \cref{eq:a_h_coercivity} yields
  \begin{align}\label{eq:error_equation_est_B_sd_zeta}
  \frac{1}{2\Delta t} (\|\zeta_u^n\|^2_{\Omega^s}&-\|\zeta_u^{n-1}\|^2_{\Omega^s})+\dfrac{\Delta t}2\|d_t\zeta_u^n\|^2_{\Omega^s}+C_a\tnorm{\boldsymbol{\zeta}_u^n}_v^2\nonumber\\
  =&\sum\limits_{K\in\mathcal{T}^s}\int_K (\partial_t\,u^n-d_t\,\Pi_Vu^n)\cdot \zeta_u^n\,\dif x +a_h(\boldsymbol{c}_h^{n-1}; \boldsymbol{\xi}^n_u, \boldsymbol{\zeta}_u^n)+ [a_h(\boldsymbol{c}^n; \boldsymbol{u}^n,\boldsymbol{\zeta}_u^n)-a_h(\boldsymbol{c}_h^{n-1}; \boldsymbol{u}^n, \boldsymbol{\zeta}_u^n)] \nonumber\\
   &+\int\limits_{\Omega^s} [f^s(c_h^{n-1})-f^s(c^{n})]\cdot \zeta_u^n \dif  x + \int\limits_{\Omega^d} [\mathbb{K}^{-1}(c_h^{n-1})f^d(c_h^{n-1})-\mathbb{K}^{-1}(c^{n})f^d(c^{n})]\cdot \zeta_u^n\dif  x \nonumber\\
   =:&I_1+\ldots+I_5.
\end{align}
  Using \cref{ineq:bound_L2_by_tnorm_vs} and employing Young's inequality for some $\epsilon>0$,
  \begin{align}\label{ineq:I1}
    I_1
    \leq (\|\partial_t\,u^n -d_t u^n\|_{\Omega^s}+\|d_t\xi_u^n\|_{\Omega^s})\|\zeta_u^n\|_{\Omega^s}
    \leq C\big(\|\partial_t\,u^n -d_t u^n\|_{\Omega^s}^2 + \|d_t\xi_u^n\|_{\Omega^s}^2\big) +\epsilon\tnorm{\boldsymbol{\zeta}_u^n}_v^2.
  \end{align}
  %
  It follows from \cref{thm:SDT_a_h_continuityvp} and Young's inequality that
\begin{align}\label{ineq:I2}
    I_2
    \leq C\tnorm{ \boldsymbol{\xi}_u^n }_{v'}\tnorm{ \boldsymbol{\zeta}_u^n }_v \le C\tnorm{ \boldsymbol{\xi}_u^n }_{v'}^2+ \epsilon\tnorm{ \boldsymbol{\zeta}_u^n }_v^2.
\end{align}
  %
Consider now $I_3$:
\begin{align}\label{ineq:I3}
    I_3
    &=\sum_{K\in\mathcal{T}^s} \int_K 2[\mu(c^n)-\mu(c^{n-1}_h)] \varepsilon(u^n) : \varepsilon(\zeta_u^n) \dif x
     -\sum_{K\in\mathcal{T}^s} \int_{\partial K} 2[\mu(c^n)-\mu(c^{n-1}_h)] \varepsilon(u^n)n^s \cdot (\zeta_u^n-\bar{\zeta}_u^n)\dif s\nonumber\\
    &\quad + \int_{\Omega^d} [\mathbb{K}^{-1}(c^n)-\mathbb{K}^{-1}(c^{n-1}_h)] u^n\cdot \zeta_u^n \dif x 
     + \sum_{j=1}^{\rm{\dim}-1}\int_{\Gamma^I}  \gamma^j [\mu(c^n)-\mu(\bar{c}^{n-1}_h)](u^{sn}\cdot \tau^j)(\bar{\zeta}^n_u\cdot \tau^j)\dif s\nonumber\\
    &=:I_{31}+\hdots +I_{34}.
\end{align}
The first term on the right side of \cref{ineq:I3} can be bounded as follows using Lipschitz continuity of $\mu$, the generalized H\"{o}lder's inequality for integrals and sums, and Young's inequality:
\begin{align}\label{ineq:I31}
    I_{31}
     &\leq 2\mu_L\sum_{K\in\mathcal{T}^s} \|c^n-c^{n-1}_h\|_{K} \|\varepsilon(u^n)\|_{0,\infty,K} \|\nabla\zeta_u^n\|_{K} 
     \leq 2\mu_L\|\nabla u^n\|_{0,\infty,\Omega^s}\Big( \sum_{K\in\mathcal{T}^s} \|c^n-c^{n-1}_h\|_{K}^2\Big)^{1/2} \tnorm{\boldsymbol{\zeta}_u^n}_{v} \nonumber\\
     &\leq C\mu_L^2\|\nabla u^n\|^2_{0,\infty,\Omega^s} \sum_{K\in\mathcal{T}^s} \|c^n-c^{n-1}_h\|^2_{K} + \epsilon\tnorm{\boldsymbol{\zeta}_u^n}_{v}^2.
\end{align}
Next we bound $I_{32}$.
By Lipschitz continuity of $\mu$, \cref{ineq:trace_cont}, \cref{lem:BDM-iv}, generalized H\"{o}lder's inequality, and Young's inequality,
\begin{align}
    &|I_{32}|\leq 2\mu_L \sum_{K\in \mathcal{T}^s}\|c^n-c_h^{n-1}\|_{\partial K}\|\nabla u^n\|_{0,\infty, K}\|\zeta_u^n-\bar{\zeta}_u^n\|_{\partial K} \nonumber\\
    &\leq 2\mu_L \|\nabla u^n\|_{0,\infty, \Omega^s}\Big(\sum_{K\in \mathcal{T}^s}h_K\|c^n-c_h^{n-1}\|^2_{\partial K}\Big)^{1/2}\Big(\sum_{K\in \mathcal{T}^s}h_K^{-1}\|\zeta_u^n-\bar{\zeta}_u^n\|_{\partial K}^2\Big)^{1/2} \nonumber\\
    &\leq 2C \mu_L\|\nabla u^n\|_{0,\infty,\Omega^s} \Big(\sum_{K\in\mathcal{T}^s} ( \|c^n-c^{n-1}_h\|_K^2
  + h_K^2 \|c^n-c^{n-1}_h\|_{1,K}^2)\Big)^{1/2}  \tnorm{\boldsymbol{\zeta}_u^n}_{v} \nonumber\\
  &\leq C\mu_L^2\|\nabla u^n\|^2_{0,\infty,\Omega^s} \sum_{K\in\mathcal{T}^s} (\|c^n-c^{n-1}_h\|_K^2
  + h_K^2 \|c^n-c^{n-1}_h\|_{1,K}^2)  +\epsilon\tnorm{\boldsymbol{\zeta}_u^n}_{v}^2. \label{ineq:I32}
 \end{align}
We bound $I_{33}$ by using the assumption on $\mathbb{K}$ given in \cref{eq:assumption_kappa}:
\begin{align}
    I_{33}
    \leq \frac{\mu_L}{\kappa_*}\|c^n-c^{n-1}_h\|_{\Omega^d} \|u^n\|_{0,\infty,\Omega^d}\tnorm{\boldsymbol{\zeta}_u^n}_{v} \leq C(\kappa_*^{-1}\mu_L)^2 \|u^n\|^2_{0,\infty,\Omega^d}\|c^n-c^{n-1}_h\|^2_{\Omega^d}
     +  \epsilon\tnorm{\boldsymbol{\zeta}_u^n}_{v}^2. \label{ineq:I33}
\end{align}
Again by the Lipschitz property of $\mu$ and H\"{o}lder's inequality,
\begin{align}
    I_{34}
    &\leq \mu_L\|c^n-\bar{c}^{n-1}_h\|_{\Gamma^I} \big(\textstyle\sum\limits_{j=1}^{\rm{\dim}-1}  \gamma^j \big)^{1/2}\|u^{sn}\|_{0,\infty,\Gamma^I}\Big(\sum_{j=1}^{\rm{\dim}-1}\gamma^j\|\bar{\zeta}^n_u\cdot \tau^j\|_{\Gamma^I}^2\Big)^{1/2}\nonumber\\
    &\leq \mu_L \|c^n-\bar{c}^{n-1}_h\|_{\Gamma^I}\big(\textstyle\sum\limits_{j=1}^{\rm{\dim}-1}  \gamma^j \big)^{1/2}\|u^n\|_{0,\infty,\Omega^s} \tnorm{\boldsymbol{\zeta}_u^n}_{v}\nonumber\\
    &\leq C\mu_L^2\big(\textstyle\sum\limits_{j=1}^{\rm{\dim}-1}  \gamma^j \big)\|u^n\|^2_{0,\infty,\Omega^s} \|c^n-\bar{c}^{n-1}_h\|^2_{\Gamma^I}+\epsilon\tnorm{\boldsymbol{\zeta}_u^n}_{v}^2.\label{ineq:I34}
\end{align}
Combining \crefrange{ineq:I3}{ineq:I34},
\begin{align}\label{ineq:I3full}
    I_3\le &C\mu_L^2\|\nabla u^n\|^2_{0,\infty,\Omega^s} \sum_{K\in\mathcal{T}^s} (\|c^n-c^{n-1}_h\|_K^2
    + h_K^2 \|c^n-c^{n-1}_h\|_{1,K}^2) + C(\kappa_*^{-1}\mu_L)^2 \|u^n\|^2_{0,\infty,\Omega^d}\|c^n-c^{n-1}_h\|^2_{\Omega^d}\nonumber\\
    &+ C\mu_L^2\big(\textstyle\sum\limits_{j=1}^{\rm{\dim}-1}  \gamma^j \big)\|u^n\|^2_{0,\infty,\Omega^s} \|c^n-\bar{c}^{n-1}_h\|^2_{\Gamma^I} + 4\epsilon\tnorm{\boldsymbol{\zeta}_u^n}_{v}^2.
\end{align}
 Since $f^s$ is Lipschitz continuous in $c$, with Lipschitz constant $L^s_f$, and recalling  \cref{ineq:bound_L2_by_tnorm_vs},
 \begin{align}\label{ineq:I4}
     I_4
     &\le CL^s_f\|c_h^{n-1}-c^n\|_{\Omega^s}\tnorm{\boldsymbol{\zeta}_u^n}_{v}\le C(L^s_f)^2\|c_h^{n-1}-c^n\|_{\Omega^s}^2+\epsilon\tnorm{\boldsymbol{\zeta}_u^n}_{v}^2.
 \end{align} 
Since $f^d$ and $\mu$ are Lipschitz continuous in $c$, with Lipschitz continuity constants $L^d_f$ and $\mu_L$, respectively,
\begin{align}\label{ineq:I5}
     I_5
     =&\int_{\Omega^d}  \big(\mathbb{K}^{-1}(c_h^{n-1})[f^d(c_h^{n-1})-f^d(c^{n})]
     + [\mu(c_h^{n-1})-\mu(c^n)]\kappa^{-1}f^d(c^{n})\big)\cdot \zeta_u^n\dif  x  \nonumber\\
     &\le \big(K_{*}^{-1}L_f^d+\mu_L\kappa_{*}^{-1}\|f^d(c^n)\|_{0,\infty,\Omega^d}\big)\|c_h^{n-1}-c^{n}\|_{\Omega^d}\|\zeta_u^n\|_{\Omega^d}\nonumber\\
     &\leq C\big(K_{*}^{-1}L_f^d+\mu_L\kappa_{*}^{-1}\|f^d(c^n)\|_{0,\infty,\Omega^d}\big)^2\|c_h^{n-1}-c^{n}\|_{\Omega^d}^2+\epsilon\tnorm{\boldsymbol{\zeta}_u^n }_{v}^2.
\end{align}
Combining the above bounds for $I_1$ to $I_5$ with \cref{eq:error_equation_est_B_sd_zeta},
letting $\epsilon=\dfrac{C_a}{16}$ ($C_a$ is the coercivity constant), 
multiplying by $2\Delta t$, summing from $1$ to $n$, noting that $\zeta_u^0=0$, and applying \cref{ineq:taylor}, we obtain:
%
 \begin{align*}
    \|\zeta_u^n\|^2_{\Omega^s}+(\Delta t)^2&\sum_{m=1}^n\|d_t\zeta_u^m\|^2_{\Omega^s}+C_a\Delta t\sum_{m=1}^n\tnorm{\boldsymbol{\zeta}_u^m}_v^2\nonumber\\
     \leq  C\Big[& (\Delta t)^2\|\partial _{tt}u\|_{L^2(0,T;L^2(\Omega^s))}^2 +\|\partial_t\xi_u\|_{L^2(0,T;L^2(\Omega^s))}^2+2\Delta t\sum_{m=1}^n\tnorm{ \boldsymbol{\xi}_u^m }_{v'}^2\\
     &+\mu_L^2\Delta t\sum_{m=1}^n\Big(\|\nabla u^m\|^2_{0,\infty,\Omega^s} \sum_{K\in\mathcal{T}^s}h_K^2 \|c^m-c^{m-1}_h\|_{1,K}^2\Big)  \nonumber\\
   & +    \Delta t\sum_{m=1}^n\Big(\kappa_*^{-1}\mu_L\|u^m\|_{0,\infty,\Omega^d}+K_{*}^{-1}L_f^d+\dfrac{\mu_L}{\kappa_{*}}\|f^d(c^m)\|_{0,\infty,\Omega^d}\Big)^2\|c^m-c^{m-1}_h\|^2_{\Omega^d}\\
   & +\mu_L^2\big(\textstyle\sum\limits_{j=1}^{\rm{\dim}-1}  \gamma^j \big)\Delta t \sum_{m=1}^n\|u^m\|^2_{0,\infty,\Omega^s}\|\bar{c}^m-\bar{c}^{m-1}_h\|_{\Gamma^I}^2 \nonumber\\
     & +\Delta t\sum_{m=1}^n\big((L^s_f)^2+\mu_L^2\|\nabla u^m\|^2_{0,\infty,\Omega^s}\big)\|c_h^{m-1}-c^{m}\|_{\Omega^s}^2\Big].
\end{align*}
Next, using \cref{lem:BDM-iv} and \cref{eq:xi_u-2},
\begin{align*}
    \|\zeta_u^n\|^2_{\Omega^s}+(\Delta t)^2&\sum_{m=1}^n\|d_t\zeta_u^m\|^2_{\Omega^s}+C_a\Delta t\sum_{m=1}^n\tnorm{\boldsymbol{\zeta}_u^m}_v^2\\
    \leq  C\Big[&(\Delta t)^2\|\partial_{tt}u\|_{L^2(0,T;[L^2(\Omega^s)]^{\dim})}^2+h^{2k_f}\big(\|\partial_{t}u\|^2_{L^2(0,T;[H^{k_f}(\Omega^s)]^{\rm \dim})}+\|u\|_{\ell^2(0,T;[H^{k_f+1}(\Omega^s)]^{\rm \dim})}^2
    \\
    &+ \|u\|_{\ell^2(0,T;[H^{k_f}(\Omega^d)]^{\dim})}^2\big)+ \|\nabla u\|_{L^{\infty}(0,T;[L^{\infty}(\Omega^s)]^{\rm \dim})}^2 \Delta t \sum_{m=1}^n\sum_{K\in\mathcal{T}^s}\big( h_K^2 \|c^m-c^{m-1}_h\|_{1,K}^2\big)\nonumber\\
    &+\big(\|u\|_{L^{\infty}(0,T;[L^{\infty}(\Omega^d)]^{\rm \dim})}^2+ \|\nabla p\|_{L^{\infty}(0,T;[L^{\infty}(\Omega^d)]^{\rm \dim})}^2+1\big)\Delta t\sum_{m=1}^n\|c^m-c^{m-1}_h\|_{\Omega^d}^2\nonumber
    \\
    &+\|u\|_{L^{\infty}(0,T;[L^{\infty}(\Omega^s)]^{\rm \dim})}^2\Delta t\sum_{m=1}^n\|\bar{c}^m-\bar{c}^{m-1}_h\|_{\Gamma^I}^2\\
    &+\big(1+\|\nabla u\|^2_{L^{\infty}(0,T; [L^{\infty}(\Omega^s)]^{\rm \dim})}\big) \Delta t\sum_{m=1}^n\|c^m-c^{m-1}_h\|_{\Omega}^2\Big],\nonumber
\end{align*}
where the constant $C>0$ depends on $\mu_L, \kappa_*, K^{*}, \gamma^j$, $L_f^s$ but is independent of $h$ and $\Delta t$.
\Cref{eq:estimate_zeta_u_n_s} follows by \cref{eq:bound_cm_chm_minus_1_1,eq:bound_cm_chm_minus_1_omega_d,eq:bound_cm_chm_minus_chm_interface}, and assumption \cref{hyp:induction}.\\

\textit{Proof of \cref{eq:estimate_zeta_u_n_d}:}\\
Let $v_h^s=0$, $\bar{v}_h=0$, and $\boldsymbol{q}_h=\boldsymbol{0}$ in \cref{eq:error_equation_0}. Then 
 \begin{multline}\label{eq:darcyerror-1}
    a_h^d(c_h^{n-1}; \zeta^n_u, v_h)
    +b_h^d(\boldsymbol{\zeta}_p^{dn}, v_h)
    =a_h^d(c_h^{n-1}; \xi^n_u, v_h)+ a_h^d(c^n; u^n,v_h)-a_h^d(c_h^{n-1}; u^n, v_h)\\
    + \sum\limits_{K\in\mathcal{T}^d}\int_K  [\mathbb{K}^{-1}(c_h^{n-1})f^d(c_h^{n-1})-\mathbb{K}^{-1}(c^{n})f^d(c^{n})]\cdot v_h\dif x.
  \end{multline}
On the other hand, letting $\boldsymbol{v}_h=\boldsymbol{0}$ and $\boldsymbol{q}_h^s=\boldsymbol{0}$  in \cref{eq:error_equation_0}, we have
\begin{equation*}
    b_h^d(\boldsymbol{q}_h^{d}, \zeta_u^n)+ b_h^{I,d}(\bar{q}_h^{d}, \bar{\zeta}_u^n)=0,
  \end{equation*}
implying that
\begin{equation}\label{eq:darcyerror-2}
    b_h^d(\boldsymbol{q}_h^{d}, \zeta_u^n)= -b_h^{I,d}(\bar{q}_h^{d}, \bar{\zeta}_u^n)=\int_{\Gamma^I} \bar{q}_h^d \bar{\zeta}_u^n\cdot n^d \dif s.
\end{equation}
From \cite[Lemma 3.2]{Girault:2014}, since $u^{sn}-u_h^{sn}\in H^{\rm div}(\Omega^s)$, there exists $w\in H^{\rm div}(\Omega^d)$ such that $\nabla \cdot w =0$ in $\Omega^d$, $w\cdot n=0$ on $\Gamma^d$ and $w\cdot n^d=(\gamma(u^{sn})-u_h^{sn})\cdot n^d$ on $\Gamma^I$. 
With this choice of $w$, and using \cref{lem:BDM-i,lem:BDM-ii,lem:BDM-iii} we observe that  
\begin{align}\label{eq:darcyerror-3}
    b_h^d(\boldsymbol{q}_h^d, \Pi_V w)=\int_{\Gamma^I}\bar{q}_h^d(u^{sn}-u_h^{sn})\cdot n^d \dif s.
\end{align}
Adding \cref{eq:darcyerror-2,eq:darcyerror-3}, and recalling \cref{eq:masscons-2,eq:masscons-3}, the definition of $\bar{\Pi}_V$, and \cref{lem:BDM-ii}, we obtain
\begin{equation}\label{eq:darcyerror-4}
    b_h^d(\boldsymbol{q}_h^d, \zeta_u^n+\Pi_V w)=0 \quad \forall \boldsymbol{q}_h^d\in \boldsymbol{Q}_h^d.
\end{equation}
This leads us to consider $v_h=\zeta_u^n+\Pi_V w$ as test function in \cref{eq:darcyerror-1}. Using
\cref{eq:assumption_K,eq:darcyerror-4} we find that
 \begin{multline}\label{eq:darcyerror-5}
    \dfrac1{K^{*}}\|\zeta_u^n\|_{\Omega^d}^2\leq -a_h^d(c_h^{n-1}; \zeta^n_u, \Pi_V w)
   + a_h^d(c_h^{n-1}; \xi^n_u, \zeta_u^n+\Pi_V w)
   \\
   + \Big(a_h^d(c^n; u^n,\zeta_u^n+\Pi_V w)-a_h^d(c_h^{n-1}; u^n, \zeta_u^n+\Pi_V w)\Big)\\
    + \sum\limits_{K\in\mathcal{T}^d}\int_K  [\mathbb{K}^{-1}(c_h^{n-1})f^d(c_h^{n-1})-\mathbb{K}^{-1}(c^{n})f^d(c^{n})]\cdot (\zeta_u^n+\Pi_V w)\dif x   =: A_1+A_2+A_3+A_4.
  \end{multline}
We will bound $A_1$ to $A_4$ by a series of Cauchy--Schwarz, H\"{o}lder's, triangle, and Young's inequalities together with the properties of $\mu$ and $\kappa$.
First,
\begin{align*}
    A_1\leq K_*^{-1}\|\zeta_u^n\|_{\Omega^d}\|\Pi_V w\|_{\Omega^d}\leq \epsilon \|\zeta_u^n\|_{\Omega^d}^2 + CK_*^{-2}\|\Pi_V w\|_{\Omega^d}^2,
\end{align*}
and
\begin{align*}
    A_2&\leq K_*^{-1}\|\xi_u^n\|_{\Omega^d}\|\zeta_u^n+\Pi_V w\|_{\Omega^d}
    \leq  CK_*^{-2}\|\xi_u^n\|_{\Omega^d}^2+\epsilon\|\zeta_u^n\|^2_{\Omega^d}+K_*^{-1}\|\xi_u^n\|_{\Omega^d}\|\Pi_V w\|_{\Omega^d}.
\end{align*}
Following the proof of \cref{ineq:I33}, 
\begin{align*}
    A_3
    \leq &  \mu_L\kappa_*^{-1}\|c^n-c_h^{n-1}\|_{\Omega^d}\|u^n\|_{0,\infty,\Omega^d}\big(\|\zeta_u^n\|_{\Omega^d}+\|\Pi_V w\|_{\Omega^d}\big)\\
    \leq &  C\mu_L^2\kappa_*^{-2}\|c^n-c_h^{n-1}\|^2_{\Omega^d}\|u^n\|^2_{0,\infty,\Omega^d}+\epsilon\|\zeta_u^n\|^2_{\Omega^d}
    +\mu_L\kappa_*^{-1}\|c^n-c_h^{n-1}\|_{\Omega^d}\|u^n\|_{0,\infty,\Omega^d}\|\Pi_V w\|_{\Omega^d},
\end{align*}
while as the proof of \cref{ineq:I5},
\begin{align*}
     A_4\le& \big(\mu_L\kappa_*^{-1}\|f^d(c^n)\|_{0,\infty,\Omega^d}+K_{*}^{-1}L_f^d\big)\|c_h^{n-1}-c^{n}\|_{\Omega^d}\big(\|\zeta_u^n\|_{\Omega^d}+\|\Pi_V w\|_{\Omega^d}\big)\nonumber\\
     \leq& C\big(\mu_L\kappa_*^{-1}\|f^d(c^n)\|_{0,\infty,\Omega^d}+K_{*}^{-1}L_f^d\big)^2\|c_h^{n-1}-c^{n}\|^2_{\Omega^d}+\epsilon\|\zeta_u^n\|^2_{\Omega^d}\nonumber\\
     &+\big(\mu_L\kappa_*^{-1}\|f^d(c^n)\|_{0,\infty,\Omega^d}+K_{*}^{-1} L_f^d\big)\|c_h^{n-1}-c^{n}\|_{\Omega^d}\|\Pi_V w\|_{\Omega^d}.\nonumber
\end{align*}
Combining the bounds of $A_1$ to $A_4$ with \cref{eq:darcyerror-5}, choosing $\epsilon=1/(8K^*)$, and applying another set of Young's inequalities,
we obtain
\begin{align}
\label{eq:boundzetaund}
 \|\zeta_u^n\|_{\Omega^d}^2\leq  C\big(\|\Pi_Vw\|_{\Omega^d}^2+\|\xi_u^n\|_{\Omega^d}^2+\|c^n-c_h^{n-1}\|_{\Omega^d}^2\big),
\end{align}
where $C$ depends on the problem parameters $\kappa_*,K_*,K^*,\mu_L, L_f^d$, and the regularity of $u, p$ and $c$. Noting that from \cref{lem:BDM-iv} and \cite[Theorem 3.3]{Girault:2014}, \cite[Lemma 10]{Cesmelioglu:2020},
\begin{equation}\label{eq:darcyerror-6}
\|\Pi_Vw\|_{\Omega^d}\leq C\|w\|_{\Omega^d}\leq Ch^{k_f}\big(\|u^n\|_{k_f+1,\Omega^s}+\|u\|_{k_f,\Omega^d}+\|g_p-g_i\|_{k_f,\Omega^d}\big).
\end{equation}
Then \cref{eq:boundzetaund,eq:darcyerror-6,lem:BDM-iv,ineq:cm-chm-1} imply
\begin{equation*}
 \|\zeta_u^n\|_{\Omega^d}^2\leq  C h^{2k_f}\big(\|u^n\|_{k_f+1,\Omega^s}^2+\|u^n\|^2_{k_f,\Omega^d}+\|g_p-g_i\|^2_{k_f,\Omega^d}
 +(\Delta t)^2\|d_tc^n\|_{\Omega^d}^2+\|c^{n-1}-c_h^{n-1}\|_{\Omega^d}^2\Big).
\end{equation*}
\Cref{eq:estimate_zeta_u_n_d} is now a consequence of the assumptions on the regularity of the exact solution and \cref{hyp:induction}.
\end{proof}
The following is a straightforward consequence of Theorem~\ref{thm:vel-error}. 
\begin{crllr} \label{cor:velerror}
Let $u^n$ and $u_h^n$ be as defined in Theorem~\ref{thm:vel-error}. Then for all $n\geq 1$,
\begin{subequations}
\begin{align}
\label{ineq:velerror}
\|u^n-u_h^n\|_{\Omega}&\leq 
C(\Delta t+h^{k_f}+h^{k_c}),
\\
\label{ineq:velerrorbdry}
\sum_{K\in \mathcal{T}_h}h_K\|u^n-u_h^n\|_{\partial K}^2&\leq 
C((\Delta t)^2+h^{2k_f}+h^{2k_c}).
\end{align}
\end{subequations}
\end{crllr}
Before moving on to the next section, we note another consequence of \cref{eq:estimates_for_zeta} that will prove useful in analysis later on. 
\begin{crllr}
Let $u$ denote the velocity solution to \cref{eq:system,eq:interface,eq:initial} satisfying the assumptions in Theorem~\ref{thm:vel-error} with $k_f, k_c\geq {\dim}-1$. Suppose $\Delta t\leq Ch_{k_f}^{{\dim}/2}$. Then for each $n\geq 1$, the discrete velocity $u_h^n$ that solves \cref{eq:hdgsd_sequence} satisfies
\begin{equation}
    \label{eq:estimate_uh_infty}
    \|u_h^n\|_{0,\infty,\Omega}\leq C,
  \end{equation}
where $C$ depends on $\mu^*, \kappa_*, \mu_L, L_f^s, L_f^d, \sum_{j=1}^{{\dim}-1}\gamma^j$, and the regularity of $u_0,u,c_0$, and $c$, but is independent of $h$, $n$ and $\Delta t$.
\end{crllr}
\begin{proof}
The proof is the same as in \cite[Lemma 1]{Cesmelioglu:2021}. Using \cref{eq:estimate_zeta_u_n_s,eq:estimate_zeta_u_n_d}, we have
\begin{equation}
   \|\zeta_u^n\|_{\Omega}^2 \leq C((\Delta t)^2 + h^{2k_f} + h^{2k_c}).
\end{equation}
Combining the above estimate with \cref{ineq:inverse-0-infty,lem:BDM-iv}, for each $K\in \mathcal{T}$, we have
\begin{align*}
    \|u_h^n\|_{0,\infty,K}
    &\leq \|\zeta_u^n\|_{0,\infty,K} + \|\Pi_V u^n\|_{0,\infty,K}
    \leq C h_K^{-{\dim}/2}\|\zeta_u^n\|_{K} + \|\Pi_V u^n\|_{0,\infty,K} 
    \\
    &\leq Ch_K^{-{\dim}/2}\big((\Delta t) + h^{k_f} + h^{k_c}\big) + \|u^n\|_{0,\infty,K}\\
     &\leq Ch_K^{-{\dim}/2}\big((\Delta t) + h^{k_f} + h^{k_c}\big) + \|u\|_{L^{\infty}(0,T;[L^{\infty}(K)]^{\dim})}.
\end{align*}
The result \cref{eq:estimate_uh_infty} follows by using the assumption on $\Delta t$ and by taking the maximum over $K\in\mathcal{T}$. 
\end{proof}
 \begin{rmrk}
    The restrictions on the polynomial degree and the time step are not necessary if we assume that $|D(u)|\leq \bar{D}$ for some positive constant $\bar{D}$ as in \cite[2.12]{Riviere:2014}. It is also possible to avoid these restrictions by using an approach involving a cutoff operator on the velocity solution, as in \cite{Sun:2002,Sun:2005}, if one is interested in lower order approximations.  
\end{rmrk}
\begin{rmrk}
    Compatibility, as defined in \cite{Dawson:2004}, can be achieved by choosing $k_c=k_f-1$ \cite{Cesmelioglu:2021}. However, the requirement that $k_c\geq {\rm dim-1}$ implies $k_f\geq {\dim}$. Therefore, when ${\dim}=2$, our theory supports compatibility only for $k_f\geq 2$ and for $k_f\geq 3$ when ${\dim}=3$.
\end{rmrk}
%
\subsection{Error estimate for the pressure}
\label{sec:error-estimates-pressure}
In this section, we briefly discuss the a priori error estimate for the pressure approximation. 
\begin{lmm} 
\label{lem:zetapmp2bound}
Suppose that the assumptions in Theorem~\ref{thm:vel-error} hold and that $p$ and $\boldsymbol{p}_h$ are the pressure solutions to \cref{eq:system,eq:interface,eq:initial} and \cref{eq:hdgsd_sequence}, respectively. Then 
\begin{align}\label{ineq:press}
    \Delta t&\sum\limits_{m=1}^n \tnorm{\boldsymbol{\zeta}_p^m}_{p}^2  \leq C
    \Big(\Delta t\sum\limits_{m=1}^n\tnorm{\boldsymbol{\zeta}_u^m}_{v}^2
    +\Delta t\sum\limits_{m=1}^n\| d_t\zeta_u^m\|_{\Omega^s}^2+\Delta t\sum\limits_{m=1}^n\tnorm{\boldsymbol{\zeta}_c^{m-1}}_{c}^2\\
    &+\Delta t \sum\limits_{m=1}^n\sum_{K\in\mathcal{T}} \|c^{m-1}-c^{m-1}_h\|_K^2
    + \Delta t \sum\limits_{m=1}^n\sum_{K\in\mathcal{T}^s}h_K^2 \|c^{m-1}-c^{m-1}_h\|_{1,K}^2+ (\Delta t)^2+h^{2k_f} + h^{2k_c+1}\Big).\nonumber
\end{align}
\end{lmm}
\begin{proof}
Setting  $\boldsymbol{q}_h = 0$ in the error equation in Theorem~\ref{thm:error_equation}, 
we obtain:
\begin{align}\label{eq:Hsum}
  \sum_{j=s,d}  \big(b_h^j(\boldsymbol{\zeta}_p^{jn}, v_h) &+ b_h^{I,j}(\bar{\zeta}_p^{jn}, \bar{v}_h)\big)
    \nonumber
    \\
    =&\big(a_h(\boldsymbol{c}^{n-1}_h;\boldsymbol{\xi}_u^n, \boldsymbol{v}_h)-a_h(\boldsymbol{c}^{n-1}_h;\boldsymbol{\zeta}_u^n, \boldsymbol{v}_h) \big) + \big(a_h(\boldsymbol{c}^n; \boldsymbol{u}^n, \boldsymbol{v}_h )-a_h(\boldsymbol{c}^{n-1}_h; \boldsymbol{u}^n,\boldsymbol{v}_h)\big) 
    \nonumber
    \\
    &+\int\limits_{\Omega^s} [f_s(c^{n-1}_h)-f_s(c^{n})]\cdot v_h \dif  x   + \int\limits_{\Omega^d} [\mathbb{K}^{-1}(c_h^{n-1})f^d(c^{n-1}_h)-\mathbb{K}^{-1}(c^{n})f^d(c^{n})]\cdot v_h\dif  x 
    \nonumber\\
    &- \sum\limits_{K\in\mathcal{T}^s}\int_K (d_tu_h^n-\partial_t\,u^n)\cdot v_h\,\dif x
    \nonumber \\
    =&: H_1+\ldots+H_5.
\end{align}
By \cref{thm:SDT_a_h_continuityvp} and Young's inequality, and using that $\tnorm{\cdot}_v$ and $\tnorm{\cdot}_{v'}$ are equivalent on $\boldsymbol{V}_h$, we have
\begin{equation}\label{eq:H1}
    H_1
    \leq C(\tnorm{\boldsymbol{\zeta}_u^n}_{v}
    +\tnorm{\boldsymbol{\xi}_u^n}_{v'})\tnorm{\boldsymbol{v}_h}_{v}.
\end{equation}
Following the proof of \cref{ineq:I3full}, we can show that
\begin{align}\label{eq:H2}
    H_2
    \leq C\Big( &\mu_L\|\nabla u^n\|_{0,\infty,\Omega^s} \Big(\sum_{K\in\mathcal{T}^s} (\|c^n-c^{n-1}_h\|_K^2
    + h_K^2 \|c^n-c^{n-1}_h\|_{1,K}^2)\Big)^{1/2} \nonumber\\
    &+ \kappa_*^{-1}\mu_L \|u^n\|_{0,\infty,\Omega^d}\|c^n-c^{n-1}_h\|_{\Omega^d}+ \mu_L\big(\textstyle\sum\limits_{j=1}^{\rm{\dim}-1}  \gamma^j \big)^{1/2}\|u^n\|_{0,\infty,\Omega^s} \|c^n-\bar{c}^{n-1}_h\|_{\Gamma^I} \Big)\tnorm{\boldsymbol{v}_h}_{v}.
\end{align}
As in \cref{ineq:I4,ineq:I5}, we find
\begin{align}
\label{eq:H3}
     H_3
     &\le L^s_f\|c_h^{n-1}-c^n\|_{\Omega^s}\tnorm{\boldsymbol{v}_h}_{v},
     \\
\label{eq:H4}
     H_4
     &\le \big(\mu_L\kappa_{*}^{-1}\|f^d(c^n)\|_{0,\infty,\Omega^d}+K_{*}^{-1}L_f^d\big)\|c_h^{n-1}-c^{n}\|_{\Omega^d}\tnorm{\boldsymbol{v}_h}_v.
\end{align}
Using Cauchy--Schwarz and triangle inequalities, and \cref{lem:BDM-iv},
\begin{align}\label{eq:H5}
    H_5
    &= -\sum\limits_{K\in\mathcal{T}^s}\int_K \big(d_t\zeta_u^n+(d_t\Pi_Vu^n-\partial_t\Pi_Vu^n)-\partial_t\xi_u^n\big)\cdot v_h\,\dif x \nonumber\\
    &\leq \big(\| d_t\zeta_u^n\|_{\Omega^s}+\|d_t\Pi_Vu^n-\partial_t\Pi_Vu^n\|_{\Omega^s}+C h^{k_f}
    \|\partial_t u^n\|_{k_f, \Omega^s}\big)\tnorm{\boldsymbol{v}_h}_v.
\end{align}
Therefore, combining \eqref{eq:Hsum}-\eqref{eq:H5}, dividing both sides by $\tnorm{\boldsymbol{v}_h}_{v}$, taking the supremum over $\boldsymbol{v}_h\in \boldsymbol{V}_h$, and using Theorem~\ref{thm:infsup}, we obtain
\begin{align*}
  c_{\inf}^{\star}\tnorm{\boldsymbol{\zeta}_p^n}_{p} 
    \leq  C\Big(&\tnorm{\boldsymbol{\zeta}_u^n}_{v}
    +\tnorm{\boldsymbol{\xi}_u^n}_{v'}
    \\
    &+\mu_L\|\nabla u^n\|_{0,\infty,\Omega^s} \Big(\sum_{K\in\mathcal{T}^s} (\|c^n-c^{n-1}_h\|_K^2
    + h_K^2 \|c^n-c^{n-1}_h\|_{1,K}^2)\Big)^{1/2} \nonumber\\
    &+ \mu_L\big(\textstyle\sum\limits_{j=1}^{\rm{\dim}-1}  \gamma^j \big)^{1/2}\|u^n\|_{0,\infty,\Omega^s} \|c^n-\bar{c}^{n-1}_h\|_{\Gamma^I}+ L^s_f\|c_h^{n-1}-c^n\|_{\Omega^s}
    \\
    &+\big(\mu_L\kappa_{*}^{-1}(\|f^d(c^n)\|_{0,\infty,\Omega^d}+\|u^n\|_{0,\infty,\Omega^d})+K_{*}^{-1}L_f^d\big)\|c_h^{n-1}-c^{n}\|_{\Omega^d}\\
    &+\| d_t\zeta_u^n\|_{\Omega^s}+\|d_t\Pi_Vu^n-\partial_t\Pi_Vu^n\|_{\Omega^s}+C h^{k_f}
    \|\partial_t u^n\|_{k_f, \Omega^s}\Big).
\end{align*}
Squaring both sides, multiplying by $(c_{\rm inf}^*)^{-2}\Delta t$, summing from $1$ to $n$, using \cref{ineq:taylor1,eq:xi_u-2},
stability of $\Pi_V$, and the regularity assumptions on $u^s$, $u^d$, and $\nabla p^d$ yields
\begin{align*}
    \Delta t\sum\limits_{m=1}^n \tnorm{\boldsymbol{\zeta}_p^m}_{p}^2  \leq& C
    \Big(\Delta t\sum\limits_{m=1}^n\tnorm{\boldsymbol{\zeta}_u^m}_{v}^2
    +\Delta t\sum\limits_{m=1}^n\| d_t\zeta_u^m\|_{\Omega^s}^2\\
    &+\Delta t \sum\limits_{m=1}^n\sum_{K\in\mathcal{T}} \|c^m-c^{m-1}_h\|_K^2
    + \Delta t \sum\limits_{m=1}^n\sum_{K\in\mathcal{T}^s}h_K^2 \|c^m-c^{m-1}_h\|_{1,K}^2
    \\
    &+\Delta t\sum\limits_{m=1}^n \|c^m-\bar{c}^{m-1}_h\|^2_{\Gamma^I}+ (\Delta t)^2\|\partial_{tt}u\|_{L^2(0,T;[L^2(\Omega^s)]^{\dim})}^2\nonumber
   \\
   &+h^{2k_f}\big(\|u\|_{\ell^2(0,T;[H^{k_f+1}(\Omega^s)]^{\dim})}^2+\|u\|_{L^2(0,T;[H^{k_f}(\Omega^d)]^{\dim})}^2
   + \|\partial_t u\|_{\ell^2(0,T;[H^{k_f}(\Omega^s)]^{\dim}}^2\big)\Big).
\end{align*}
Therefore, the result follows by \cref{eq:bound_cm_chm_minus_1_1,eq:bound_cm_chm_minus_1_omega_d,eq:bound_cm_chm_minus_chm_interface} under the assumptions on the exact solution given in Theorem~\ref{thm:vel-error}.
\end{proof}
\begin{rmrk}
An immediate consequence Lemma~\ref{lem:zetapmp2bound} and Theorem~\ref{thm:vel-error} is 
\begin{align*}
    \Delta t\sum\limits_{m=1}^n \tnorm{\boldsymbol{\zeta}_p^m}_{p}^2  \leq
    C\Big(&\Delta t + (\Delta t)^{-1}(h^{2k_f}+h^{2k_c})
    + \Delta t\sum\limits_{m=1}^n\tnorm{\boldsymbol{\zeta}_c^{m-1}}_{c}^2\\
    &+\Delta t \sum\limits_{m=1}^n\sum_{K\in\mathcal{T}} \|c^{m-1}-c^{m-1}_h\|_K^2
    + \Delta t \sum\limits_{m=1}^n\sum_{K\in\mathcal{T}^s}h_K^2 \|c^{m-1}-c^{m-1}_h\|_{1,K}^2\Big).\nonumber
\end{align*}
This loss of $\Delta t$ is due to $\Delta t\sum\limits_{m=1}^n\| d_t\zeta_u^m\|_{\Omega^s}^2$ in \cref{ineq:press}. However, an improved estimate can be obtained by bounding this term as follows: testing \cref{eq:error_equation_0} with $\boldsymbol{v}_h=d_t\boldsymbol{\zeta}_u^m$, multiplying by $\Delta t$, summing from $m=1$ to $n$, employing a summation-by-parts formula on the terms on the right hand side that are contained in $a_h(\cdot;\cdot, d_t\boldsymbol{\zeta}_u^m)$ to transfer the discrete time derivative on $d_t\boldsymbol{\zeta}_u^m$ to the other terms, and assuming that the exact solution is sufficiently smooth in time, leads to 
\begin{equation*}
    \Delta t\sum\limits_{m=1}^n\| d_t\zeta_u^m\|_{\Omega^s}^2\leq C(h^{2k_f}+h^{2k_c}+(\Delta t)^{2}).    
\end{equation*}
We do not provide the details of this proof here, but instead refer to \cite[p.42]{Chaabane:2017}.
\end{rmrk}
%
\subsection{Existence and uniqueness of the concentration solution}
\label{ss:existuniq_ch}
In this section, we will prove existence and uniqueness of the discrete concentration solution $\boldsymbol{c}_h^n\in \boldsymbol{C}_h$ to \cref{eq:hdgtr_sequence}. First observe that assumption \cref{eq:upperboundDuabs} on $D$ implies that for $v \in [L^\infty(\Omega^d)]^{\dim}$,
\begin{equation}\label{ineq:Dv-1}
    \|D(v)\|_{0,\infty,\Omega^d} \leq C(1 + \|v\|_{0,\infty,\Omega^d}),
\end{equation}
and together with \cref{ineq:trace_cont_inf} that 
\begin{equation}\label{ineq:Dv-2}
    \|D(v)\|_{0,\infty,\partial K}\le C(1 + \|v\|_{0,\infty,\partial K}) \leq C \quad \forall v \in [W^{1,\infty}(K)]^{\dim}, K \in \mathcal{T}^d,
\end{equation}
where $C$ depends on $\|v\|_{1,\infty,K}$.
Therefore by \cref{eq:estimate_uh_infty,ineq:Dv-1,ineq:Dv-2}, 
there exists a constant $\widetilde{D}_{\max}>0$ that depends on $d$ and the upper bound in \cref{eq:estimate_uh_infty} such that
\begin{equation}\label{ineq:Dtildemax}
\|\widetilde{D}(u_h^n)\|_{0,\infty, \Omega}\leq \widetilde{D}_{\max},\quad \|\widetilde{D}(u_h^n)\|_{0,\infty,\partial K}\leq \widetilde{D}_{\max}.
\end{equation}
With \cref{ineq:Dtildemax}, the following coercivity result can be proved following the same steps as the proofs of \cite[Lemmas 2 and 3]{Cesmelioglu:2021}.
\begin{thrm}[Coercivity of $B_h^{tr}(u_h^n; \boldsymbol{w}_h, \boldsymbol{w}_h)$]\label{thm:B_h_tr_coercivity} 
    There exists a constant $\beta_0 > 0$ such that if $\beta_{tr} > \beta_0^{tr}$, then for all $\boldsymbol{w}_h \in \boldsymbol{C}_h$,
 \begin{equation}  \label{eq:B_h_a_d_coercivity}
  B_h^{tr}(u_h^n;\boldsymbol{w_h}, \boldsymbol{w_h})\geq C_{tr}\tnorm{\boldsymbol{w}_h}_c^2+\frac{1}{2} \sum_{K\in \mathcal{T}^d} \int_K \nabla \cdot u_h^n w_h^2\,\dif x,
\end{equation}
where $C_{tr}>0$ is a constant that depends on $d, D_{\min}$, and the upper bound in \cref{eq:estimate_uh_infty}.
\end{thrm}
Now that we have coercivity, we proceed with the existence and stability proof for the discrete concentration.
\begin{thrm}\label{thm:stable_SDT_1} 
Let $c_0\in L^2(\Omega)$ and $g_i-g_p\in L^{\infty}(0,T;L^{\infty}(\Omega^d))$. Let $n\geq 1$ and let $\boldsymbol{u}_h^n$ be the solution to \cref{eq:hdgsd_sequence} that satisfies \cref{eq:estimate_zeta_u_n_s,eq:estimate_zeta_u_n_d}. If $d_n\Delta t<1$, where $d_n=:\frac1{\phi_*}(1 +C\|g_i^n-g_p^n\|_{0,\infty,\Omega^d})$, then there exists  a unique solution $\boldsymbol{c}_h^n\in \boldsymbol{C}_h$ to \cref{eq:hdgtr_sequence}.
Furthermore, if $K:=\sum\limits_{m=1}^n d_m/(1-\Delta t d_m)$, then
\begin{equation}\label{ineq:chbd}
     \phi_* \|c_h^n\|_{\Omega}^2
  + C_{tr}\Delta t\sum\limits_{m=1}^n\tnorm{\boldsymbol{c}_h^m}_c^2 \leq e^{K\Delta t}\big(\phi_*\|c_{0}\|_{\Omega}^2 +\|g_i\|_{\ell^2(0,T;L^2(\Omega^d))}^2\big).
\end{equation}
\end{thrm}
\begin{proof}
Let $\boldsymbol{w}_h=\boldsymbol{c}_h^n$ in \cref{eq:hdgtr_sequence}.
From the algebraic inequality $(a-b)a\geq \frac{1}{2}(a^2-b^2)$, \cref{eq:B_h_a_d_coercivity}, and the assumption that $0\leq c_I\leq 1$ a.e., we have
\begin{align*}
\frac{\phi_*}{2\Delta t}\, \big(\|c_h^n\|_{\Omega}^2-\|c_h^{n-1}\|_{\Omega}^2\big) 
  + C_{tr}\tnorm{\boldsymbol{c}_h^n}_c^2 +\int_{\Omega^d}g_p^n(c_h^n)^2\,\dif x\nonumber
&\leq  \sum_{K\in \mathcal{T}^d} \int_K c_{I} g_{i}^n  c_h^n \dif  x -\frac{1}{2} \sum_{K\in \mathcal{T}^d} \int_K \nabla \cdot u_h^n (c_h^n)^2\,\dif x \nonumber\\
  &\le \|g_i^n\|_{\Omega^d}\|c_h^n\|_{\Omega^d} +\frac{1}{2}\|\nabla\cdot u_h^n\|_{0,\infty,\Omega^d}\|c_h^n\|_{\Omega^d}^2.\nonumber
\end{align*}
Multiplying this inequality by $2\Delta t$, summing from $m=1$ to $n$, noting that $g_p\geq 0$, 
and recalling \cref{eq:masscons-1} with stability of the $L^2$-projections $\Pi_C$ and $\Pi_Q$, we obtain
\begin{align*}
\phi_* \|c_h^n\|_{\Omega}^2 
  + & 2C_{tr}\Delta t\sum\limits_{m=1}^n\tnorm{\boldsymbol{c}_h^m}_c^2 
  \\
  &\le \phi_*  \|c_h^{0}\|_{\Omega}^2 + \Delta t\sum\limits_{m=1}^n (\|g_i^m\|_{\Omega^d}^2+\|c_h^m\|_{\Omega^d}^2) +\Delta t\sum\limits_{m=1}^n\|\Pi_Q(g_i^m-g_p^m)\|_{0,\infty,\Omega^d}\|c_h^m\|_{\Omega^d}^2
  \nonumber
  \\
  &\le \phi_*  \|c_{0}\|_{\Omega}^2 +\|g_i^m\|_{\ell^2(0,T;L^2(\Omega^d))}^2 +\Delta t\sum\limits_{m=1}^n(1 +C\|g_i^m-g_p^m\|_{0,\infty,\Omega^d})\|c_h^m\|_{\Omega}^2.
  \nonumber
\end{align*}
\Cref{ineq:chbd} follows after applying Gr\"{o}nwall's inequality \cite[Lemma 27]{Layton:book}.
This stability bound then implies the existence of a unique solution since the system is finite dimensional and linear.
\end{proof}
\subsection{Error estimate for the discrete concentration}
This section is devoted to proving an error estimate for the discrete concentration.
\begin{lmm}[Error equation for \cref{eq:hdgtr_sequence}]
\label{lem:erroreq-c}
\begin{align}
\sum_{K\in \mathcal{T}}\int_K \phi\, d_t \zeta_c^n w_h \dif  x 
  +& B_h^{tr}(u_h^n; \boldsymbol{\zeta}_c^n, \boldsymbol{w}_h)
+ \int_{\Omega^d}g_p^n\zeta_c^n w_h dx\\
  =&\sum_{K\in \mathcal{T}}\int_K \phi\, d_t \xi_c^n w_h \dif  x-\sum_{K\in \mathcal{T}}\int_K \phi\,(d_tc^n-\partial_t c^n) w_h \dif  x +B_h^{tr}(u_h^n;\boldsymbol{\xi}_c^n, \boldsymbol{w}_h)\nonumber\\
  &+B_h^{a}(u^n-u_h^n;  \boldsymbol{c}^n, \boldsymbol{w}_h)+B_h^d(u^n; \boldsymbol{c}^n, \boldsymbol{w}_h)-B_h^d(u_h^n; \boldsymbol{c}^n, \boldsymbol{w}_h)+ \int_{\Omega^d}g_p^n\xi_c^n w_h dx.\nonumber
\end{align}
\end{lmm}
\begin{proof}
By Lemma~\ref{lem:consistency}, for $t=t^n$, we have
\begin{equation}\label{eq:c-consistency_at_t1}
      \sum_{K\in \mathcal{T}}\int_K \phi\, \partial_t c^n w_h \,\dif x 
  + B_h^{tr}(u^n; \boldsymbol{c}^n, \boldsymbol{w}_h)+\int_{\Omega^d}g_p^nc^nw_h\,\dif x
  = \sum_{K\in \mathcal{T}^d} \int_K c_{I}\,g_{i}^n\, w_h \dif x\quad
  \forall \boldsymbol{w}_h\in \boldsymbol{C}_h,
\end{equation}
where $\boldsymbol{c}^n = (c^n, \gamma(c^n))$. 
Subtracting \cref{eq:c-consistency_at_t1} from \cref{eq:hdgtr_sequence} yields that for all $\boldsymbol{w}_h\in \boldsymbol{C}_h$,
 \begin{equation}\label{eq:error_hdgtr_sequence_1}
  \sum_{K\in \mathcal{T}}\int_K \phi\, (d_t c_h^n-\partial_t c^n) w_h \dif  x 
  + B_h^{tr}(u_h^{n}; \boldsymbol{c}_h^{n}, \boldsymbol{w}_h)-B_h^{tr}(u^n; \boldsymbol{c}^n, \boldsymbol{w}_h) 
  + \int_{\Omega^d}g_p^n(c_h^n-c^n) w_h \dif x  = 0.
\end{equation}
Next, we rewrite the $B_h^{tr}$ terms in \cref{eq:error_hdgtr_sequence_1} by observing that $B_h^{tr}$ is linear in the second slot and that $B_h^a$ is linear in the first slot:
\begin{align}\label{eq:c_splitting_error}
B_h^{tr}(u_h^{n}; \boldsymbol{c}_h^{n}, \boldsymbol{w}_h)&-B_h^{tr}(u^n; \boldsymbol{c}^n, \boldsymbol{w}_h)\\
=&B_h^{tr}(u_h^{n}; \boldsymbol{c}_h^{n}-\boldsymbol{c}^{n}, \boldsymbol{w}_h)
+B_h^{a}(u_h^n-u^n; \boldsymbol{c}^n, \boldsymbol{w}_h)+\big[B_h^{d}(u_h^n; \boldsymbol{c}^n, \boldsymbol{w}_h)-B_h^{d}(u^n; \boldsymbol{c}^n, \boldsymbol{w}_h)\big].\nonumber
\end{align}
Using \cref{eq:err-split-c}, again the linearity of $B_h^{tr}$ in the second slot, and  \cref{eq:error_hdgtr_sequence_1,eq:c_splitting_error} completes the proof.
\end{proof}

\begin{thrm}
In addition to the assumptions in Theorem~\ref{thm:vel-error}, suppose that 
\begin{align*}
    c_0\in H^{k_c}(\Omega), \quad &c\in L^2(0,T;H^{k_c+1}(\Omega))\cap L^{\infty}(0,T;W^{1,\infty}(\Omega)),\quad
\partial_tc\in L^2(0,T;H^{k_c}(\Omega)),\\
&\partial_{tt}c \in L^2(0,T;L^2(\Omega)), \quad g_i, g_p\in L^{\infty}(0,T;L^{\infty}(\Omega^d)).
\end{align*}
Then for sufficiently small $\Delta t$,
  \begin{equation}
      \|\zeta_c^n\|_{\Omega}^2+\Delta t \sum\limits_{m=1}^{n}\tnorm{\boldsymbol{\zeta}_c^m}_c^2 \leq C(h^{2k_f} + h^{2k_c} + (\Delta t)^2),
  \end{equation}
where $C$ depends on $\phi_*, \phi^*, d, D$ and the regularity of the solution but is independent of $h$ and $\Delta t$.
\end{thrm}
\begin{proof}
Setting $\boldsymbol{w}_h = \boldsymbol{\zeta}_c^n$ in Lemma~\ref{lem:erroreq-c}, using the inequality $a(a-b)\ge \frac{a^2-b^2}{2}$, and Theorem~\ref{thm:B_h_tr_coercivity}, 
\begin{align*}
\frac{\phi_*}{2\Delta t}(\|\zeta_c^n\|_{\Omega}^2&-\|\zeta_c^{n-1}\|_{\Omega}^2)+C_{tr}\tnorm{\boldsymbol{\zeta}_c^n}_c^2 +\int_{\Omega^d}g_p^n(\zeta_c^n)^2  dx \nonumber\\
&\leq\sum_{K\in \mathcal{T}}\int_K \phi\, d_t \xi_c^n \zeta_c^n \dif  x+\sum_{K\in \mathcal{T}}\int_K \phi\,(\partial_t c^n-d_tc^n) \zeta_c^n \dif  x +B_h^{tr}(u_h^n;\boldsymbol{\xi}_c^n, \boldsymbol{\zeta}_c^n)\nonumber\\
  &\quad+B_h^{a}(u^n-u_h^n;  \boldsymbol{c}^n, \boldsymbol{\zeta}_c^n)+\Big[B_h^d(u^n; \boldsymbol{c}^n, \boldsymbol{\zeta}_c^n)-B_h^d(u_h^n; \boldsymbol{c}^n, \boldsymbol{\zeta}_c^n)\Big]\nonumber\\
  &\quad-\frac{1}{2} \sum_{K\in \mathcal{T}^d} \int_K \nabla \cdot u_h^n (\zeta_c^n)^2\,\dif x + \int_{\Omega^d}g_p^n\xi_c^n\zeta_c^n \dif x\nonumber\\
  &:= I_1+\ldots+I_7.
\end{align*}
Using \cref{eq:assumption_phi}, the Cauchy-Schwarz inequality, Young’s inequality with constant $\gamma$, and  \cref{eq:est_interpolant_brenner},
\begin{align*}
    I_1
    &\leq \phi^*\Big\|\frac{1}{\Delta t}\int_{t^{n-1}}^{t^n}\partial_t\xi_c\,{\rm d}t\Big\|_{\Omega}\|\zeta_c^n \|_{\Omega}
    \leq \phi^*\frac{1}{\sqrt{\Delta t}}\|\partial_t\xi_c\|_{L^2(t^{n-1},t^{n},L^2(\Omega))}\|\zeta_c^n \|_{\Omega}\\
    &\le  C\frac{(\phi^*)^2 h^{2k_c}}{\phi_*\Delta t}\|\partial_t c\|_{L^2(t^{n-1},t^n,H^{k_c}(\Omega))}^2+\gamma(\phi_*\|\zeta_c^n \|_{\Omega}^2).\nonumber
\end{align*}
Again by \cref{eq:assumption_phi}, this time using Taylor's theorem in integral form, and applying Young's inequality,
\begin{align*}
    I_2
    =\sum_{K\in \mathcal{T}}\int_K \phi\,(\partial_t c^n-d_tc^n) \zeta_c^n \dif  x
    &\le \phi^* \Big\|\frac{1}{\Delta t}\int_{t^{n-1}}^{t^n}(t-t^{n-1})\partial_{tt} c\,{\rm d}t\Big\|_{\Omega}\|\zeta_c^n\|_{\Omega}
    \\
    &\le C\frac{(\phi^*)^2\Delta t}{\phi_*}\|\partial_{tt} c\|_{L^2(t^{n-1},t^n;L^2(\Omega))}^2+\gamma(\phi_*\|\zeta_c^n\|_{\Omega}^2).
\end{align*}
The following series of inequalities is dedicated to finding an upper bound for $I_3$. By definition of $B_h^{tr}$,
\begin{equation}\label{eq:I3c}
    I_3=B_h^{tr}(u_h^n;\boldsymbol{\xi}_c^n, \boldsymbol{\zeta}_c^n)= B_h^a(u_h^n;\boldsymbol{\xi}_c^n,\boldsymbol{\zeta}_c^n) +  B_h^d(u_h^n;\boldsymbol{\xi}_c^n,\boldsymbol{\zeta}_c^n).
\end{equation}
We will bound $B_h^a$ and $B_h^d$ separately, starting with $B_h^a$. Noting that $\xi_c^n-\bar{\xi}_c^n$ vanishes on facets, we have by \cref{eq:bilinearform_Bha},
\begin{multline}\label{ineq:I31c}
    B_h^a(u_h^n; \boldsymbol{\xi}_c^n, \boldsymbol{\zeta}_c^n)
    =- \sum_{K\in \mathcal{T}} \int_K \xi_c^n \, u_h^n \cdot\nabla \zeta_c^n\,\dif x
    + \sum_{K\in \mathcal{T}} \int_{\partial K} \xi_c^n\, (u_h^n-u^n) \cdot n  (\zeta_c^n-\bar{\zeta}_c^n)\,\dif s\\
    +\sum_{K\in \mathcal{T}} \int_{\partial K} \xi_c^n\, u^n \cdot n  (\zeta_c^n-\bar{\zeta}_c^n)\,\dif s =:  I_{311} + I_{312}+I_{313}.
\end{multline}
The term $I_{311}$ can be bounded by H\"{o}lder's inequality, and \cref{eq:estimate_uh_infty,eq:est_interpolant_brenner}:
\begin{equation}
\label{ineq:I311c}
      I_{311}\leq C \|u_h^n\|_{0,\infty, \Omega} \|\xi_c^n\|_{\Omega} 
      \|\nabla \zeta_c^n\|_{\Omega}
      \leq Ch^{k_c} \|c^n\|_{k_c,\Omega} 
      \tnorm{\boldsymbol{\zeta}_c^n}_c.
\end{equation}
By H\"{o}lder's inequality and using \cref{ineq:trace_cont_inf}, 
\begin{align}\label{ineq:I312c}
    I_{312}
    &\leq
    \sum_{K\in \mathcal{T}}\|\xi_c^n\|_{0,\infty,\partial K}\|u_h^n-u^n\|_{\partial K}\|\zeta_c^n-\bar{\zeta}_c^n\|_{\partial K}\nonumber\\
& \leq C\Big(\sum_{K\in \mathcal{T}}h_K\|u_h^n-u^n\|_{\partial K}^2\Big)^{1/2} \tnorm{\boldsymbol{\zeta}_c^n}_c.
\end{align}
Using H\"{o}lder's inequality and this time employing \cref{ineq:trace_cont_inf,ineq:trace_cont,eq:est_interpolant_brenner},
\begin{align}\label{ineq:I313c}
    I_{313}&=\sum_{K\in \mathcal{T}}\|\xi_c^n\|_{\partial K}\|u^n\|_{0,\infty,\partial K}\|\zeta_c^n-\bar{\zeta}_c^n\|_{\partial K}\nonumber\\
    &\leq C  \Big(\sum_{K\in \mathcal{T}}(\|\xi_c^n\|_{K}^2+h_K^2\|\xi_c^n\|_{1,K}^2)\Big)^{1/2}\Big(\sum_{K\in \mathcal{T}}h_K^{-1}\|\zeta_c^n-\bar{\zeta}_c^n\|_{\partial K}^2\Big)^{1/2}\\   
    & \leq C h^{k_c}\|c^n\|_{k_c,\Omega}\tnorm{\boldsymbol{\zeta}_c^n}_c. \nonumber
\end{align}
Putting \crefrange{ineq:I31c}{ineq:I313c} together and using Young's inequality, we find
\begin{equation}
    \label{ineq:Bhan}
   B_h^a(u_h^n; \boldsymbol{\xi}_c^n, \boldsymbol{w}_h)
   \leq C\Big(\sum_{K\in \mathcal{T}}h_K\|u_h^n-u^n\|^2_{\partial K}+h^{2k_c}\|c^n\|^2_{k_c,\Omega}\Big)+\epsilon\tnorm{\boldsymbol{\zeta}_c^n}_c^2
\end{equation}
We now bound $B_h^d$ in \cref{eq:I3c}. 
Since $\xi_c^n-\bar{\xi}_c^n=0$ on $\partial K$,
\begin{align}\label{eq:I32c}
      B_h^d(u_h^n;\boldsymbol{\xi}_c^n,\boldsymbol{\zeta}_c^n)=&  \sum_{K\in \mathcal{T}}\int_{ K} \widetilde{D}(u_h^n)\nabla \xi_c^n \cdot \nabla \zeta_c^n\,\dif x
    - \sum_{K\in \mathcal{T}}\int_{\partial K} [\widetilde{D}(u_h^n)\nabla \xi_c^n] \cdot n \,(\zeta_c^n-\bar{\zeta}_c^n)\,\dif s=:I_{321}+I_{322}.
\end{align}
By H\"{o}lder's inequality, \cref{ineq:Dtildemax}, and \cref{eq:est_interpolant_brenner},
\begin{equation}\label{ineq:I321c}
     I_{321} \leq  \widetilde{D}_{\max}\|\nabla \xi_c^n\|_{\Omega}\|\nabla \zeta_c^n\|_{\Omega}\leq Ch^{k_c} \|c\|_{k_c+1,\Omega}\tnorm{\boldsymbol{\zeta}_c^n}_c.
\end{equation}
Again by H\"{o}lder's inequality and this time using \cref{ineq:Dtildemax,ineq:trace_cont,eq:est_interpolant_brenner},
\begin{align}\label{ineq:I322c}
|I_{322}|&\leq  \widetilde{D}_{\max}\Big(\sum_{K\in \mathcal{T}} h_K\|\nabla \xi_c^n\|_{\partial K}^2\Big)^{1/2}\nonumber
      \Big(\sum_{K\in \mathcal{T}}h_K^{-1}\|\zeta_c^n-\bar{\zeta}_c^n\|_{\partial K}^2\Big)^{1/2}\nonumber\\
      &\leq   C\Big(\sum_{K\in \mathcal{T}} (\|\nabla \xi_c^n\|_{K}^2+h_K^2\|\nabla \xi_c^n\|_{1,K}^2)\Big)^{1/2}\tnorm{\boldsymbol{\zeta}_c^n}_{c} 
    \leq  C  h^{k_c} \|c^n\|_{k_c+1,\Omega} \tnorm{\boldsymbol{\zeta}_c^n}_{c}.
\end{align} 
Hence, the combination of \cref{eq:I32c,ineq:I321c,ineq:I322c} and using Young's inequality results in:
\begin{align}\label{ineq:Bhdn}
    B_h^d(u_h^n;\boldsymbol{\xi}_c^n,\boldsymbol{\zeta}_c^n)\le C h^{k_c} \|c^n\|_{k_c+1,\Omega} \tnorm{\boldsymbol{\zeta}_c^n}_{c}
    \le Ch^{2k_c} \|c^n\|_{k_c+1,\Omega}^2+\epsilon\tnorm{\boldsymbol{\zeta}_c^n}_{c}^2.
\end{align}
Therefore, from \cref{ineq:Bhan,ineq:Bhdn},
\begin{equation*}
	I_3\leq Ch^{2k_c}\|c^n\|_{k_c+1,\Omega}^2+\sum_{K\in \mathcal{T}}h_K\|u_h^n-u^n\|^2_{\partial K}+2\epsilon\tnorm{\boldsymbol{\zeta}_c^n}_{c}^2.
\end{equation*}
Since $c = \gamma(c)$ on $\partial K$,
\begin{align}\label{ineq:I4c}
   I_4
   &=  - \sum_{K\in \mathcal{T}} \int_K c^n \, (u^n-u_h^n) \cdot\nabla \zeta_c^n\,\dif x
  + \sum_{K\in \mathcal{T}} \int_{\partial K} c^n\, (u^n-u_h^n) \cdot n  (\zeta_c^n-\bar{\zeta}_c^n)\,\dif s=: I_{41} + I_{42}. 
\end{align}  
H\"{o}lder's and Young's inequalities give
\begin{align}\label{ineq:I41c}
    I_{41}
    &\leq \|u^n-u_h^n\|_{\Omega}\|c^n\|_{0,\infty, \Omega}\tnorm{\boldsymbol{\zeta}_c^n}_c\leq C\|c^n\|_{0,\infty, \Omega}^2\|u^n-u_h^n\|_{\Omega}^2+\epsilon\tnorm{\boldsymbol{\zeta}_c^n}_c^2,
\end{align}
and
\begin{align}\label{ineq:I42c}
      I_{42}
      &\leq C \|c^n\|_{0,\infty,\Omega}
      \Big(\sum_{K\in \mathcal{T}_h}h_K \|u^n-u_h^n\|^2_{\partial K}\Big)^{1/2}\tnorm{\boldsymbol{\zeta}_c^n}_c
      \leq C \|c^n\|_{0,\infty,\Omega}^2
     \sum_{K\in \mathcal{T}_h}h_K \|u^n-u_h^n\|^2_{\partial K}+\epsilon \tnorm{\boldsymbol{\zeta}_c^n}_c^2.
\end{align}    
Collecting \cref{ineq:I4c,ineq:I41c,ineq:I42c} leads to
\begin{align}
    I_4\leq 2\epsilon\tnorm{\boldsymbol{\zeta}_c^n}_c^2 +C\|c^n\|_{0,\infty,\Omega}^2\Big(\|u^n-u_h^n\|_{\Omega}^2+ 
     \sum_{K\in \mathcal{T}_h}h_K \|u^n-u_h^n\|^2_{\partial K}\Big).\nonumber
\end{align}
Since $c = \gamma(c)$ on element boundaries and $\widetilde{D}(u^n)-\widetilde{D}(u_h^n)=0$ in $\Omega^s$,
\begin{align} \label{ineq:I5csplit}
 I_5
    =&\sum_{K\in \mathcal{T}^d}\int_{ K} [D(u^n)-D(u_h^n)]\nabla c^n \cdot \nabla \zeta_c^n\,\dif x
    - \sum_{K\in \mathcal{T}^d}\int_{\partial K} [(D(u^n)-D(u_h^n))\nabla c^n] \cdot n \,(\zeta_c^n-\bar{\zeta}_c^n)\,\dif s \nonumber\\
    =:& I_{51}+I_{52}.    
\end{align}
Using the Lipschitz property of $D$ \cref{eq:DLipschitz}, H\"{o}lder's and Young's inequalities,
\begin{align}\label{ineq:I51c}
     I_{51} &\leq C\|u^n-u_h^n\|_{\Omega^d} \|\nabla c^n\|_{0,\infty, \Omega^d}\|\nabla \zeta_c^n\|_{\Omega^d}\leq C\|\nabla c^n\|_{0,\infty, \Omega^d}^2\|u^n-u_h^n\|_{\Omega^d}^2+\epsilon\tnorm{\boldsymbol{\zeta}_c^n}_c^2.
\end{align}
%
Similarly,
\begin{align}\label{ineq:I52}
    I_{52}
    &\leq \sum_{K\in \mathcal{T}^d}C\|u^n-u_h^n\|_{\partial K}\|\nabla c^n\|_{0,\infty,K}\|\zeta_c^n-\bar{\zeta}_c^n\|_{\partial K}
    \leq C \|c^n\|_{1,\infty,\Omega^d}\Big(\sum_{K\in \mathcal{T}^d}h_K\norm[0]{u^n-u_h^n}^2_{\partial K}\Big)^{1/2}\tnorm{\boldsymbol{\zeta}_c^n}_c\nonumber\\
    &\leq C\Big(\|c^n\|^2_{1,\infty, \Omega^d}\sum_{K\in \mathcal{T}^d}h_K\norm[0]{u^n-u_h^n}^2_{\partial K}\Big)+\epsilon\tnorm{\boldsymbol{\zeta}_c^n}^2_c.
\end{align}
Therefore, substituting \cref{ineq:I51c,ineq:I52} in \cref{ineq:I5csplit}  yields
\begin{align*}\label{ineq:I5c}
    I_5 \leq 2\epsilon\tnorm{\boldsymbol{\zeta}_c^n}_{c}^2+C\|c^n\|^2_{1,\infty, \Omega^d}\Big(\sum_{K\in \mathcal{T}^d}h_K\norm[0]{u^n-u_h^n}^2_{\partial K}+\|u^n-u_h^n\|_{\Omega^d}^2\Big).
\end{align*}
By \cref{eq:masscons-1},  the stability of the $L^2$-projection $\Pi_Q$ and H\"{o}lder's inequality,
\begin{align*}
    I_6
    =-\frac{1}{2} \sum_{K\in \mathcal{T}^d} \int_K \Pi_Q(g_i^n-g_p^n) (\zeta_c^n)^2\,\dif x
    \leq \frac1{2\phi_*} \|g_i^n-g_p^n\|_{0,\infty,\Omega^d}(\phi_*\| \zeta_c^n\|_{\Omega}^2).
\end{align*}
Finally, using H\"{o}lder's inequality, \cref{eq:est_interpolant_brenner}, and Young's inequality,
\begin{align*}
    I_7 & \leq \|g_p^n\|_{0,\infty,\Omega^d}\|\xi_c^n\|_{\Omega^d}\|\zeta_c^n\|_{\Omega^d} \leq Ch^{k_c}\|g_p^n\|_{0,\infty,\Omega^d}\|c^n\|_{k_c,\Omega^d}\|\zeta_c^n\|_{\Omega^d}
    \nonumber\\
    &\leq Ch^{2k_c}\phi_*^{-1}\|g_p^n\|^2_{0,\infty,\Omega^d}\|c^n\|_{k_c,\Omega^d}^2+\gamma(\phi_*\|\zeta_c^n\|_{\Omega}^2).
\end{align*}
Collecting all bounds, choosing $\epsilon=C_{tr}/12$ ($C_{tr}$ is the coercivity constant), $\gamma=1/6$, and recalling that $g_p^n\geq 0$, we find:
\begin{align*}
   \frac{\phi_*}{2\Delta t} &(\|\zeta_c^n\|_{\Omega}^2-\|\zeta_c^{n-1}\|_{\Omega}^2)+\dfrac{C_{tr}}2\tnorm{\boldsymbol{\zeta}_c^n}_c^2  \nonumber\\
    \leq& C\Big(\dfrac{(\phi^*)^2\phi_*^{-1} h^{2k_c}}{\Delta t}\|\partial_t c\|_{L^2(t^{n-1},t^n;H^{k_c}(\Omega))}^2+ (\phi^*)^2\phi_*^{-1}\Delta t\|\partial_{tt} c\|_{L^2(t^{n-1},t^n;L^2(\Omega))}^2\\
     &+ h^{2k_c}\|c^n\|_{k_c+1,\Omega}^2+\|c^n\|_{1,\infty,\Omega}^2\|u^n-u_h^n\|_{\Omega}^2+(\|c^n\|_{1,\infty,\Omega}^2+1) \sum_{K\in \mathcal{T}}h_K \|u^n-u_h^n\|^2_{\partial K}\Big)\nonumber\\
        &+h^{2k_c}\phi_*^{-1}\|g_p^n\|^2_{0,\infty,\Omega^d}\|c^n\|_{k_c,\Omega^d}^2\Big)+\frac12\Big(\phi_*^{-1}\|g_i^n-g_p^n\|_{0,\infty,\Omega^d}+1\Big)(\phi_*\| \zeta_c^n\|_{\Omega}^2).\nonumber
\end{align*}
Multiplying by $2\Delta t$, summing over $m$, and using Corollary~\ref{cor:velerror},
\begin{align*}
    \phi_*\|\zeta_c^n\|_{\Omega}^2&+C_{tr}\Delta t\sum\limits_{m=1}^n\tnorm{\boldsymbol{\zeta}_c^m}_c^2  \nonumber\\
    &\leq \phi_*\|\zeta_c^{0}\|_{\Omega}^2+ C\Big[h^{2k_c}\big(\dfrac{(\phi^*)^2}{\phi_*} \|\partial_t c\|_{L^2(0,T;H^{k_c}(\Omega))}^2+\|c\|_{\ell^2(0,T;H^{k_c+1}(\Omega))}^2\big)\\
    &\quad+ (\phi^*)^{2}\phi_*^{-1}(\Delta t)^2\|\partial_{tt} c\|_{L^2(0,T;L^2(\Omega))}^2+\Delta t\sum\limits_{m=1}^n\Big(\phi_*^{-1}\|g_i^m-g_p^m\|_{0,\infty,\Omega^d}+1\Big)(\phi_*\| \zeta_c^m\|_{\Omega}^2)\nonumber\\
    &\quad+\Delta t\sum\limits_{m=1}^n(\|c^m\|_{1,\infty,\Omega}^2+1)\big((\Delta t)^2+h^{2k_f}+h^{2k_c}\big)+h^{2k_c}\Delta t\sum\limits_{m=1}^n\phi_*^{-1}\|g_p^m\|^2_{0,\infty,\Omega^d}\|c^m\|_{k_c,\Omega^d}^2\Big].\nonumber
\end{align*}
Using \cite[Lemma 1.58]{DiPietro:book} and \cref{eq:est_interpolant_brenner},
\begin{align*}
    \|\zeta_c^0\|_{\Omega} =\|c_h^0-\mathcal{I} c_0\|_{\Omega}
    &\leq \|\Pi_Cc_0-c_0\|_{\Omega} + \|c_0-\mathcal{I} c_0\|_{\Omega}\leq Ch^{k_c}\|c_0\|_{k_c,\Omega}.   
\end{align*}
Therefore, the result follows by Gr\"{o}nwall's inequality \cite[Lemma 27]{Layton:book} assuming that $\Delta t$ is sufficiently small.
\end{proof}
By the triangle inequality and \cref{eq:est_interpolant_brenner}, we immediately have
\begin{equation}
\label{eq:L2-c}
 \|c_h^n-c^n\|_{\Omega}\leq C(\Delta t+h^{k_f}+h^{k_c}), \quad \forall n\geq 1.
\end{equation}
%
\section{Numerical examples}
\label{sec:numerical_examples}
\Cref{alg:uhfirst} is implemented in the higher-order finite element library Netgen/NGSolve \cite{Schoberl:1997,Schoberl:2014}.
In all numerical examples we choose $\bar{\Omega}=[0,1]^2$ with subregions $\bar{\Omega}^d=[0,1]\times [0,0.5]$ and $\bar{\Omega}^s=[0,1]\times [0.5,1]$. We furthermore choose the penalty parameters as $\beta_s=6k_f^2$ and $\beta_{tr}=6k_c^2$ \cite[Lemma 1, Section 5]{Ainsworth:2018}.
\subsection{Example 1} 
\label{sec:example1}
We first consider the constant coefficient case, i.e., the time-dependent one-way coupled problem in which the numerical solution to the Stokes/Darcy problem is unaffected by the concentration. Let $\alpha = \tfrac{1}{2}(1+4\pi^2)\sqrt{\kappa}$, 
$\widetilde{D} = D =\left[\begin{array}{cc} 0.01&0.005\\ 0.005 & 0.02\end{array}\right]$ on $\Omega$, and $T=0.1$.
The source terms and boundary conditions for the Stokes/Darcy--transport problem are chosen such that the exact solution is given by
\begin{subequations}\label{eq:exact_solution_SDT}
\begin{align}
    u^s &= (-\dfrac1{2\pi^2}\sin(\pi x_1+t)e^{(x_2+t)/2}, \dfrac1{\pi}\cos(\pi x_1+t)e^{(x_2+t)/2})^T\\
    u^d &= (-2\sin(\pi x_1+t)e^{(x_2+t)/2}, \dfrac1{\pi}\cos(\pi x_1+t)e^{(x_2+t)/2})^T,\\
    p^s&=\dfrac{\kappa\mu-2}{\kappa\pi}\cos(\pi x_1+t)e^{(x_2+t)/2},\\
    p^d&=-\dfrac2{\kappa\pi}\cos(\pi x_1+t)e^{(x_2+t)/2},\\
    c&=\sin(2\pi(x_1-t))\cos(2\pi(x_2-t)).
\end{align}
\end{subequations}
Note that this solution satisfies all the interface conditions and that $\nabla \cdot u^s=0$ in $\Omega^s$. 

We present our numerical results for a wide range of values for $\kappa$ and $\mu$: $\kappa=\mu=1$; $\kappa=10^3$, $\mu=10^{-6}$; $\kappa=1$, $\mu=10^{-6}$; and $\kappa=10^{-3}$, $\mu=10^{-6}$. Since we are primarily interested in the spatial error, to minimize the temporal error as much as possible, we use the third order backward differentiation formulae (BDF3) as time stepping method even though the sequential algorithm \ref{alg:uhfirst} is only first order accurate in time. We choose $\Delta t=0.1h^{k_f}/(k_f+1)$ and present errors and rates of convergence using $k_f=2$, $k_c=1$ in \Cref{tab:E1_Stokes-k2,tab:E1_Darcy-k2,tab:E1_concentration-k2} and using $k_f=3, k_c=2$ in \Cref{tab:E1_Stokes-k3,tab:E1_Darcy-k3,tab:E1_concentration-k3}.

\begin{table}[!ht]
  \caption{Errors and rates of convergence at final time $T=0.1$ for $u_h$ and $p_h$ in the Stokes region $\Omega^s$ for the test case in \Cref{sec:example1} using $k_f=2$, $k_c=k_f-1$, and BDF3 time stepping with $\Delta t= 0.1h^2/3$.}  
  \label{tab:E1_Stokes-k2}
 \centering
   \begin{tabular}{ccccccc}
  \textit{h} & dofs   &  $\|u_h-u\|_{\Omega^s}$ &    rate   &   $\|p_h-p\|_{\Omega^s}$  &   rate    &  $\|\nabla\cdot u_h\|_{\Omega^s}$ \\ \hline
    \multicolumn{7}{l}{$\kappa=1, \mu=1$} \\
    \hline
        1/4 &  745 & 2.6e-04 & -- & 1.2e-02 & -- & 1.4e-16 \\
        1/8 &  3811 & 2.0e-05 & 3.7 & 1.9e-03 & 2.6 & 1.8e-16\\
        1/16 & 14167 & 2.2e-06 & 3.2 & 4.5e-04 & 2.1 & 1.6e-16\\
        1/32 & 57181 & 2.7e-07 & 3.1 & 1.1e-04 & 2.1 & 1.7e-16\\
    \hline
    \multicolumn{7}{l}{$\kappa=10^3$, $\mu=10^{-6}$} \\
    \hline
        1/4 & 745 & 2.5e-04 & -- & 1.3e-05 & -- & 8.5e-17\\
        1/8 & 3811 & 2.0e-05 & 3.6 & 1.9e-06 & 2.8 & 9.7e-17\\
        1/16& 14167 & 2.2e-06 & 3.2 & 4.5e-07 & 2.1 & 9.4e-17\\
        1/32& 57181 & 2.6e-07 & 3.1 & 1.1e-07 & 2.1 & 9.0e-17\\
    \hline
    \multicolumn{7}{l}{$\kappa=1, \mu=10^{-6}$} \\
    \hline
        1/4 & 745 & 6.4e-04 & -- & 2.7e-02 & -- & 8.6e-16  \\
        1/8 & 3811 & 4.8e-05 & 3.7 & 5.3e-03 & 2.3 & 1.8e-15    \\
        1/16 & 14167 & 4.4e-06 & 3.4 & 1.2e-03 & 2.1 & 3.0e-15  \\
        1/32 &  57181 & 5.2e-07 & 3.1 & 3.0e-04 & 2.0 & 6.0e-15\\
    \hline 
    \multicolumn{7}{l}{$\kappa=10^{-3}, \mu=10^{-6}$} \\
     \hline
        1/4 & 745 &  7.3e-04 & -- & 1.2e+01 & -- & 1.2e-13 \\
        1/8 &  3811 & 5.0e-05 & 3.9 & 1.9e+00 & 2.6 & 1.6e-13   \\
        1/16 & 14167 & 4.1e-06 & 3.6 & 4.5e-01 & 2.1 & 1.4e-13  \\
        1/32 &  57181 & 3.8e-07 & 3.4 & 1.1e-01 & 2.1 & 1.5e-13 \\
  \end{tabular}
\end{table}
\begin{table}[!ht]
  \caption{Errors and rates of convergence at final time $T=0.1$ for $u_h$ and $p_h$ in the Darcy region $\Omega^d$ for the test case in \Cref{sec:example1} using $k_f=2$, $k_c=k_f-1$, and BDF3 time stepping with $\Delta t= 0.1h^2/3$.}  
  \label{tab:E1_Darcy-k2}
  \centering
    \begin{tabular}{ccccccc}
 \textit{h} & dofs   &  $\|u_h-u\|_{\Omega^d}$ &    rate   &   $\|p_h-p\|_{\Omega^d}$  &   rate    &  $\|\Pi_Q(\nabla\cdot (u_h -u))\|_{\Omega^d}$ \\ \hline
    \multicolumn{7}{l}{$\kappa=1, \mu=1$} \\
    \hline
        1/4 & 745 & 3.1e-03 & -- & 9.1e-03 & -- & 5.9e-09 \\
        1/8 &  3811 & 1.9e-04 & 4.0 & 1.4e-03 & 2.7 & 9.1e-11\\
        1/16 & 14167 & 2.5e-05 & 2.9 & 3.7e-04 & 1.9 & 1.8e-12\\
        1/32 & 57181 & 2.7e-06 & 3.2 & 8.4e-05 & 2.1 & 3.8e-12\\
    \hline
    \multicolumn{7}{l}{$\kappa=10^3$, $\mu=10^{-6}$} \\
    \hline
        1/4 & 745 &  3.1e-03 & -- & 9.1e-06 & -- & 5.9e-09\\
        1/8 &  3811 & 1.9e-04 & 4.0 & 1.4e-06 & 2.7 & 9.1e-11\\
        1/16& 14167 & 2.5e-05 & 2.9 & 3.7e-07 & 1.9 & 1.7e-12\\
        1/32& 57181 & 2.7e-06 & 3.2 & 8.4e-08 & 2.1 & 3.8e-12\\
    \hline
    \multicolumn{7}{l}{$\kappa=1, \mu=10^{-6}$} \\
    \hline
        1/4 & 745 &  3.1e-03 & -- & 9.1e-03 & -- & 5.9e-09  \\
        1/8 & 3811 &1.9e-04 & 4.0 & 1.4e-03 & 2.7 & 9.1e-11   \\
        1/16 & 14167 & 2.5e-05 & 2.9 & 3.7e-04 & 1.9 & 1.8e-12  \\
        1/32 & 57181 &   2.7e-06 & 3.2 & 8.4e-05 & 2.1 & 3.8e-12 \\
    \hline 
    \multicolumn{7}{l}{$\kappa=10^{-3}, \mu=10^{-6}$} \\
     \hline
        1/4 & 745 &  3.1e-03 & -- & 9.1e+00 & -- & 5.9e-09 \\
        1/8 & 3811 & 1.9e-04 & 4.0 & 1.4e+00 & 2.7 & 9.1e-11   \\
        1/16 & 14167 & 2.5e-05 & 2.9 & 3.7e-01 & 1.9 & 1.8e-12 \\
        1/32 & 57181 &  2.6e-06 & 3.2 & 8.4e-02 & 2.1 & 4.7e-12\\
  \hline 
  \end{tabular}
\end{table}
\begin{table}[!ht]
  \caption{Errors and rates of convergence at final time $T=0.1$ for $c_h$ in $\Omega$, on a mesh with $h=1/4, 1/8, 1/16, 1/32$, for the test case in \Cref{sec:example1} using $k_f=2$, $k_c=k_f-1$, and BDF3 time stepping with $\Delta t= 0.1h^2/3$.} 
      \centering
  \label{tab:E1_concentration-k2}
    \begin{tabular}{c|cc|cc|cc|cc}
\multicolumn{1}{c}{}&\multicolumn{2}{c|}{$\kappa=1$, $\mu=1$} &\multicolumn{2}{c|}{$\kappa=10^3$, $\mu=10^{-6}$} &\multicolumn{2}{c|}{$\kappa=1$, $\mu=10^{-6}$} &\multicolumn{2}{c}{$\kappa=10^{-3}$, $\mu=10^{-6}$}   \\ \hline
        dofs & $||c-c_h||_{\Omega}$ & rate & $||c-c_h||_{\Omega}$& rate & $||c-c_h||_{\Omega}$ & rate & $||c-c_h||_{\Omega}$ & rate\\
    \hline 
        184  & 9.7e-02 & -- & 9.7e-02 & -- & 9.7e-02 & -- & 9.8e-02 & --\\
        944  & 2.2e-02 & 2.1 & 2.2e-02 & 2.1 & 2.2e-02 & 2.1 & 2.2e-02 & 2.1\\
        3520  & 5.4e-03 & 2.0 & 5.4e-03 & 2.0 & 5.4e-03 & 2.0 & 5.0e-03 & 2.2\\
        14216  & 1.1e-03  & 2.3 & 1.1e-03 & 2.3 & 1.1e-03 & 2.3 & 1.1e-03 & 2.2\\
  \end{tabular}
\end{table}
\begin{table}[!ht]
  \caption{Errors and rates of convergence at final time $T=0.1$ for $u_h$ and $p_h$ in the Stokes region $\Omega^s$ for the test case in \Cref{sec:example1} using $k_f=3$, $k_c=k_f-1$, and BDF3 time stepping with $\Delta t= 0.1h^3/4$.}  
  \label{tab:E1_Stokes-k3}
    \centering
   \begin{tabular}{ccccccc}
 \textit{h} & dofs   &  $\|u_h-u\|_{\Omega^s}$ &    rate   &   $\|p_h-p\|_{\Omega^s}$  &   rate    &  $\|\nabla\cdot u_h\|_{\Omega^s}$ \\ \hline
    \multicolumn{7}{l}{$\kappa=1, \mu=1$} \\
    \hline
        1/4 &  1161 & 5.6e-05 & -- & 4.6e-03 & -- & 2.0e-15 \\
        1/8 &  5993 & 1.4e-06 & 5.4 & 2.2e-04 & 4.4 & 3.7e-15 \\
        1/16 & 22081 & 7.1e-08 & 4.3 & 2.1e-05 & 3.3 & 4.5e-15\\
        1/32 & 90241 & 3.7e-09 & 4.3 & 2.3e-06 & 3.2 & 6.3e-15\\
    \hline
    \multicolumn{7}{l}{$\kappa=10^3$, $\mu=10^{-6}$} \\
    \hline
        1/4 & 1161 & 2.2e-05 & -- & 7.3e-07 & -- & 1.5e-16\\
        1/8 & 5993 & 6.0e-07 & 5.2 & 5.7e-08 & 3.7 & 1.3e-16 \\
        1/16& 22081 & 3.3e-08 & 4.2 & 6.4e-09 & 3.2 & 1.2e-16 \\
        1/32& 90241& 1.8e-09 & 4.2 & 7.7e-10 & 3.1 & 1.2e-16\\
    \hline
    \multicolumn{7}{l}{$\kappa=1, \mu=10^{-6}$} \\
    \hline
        1/4 & 1161 & 2.2e-05 & --  & 6.9e-04 & -- & 1.5e-16  \\
        1/8 & 5993 & 6.1e-07 & 5.2 & 5.6e-05 & 3.6 & 1.5e-16   \\
        1/16 & 22081 & 3.3e-08 & 4.2 & 6.4e-06 & 3.1 & 1.2e-16  \\
        1/32 &  90241&  1.8e-09 & 4.2 & 7.7e-07 & 3.1 & 1.2e-16 \\
    \hline 
    \multicolumn{7}{l}{$\kappa=10^{-3}, \mu=10^{-6}$} \\
     \hline
        1/4 & 1161 & 2.6e-05 & -- & 6.9e-01 & -- & 5.5e-14 \\
        1/8 &  5993 & 9.4e-07 & 4.8 & 5.4e-02 & 3.7 & 4.2e-14    \\
        1/16 & 22081 & 4.5e-08 & 4.4 & 6.3e-03 & 3.1 & 1.9e-14  \\
        1/32 &  90241& 2.2e-09 & 4.3 & 7.7e-04 & 3.0 & 1.0e-14\\
  \end{tabular}
\end{table}
\begin{table}[!ht]
  \caption{Errors and rates of convergence at final time $T=0.1$ for $u_h$ and $p_h$ in the Darcy region $\Omega^d$ for the test case in \Cref{sec:example1} using $k_f=3$, $k_c=k_f-1$, and BDF3 time stepping with $\Delta t= 0.1h^3/4$.}  
  \label{tab:E1_Darcy-k3}
    \centering
    \begin{tabular}{ccccccc}
 \textit{h} & dofs   &  $\|u_h-u\|_{\Omega^d}$ &    rate   &   $\|p_h-p\|_{\Omega^d}$  &   rate    &  $\|\Pi_Q(\nabla\cdot (u_h -u))\|_{\Omega^d}$ \\ \hline
    \multicolumn{7}{l}{$\kappa=1, \mu=1$} \\
    \hline
        1/4 & 1161 & 1.2e-04 & -- & 5.4e-04 & -- & 4.1e-12  \\
        1/8 &  5993 & 3.6e-06 & 5.1 & 4.0e-05 & 3.8 & 8.3e-13 \\
        1/16 & 22081 & 2.2e-07 & 4.0 & 5.1e-06 & 3.0 & 2.9e-12 \\
        1/32 & 90241 & 1.3e-08 & 4.1 & 6.1e-07 & 3.1 & 1.2e-11 \\
    \hline
    \multicolumn{7}{l}{$\kappa=10^3$, $\mu=10^{-6}$} \\
    \hline
        1/4 & 1161 & 1.2e-04 & -- & 5.4e-07 & -- & 4.1e-12 \\
        1/8 &  5993 & 3.6e-06 & 5.1 & 4.0e-08 & 3.8 & 8.2e-13\\
        1/16& 22081 & 2.2e-07 & 4.0 & 5.1e-09 & 3.0 & 3.1e-12\\
        1/32& 90241 & 1.3e-08 & 4.1 & 6.1e-10 & 3.1 & 1.2e-11\\
    \hline
    \multicolumn{7}{l}{$\kappa=1, \mu=10^{-6}$} \\
    \hline
        1/4 & 1161 & 1.2e-04 & -- & 5.4e-04 & -- & 4.1e-12  \\
        1/8 & 5993 & 3.5e-06 & 5.1 & 4.0e-05 & 3.8 & 7.2e-13  \\
        1/16 & 22081 & 2.2e-07 & 4.0 & 5.1e-06 & 3.0 & 3.0e-12  \\
        1/32 & 90241 & 1.3e-08 & 4.1 & 6.1e-07 & 3.1 & 1.2e-11 \\
    \hline 
    \multicolumn{7}{l}{$\kappa=10^{-3}, \mu=10^{-6}$} \\
     \hline
        1/4 & 1161 & 1.2e-04 & -- & 5.4e-01 & -- & 4.1e-12  \\
        1/8 & 5993 & 3.5e-06 & 5.1 & 4.0e-02 & 3.8 & 8.9e-13    \\
        1/16 & 22081 & 2.2e-07 & 4.0 & 5.1e-03 & 3.0 & 3.3e-12 \\
        1/32 & 90241 & 1.3e-08 & 4.1 & 6.1e-04 & 3.1 & 1.3e-11 \\
  \end{tabular}
\end{table}
\begin{table}[!ht]
  \caption{Errors and rates of convergence at final time $T=0.1$ for $c_h$ in $\Omega$, on a mesh with $h=1/4, 1/8, 1/16, 1/32$, for the test case in \Cref{sec:example1} using $k_f=3$, $k_c=k_f-1$, and BDF3 time stepping with $\Delta t= 0.1h^3/4$.}  
  \label{tab:E1_concentration-k3}
  \centering
    \begin{tabular}{c|cc|cc|cc|cc}
\multicolumn{1}{c}{}&\multicolumn{2}{c|}{$\kappa=1$, $\mu=1$} &\multicolumn{2}{c|}{$\kappa=10^3$, $\mu=10^{-6}$} &\multicolumn{2}{c|}{$\kappa=1$, $\mu=10^{-6}$} &\multicolumn{2}{c}{$\kappa=10^{-3}$, $\mu=10^{-6}$}   \\ \hline
        dofs & $||c-c_h||_{\Omega}$ & rate & $||c-c_h||_{\Omega}$& rate & $||c-c_h||_{\Omega}$ & rate & $||c-c_h||_{\Omega}$ & rate\\
    \hline 
        318 & 2.1e-02 & -- &  2.1e-02 & -- & 2.1e-02 & -- & 2.1e-02 & --\\
       1644  & 1.8e-03 & 3.6 &1.8e-03 & 3.6 & 1.8e-03 & 3.6& 1.8e-03 & 3.6 \\
        6144  & 2.2e-04 & 3.0  & 2.2e-04 & 3.0 & 2.2e-04 & 3.0 & 2.2e-04 & 3.0\\
        24846  & 2.5e-05 & 3.2 & 2.5e-05 & 3.1 & 2.5e-05 & 3.1 & 2.5e-05 & 3.2\\
  \end{tabular}
\end{table}
\Cref{tab:E1_Stokes-k2,tab:E1_Darcy-k2} for $k_f=2$, $k_c=1$, and \Cref{tab:E1_Stokes-k3,tab:E1_Darcy-k3} for $k_f=3$, $k_c=2$ show that in the Stokes and Darcy regions $u_h$ and $p_h$ both converge optimally in the $L^2$-norm with orders $k_f+1$ and $k_f$, respectively. This is consistent with our theoretical convergence rate in Corollary~\ref{cor:velerror} that predicts at least suboptimal rates for the velocity. The right most columns in these tables demonstrate pointwise mass conservation. 

Furthermore, even though the magnitude of the pressure error changes dramatically as we change the values of $\kappa$ and $\mu$, there is no significant change in the velocity errors. This is more pronounced in the case where both the permeability and the viscosity are small ($10^{-3}$ and $10^{-6}$).
This confirms that the velocity error bounds in Theorem~\ref{thm:vel-error} and Corollary~\ref{cor:velerror} are independent of the pressure error.

We furthermore observe from \Cref{tab:E1_concentration-k2} and \Cref{tab:E1_concentration-k3} that $c_h$ converges optimally in the $L^2$-norm with order $k_c+1$. This supports our result \cref{eq:L2-c} that shows at least suboptimal convergence. 

\subsection{Example 2} 
\label{sec:example2}
We now consider the fully coupled problem by incorporating the influence of the velocity solution on the dispersion/diffusion tensor and the dependence of the viscosity on the concentration solution. The source terms and boundary conditions for the Stokes/Darcy-transport problem \eqref{eq:system} are chosen such that the exact solution is given by \cref{eq:exact_solution_SDT}. We define the diffusion dispersion tensor in $\Omega$ and the viscosity according to
\begin{equation}
\label{eq:quarterpowermixing}   
  \widetilde{D}(u) = \begin{bmatrix}
    1+u_1^2 & 0\\
    0 & 1+u_2^2
    \end{bmatrix},
    \qquad
    \mu(c) = \mu_0\Big[\Big(\frac{\mu_0}{\mu_1}\Big)^{1/4}c+(1-c) \Big]^{-4},
\end{equation}
where we remark the the viscosity is defined as the quarter-power mixing rule \cite{Lohrenz:1964} with
where $\mu_0 = 0.9$, $\mu_1 = 1.3$. 

We use $k_f=3$, $k_c=2$, and BDF3 time stepping with $\Delta t = 0.1h^3/4$. We present numerical results for $\kappa=10^3, 1, 10^{-3}$. We observe from \Cref{tab:E2_Stokes-k3,tab:E2_Darcy-k3} that when $\mu$ is changed from a constant to a concentration dependent function the rate of convergence reduces from $k_f+1$ to a value between $k_f$ and $k_f + 1$. 
 This is due to our choice $k_c=k_f-1$ to achieve compatibility and is consistent with our a priori estimates \cref{eq:estimate_zeta_u_n_s,eq:estimate_zeta_u_n_d} in the energy norm which imply that the rate of convergence of the velocity approximation is polluted by the concentration approximation. Indeed, from \Cref{tab:E2_concentration-k3} we observe that $c_h$ converges in the $L^2$-norm with order $k_f$.
 Therefore, for the velocity we expect an order of at least $k_f-1$ in the energy norm and $k_f$ in the $L^2$-norm. From the right most columns in \Cref{tab:E2_Stokes-k3,tab:E2_Darcy-k3} we observe that the discretization is exactly mass conserving. 
%
\begin{table}[!ht]
  \caption{Errors and rates of convergence at final time $T=0.1$ for $u_h$ and $p_h$ in the Stokes region $\Omega^s$ for the test case in \Cref{sec:example2} using $k_f=3$, $k_c=k_f-1$, and BDF3 time stepping with $\Delta t= 0.1h^3/4$.}  
  \label{tab:E2_Stokes-k3}
  \centering
   \begin{tabular}{ccccccc}
 \textit{h} & dofs   &  $\|u_h-u\|_{\Omega^s}$ &    rate   &   $\|p_h-p\|_{\Omega^s}$  &   rate    &  $\|\nabla\cdot u_h\|_{\Omega^s}$ \\ \hline
    \multicolumn{7}{l}{$\kappa=10^{3}$} \\
    \hline
        1/4 & 1161 &  9.2e-05 & --& 8.1e-03 & --& 2.3e-15 \\
        1/8 & 5993 & 3.7e-06 & 4.6 & 3.9e-04 & 4.4 & 5.6e-15 \\
        1/16 & 22081 &  2.4e-07 & 4.0 & 4.1e-05 & 3.2 & 6.1e-15\\
        1/32 &  90241& 2.1e-08 & 3.5 & 4.7e-06 & 3.1 & 8.5e-15  \\
    \hline 
    \multicolumn{7}{l}{$\kappa=1$} \\
    \hline
        1/4 &  1161 & 9.6e-05 & -- & 8.1e-03 & -- & 2.6e-15 \\
        1/8 &  5993 & 3.8e-06 & 4.7 & 3.9e-04 & 4.4 & 5.5e-15 \\
        1/16 & 22081 & 2.3e-07 & 4.1 & 4.1e-05 & 3.2 & 6.1e-15 \\
        1/32 & 90241 & 1.8e-08 & 3.6 & 4.7e-06 & 3.1 & 8.6e-15\\
    \hline
    \multicolumn{7}{l}{$\kappa=10^{-3}$} \\
    \hline
        1/4 & 1161 & 1.8e-03 & -- & 6.9e-01 & -- & 6.3e-15\\
        1/8 & 5993 &  5.0e-05 & 5.1 & 5.4e-02 & 3.7 & 7.1e-15\\
        1/16& 22081 &  3.6e-06 & 3.8 & 6.3e-03 & 3.1 & 7.4e-15 \\
        1/32& 90241&4.0e-07 & 3.2 & 7.7e-04 & 3.0 & 9.2e-15  \\
    \end{tabular}
\end{table}
\begin{table}[!ht]
  \caption{Errors and rates of convergence at final time $T=0.1$ for $u_h$ and $p_h$ in the Darcy region $\Omega^d$ for the test case in \Cref{sec:example2} using $k_f=3$, $k_c=k_f-1$, and BDF3 time stepping with $\Delta t= 0.1h^3/4$.}  
  \label{tab:E2_Darcy-k3}
  \centering
    \begin{tabular}{ccccccc}
      \hline
 \textit{h} & dofs   &  $\|u_h-u\|_{\Omega^d}$ &    rate   &   $\|p_h-p\|_{\Omega^d}$  &   rate    &  $\|\Pi_Q(\nabla\cdot (u_h -u))\|_{\Omega^d}$ \\ \hline
    \multicolumn{7}{l}{$\kappa=10^3$} \\
    \hline
        1/4 & 1161 & 3.5e-03 & -- & 6.6e-07 & -- & 4.2e-12\\
        1/8 &  5993 & 3.4e-04 & 3.4 & 4.2e-08 & 4.0 & 8.4e-13\\
        1/16 & 22081 & 3.8e-05 & 3.1 & 5.3e-09 & 3.0 & 2.9e-12\\
        1/32 & 90241 & 4.4e-06 & 3.1 & 6.3e-10 & 3.1 & 1.2e-11 \\
    \hline
    \multicolumn{7}{l}{$\kappa=1$} \\
    \hline
        1/4 & 1161 & 3.5e-03 & --  & 6.6e-04 & -- & 4.2e-12\\
        1/8 &  5993 & 3.4e-04 & 3.4 & 4.2e-05 & 4.0 & 9.0e-13\\
        1/16& 22081 & 3.8e-05 & 3.1 & 5.3e-06 & 3.0 & 3.0e-12\\
        1/32& 90241 & 4.4e-06 & 3.1 & 6.3e-07 & 3.1 & 1.2e-11\\
    \hline
    \multicolumn{7}{l}{$\kappa=10^{-3}$} \\
    \hline
        1/4 & 1161 & 3.9e-03 & -- & 6.0e-01 & -- & 4.2e-12  \\
        1/8 & 5993 & 3.4e-04 & 3.5 & 4.1e-02 & 3.9 & 9.6e-13 \\
        1/16 & 22081 &  3.8e-05 & 3.1 & 5.2e-03 & 3.0 & 3.5e-12 \\
        1/32 & 90241 & 4.4e-06 & 3.1 & 6.2e-04 & 3.1 & 1.5e-11 \\
  \end{tabular}
\end{table}
%


%
\begin{table}[!ht]
  \caption{Errors and rates of convergence at final time $T=0.1$ for $c_h$ in $\Omega$, on a mesh with $h=1/4, 1/8, 1/16, 1/32$, for the test case in \Cref{sec:example2} using $k_f=3$, $k_c=k_f-1$, and BDF3 time stepping with $\Delta t= 0.1h^3/4$.}  
  \label{tab:E2_concentration-k3}
  \centering
    \begin{tabular}{c|cc|cc|cc}
\multicolumn{1}{c}{}&\multicolumn{2}{c|}{$\kappa=10^3$} &\multicolumn{2}{c|}{$\kappa=1$} &\multicolumn{2}{c}{$\kappa=10^{-3}$}   \\\hline
        dofs & $||c-c_h||_{\Omega}$ & rate & $||c-c_h||_{\Omega}$& rate & $||c-c_h||_{\Omega}$ & rate \\\hline
        318 & 1.4e-02 & --& 1.4e-02 & -- & 1.4e-02 & --\\
       1644  & 1.2e-03 & 3.6 & 1.2e-03 & 3.6 &  1.2e-03 & 3.6\\
        6144  & 1.3e-04 & 3.2 & 1.3e-04 & 3.2 & 1.3e-04 & 3.2\\
        24846  & 1.3e-05 & 3.3 & 1.3e-05 & 3.3 & 1.3e-05 & 3.3\\
  \end{tabular}
\end{table}

\subsection{Example 3}
\label{sec:example3}
\begin{figure}[ht!]  
    \centering
       \includegraphics[width=0.4\textwidth]{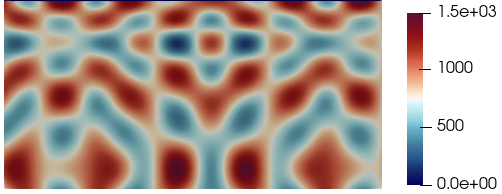}
        \caption{The permeability field in $\Omega^d=[0,1]\times [0,0.5]$ defined by \cref{eq:permeability}.}
            \label{fig:perm}
\end{figure}

In this final example, we simulate a more realistic problem similar to \cite[Section 6.2]{Cesmelioglu:2020} in which the permeability field in the Darcy region is highly heterogeneous. For this, let the boundary of the Stokes region be partitioned as
$\Gamma^s = \Gamma^s_1\cup\Gamma_2^s\cup\Gamma_3^s$ where
\[
    \Gamma^s_1 :=\cbr{ x \in\Gamma^s:\ x_1=0 }, \quad 
    \Gamma^s_2 :=\cbr{ x \in\Gamma^s:\ x_1=1 }, \quad 
    \Gamma^s_3 :=\cbr{ x \in\Gamma^s:\ x_2=1 }.
\] 
Similarly, let
$\Gamma^d = \Gamma^d_1\cup\Gamma_2^d$ where
\[
    \Gamma^d_1 :=\{ x \in\Gamma^d\:\ x_1=0 \ \text{or} \ x_1=1\}, \quad
    \Gamma^d_2 :=\{ x \in\Gamma^d\:\ x_2=0\}.
\]
We impose the
following boundary conditions:
\begin{align*}
  u &= (x_2(3/2-x_2)/5, 0) & \text{on }\, \Gamma_1^s, \\
  (-2\mu\varepsilon(u) + p\mathbb{I})n &= 0 & \text{on }\, \Gamma_2^s, \\
  u\cdot n &= 0 \ \text{and}\, (-2\mu\varepsilon(u) + p\mathbb{I})^t = 0 & \text{on }\, \Gamma_3^s, \\
  u\cdot n &= 0 & \text{on }\, \Gamma_1^d, \\
  p &= -0.05 & \text{on }\, \Gamma_2^d.
\end{align*}
The first boundary condition on the left boundary $\Gamma_1^s$ of $\Omega^s$ imposes a parabolic velocity profile.
We set the permeability to
\begin{equation}
\label{eq:permeability}
    \kappa = 700 (1 + 0.5(\sin(10\pi x_1) \cos(20\pi x_2^2) + \cos^2(6.4\pi x_1)\sin(9.2\pi x_2 ))) + 100,
\end{equation}
a plot of which is given in \Cref{fig:perm}. 
The viscosity is defined by the quarter-power mixing rule as in \cref{eq:quarterpowermixing}.
The other parameters are set as $\mu = 0.1$, $\alpha = 0.5$, $k_f = 3$, $h = 1/80$, $\Delta t = 10^{-3}$, $T = 15$, and the source/sink terms are set to zero.
In the Darcy region $\Omega^d$ the porosity is set to $\phi=0.4$. The dispersion/diffusion tensor is defined as
\begin{equation*}
  \widetilde{D}(u) =
  \begin{cases}
    \delta I, & \text{ in } \Omega^s,
    \\
    \phi d_m \mathbb{I} + d_l|u|\mathbb{T} + d_t|u|(\mathbb{I} - uu^T/|u|^2),
    & \text{ in } \Omega^d,    
  \end{cases}
\end{equation*}
where $d_l, d_t$, and $d_m$ represent longitudinal and transverse dispersivities and
the molecular diffusivity, respectively, 
and $u^T$ is the
transpose of the vector $u$. 
Under the condition $d_l\geq d_t$ (which is usually the
case), $D(u)$ satisfies the assumptions \cref{eq:D_min,eq:upperboundDuabs,eq:DLipschitz} (see, for example, \cite{Douglas:1983}, \cite[Lemmas 4.3, 4.4]{Sun:2002}, and \cite[Lemma 1.3]{Riviere:2011}).
In this numerical experiment, we choose $\delta = 10^{-6}$, $d_m = 10^{-5}$, $d_l = 10^{-5}$, and $d_t = 10^{-5}$. The initial velocity is set to zero while the initial concentration of the plume of contaminant is defined as
\begin{equation*}
  c_0(x) =
  \begin{cases}
    0.95 & \text{if}\ \sqrt{ (x_1-0.2)^2 + (x_2-0.7)^2 } < 0.1, \\
    0.05 & \text{otherwise}.
  \end{cases}
\end{equation*}
\begin{figure}  \centering
            \includegraphics[width=0.35\textwidth]{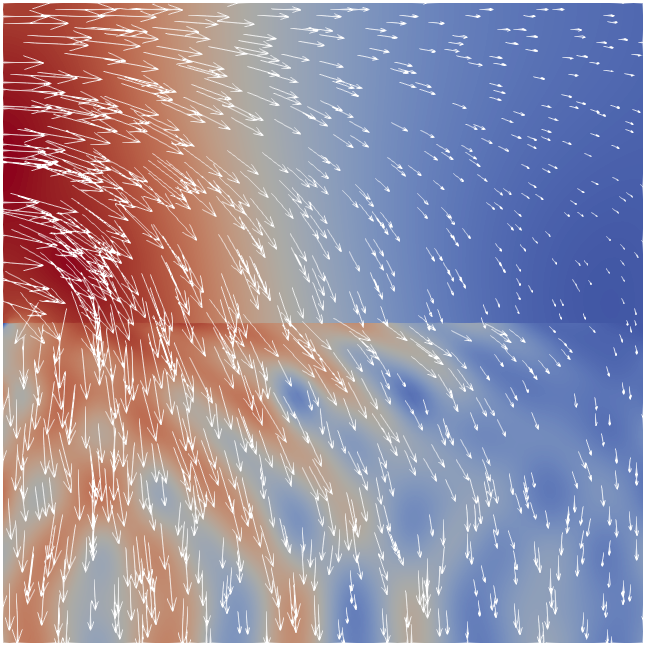}
            \includegraphics[width=0.35\textwidth]{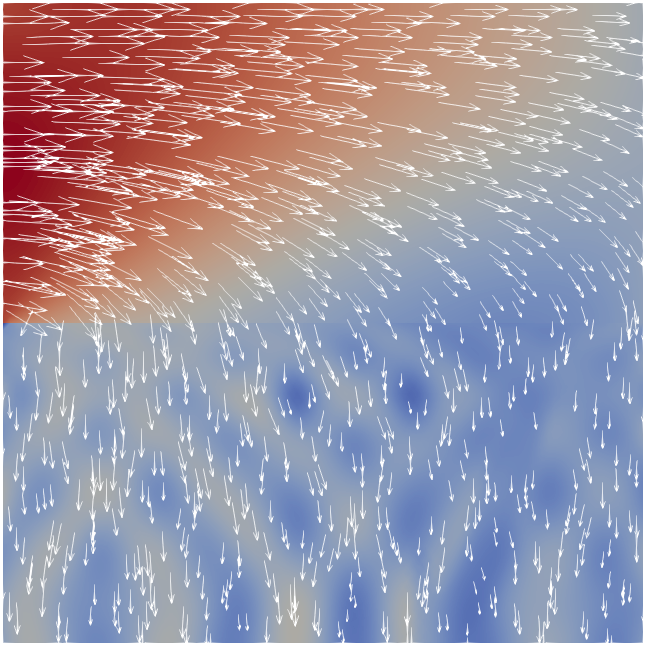}
            \includegraphics[width=0.1\textwidth]{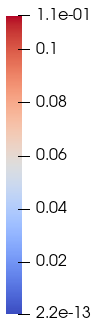}
 \caption{The velocity field after one time step (left) and at the final time (right)  for the example in \Cref{sec:example3}. The color represents the magnitude of the velocity.}
    \label{fig:initfinalvel}
\end{figure}
\begin{figure}
  \centering \subfloat[$t=\Delta t$]{\includegraphics[width=0.34\textwidth]{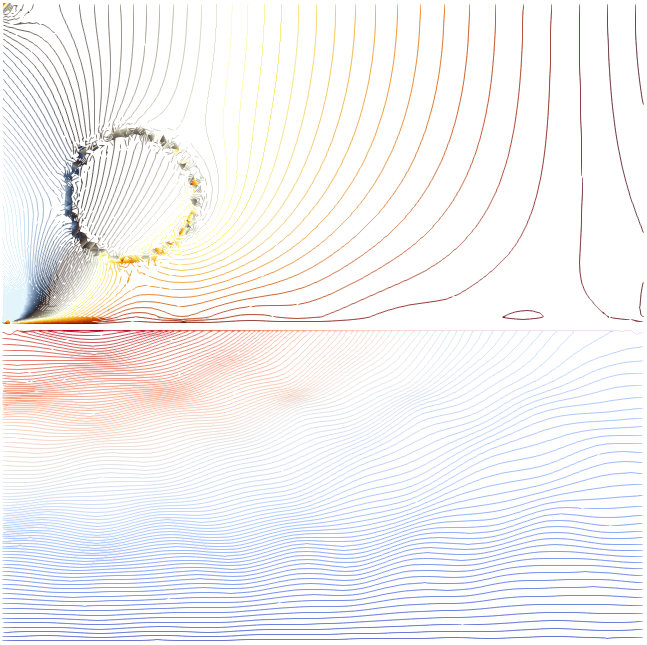}}
  \quad
  \subfloat[$t=3$ ]{\includegraphics[width=0.34\textwidth]{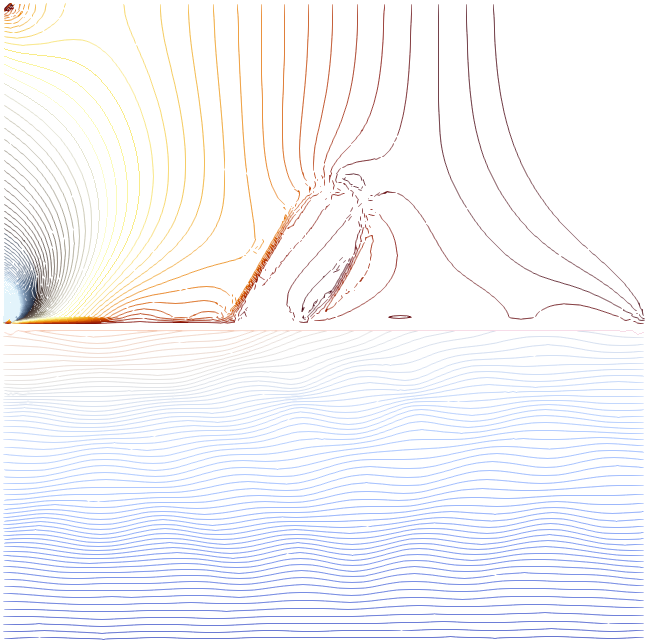}}
  \subfloat{\includegraphics[width=0.06\textwidth]{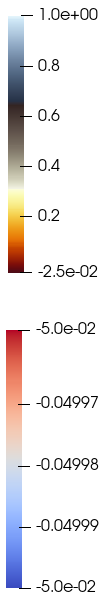}}
  \\
  \subfloat[$t=6$]{\includegraphics[width=0.34\textwidth]{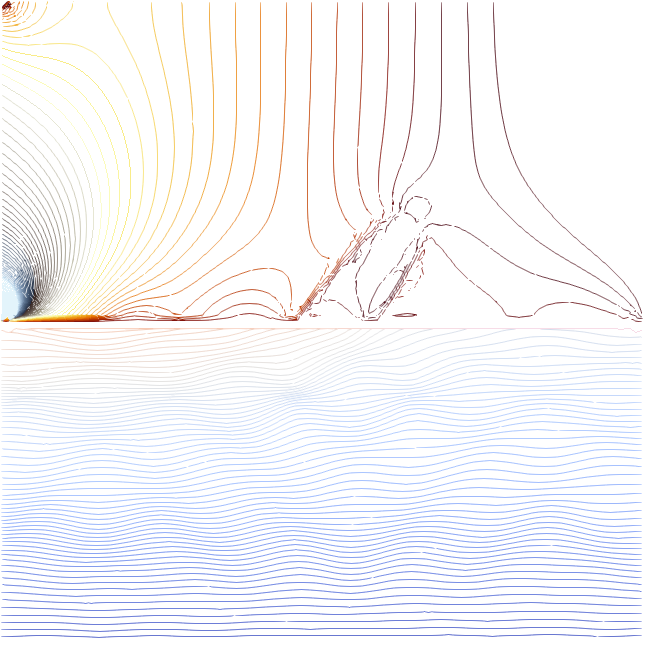}}
  \quad \subfloat[$t=9$ ]{\includegraphics[width=0.34\textwidth]{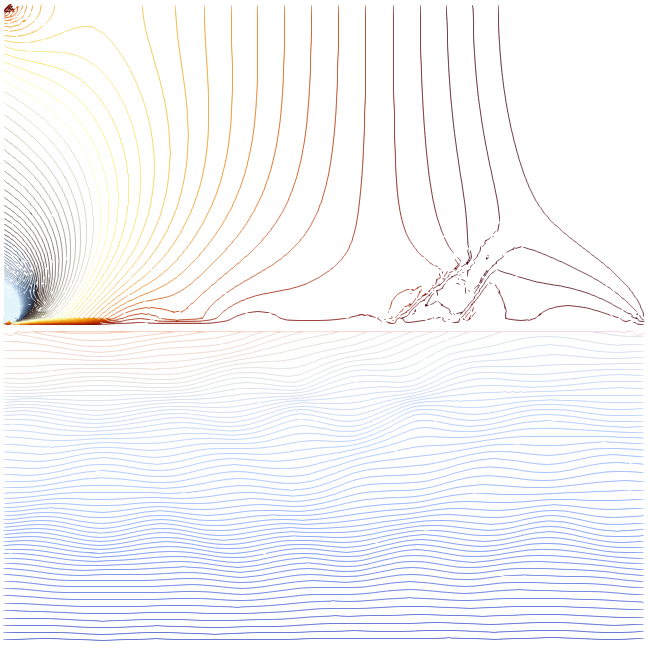}}
  \subfloat{\includegraphics[width=0.06\textwidth]{pressure-scale.png}}
  \\
  \subfloat[$t=12$. ]{\includegraphics[width=0.34\textwidth]{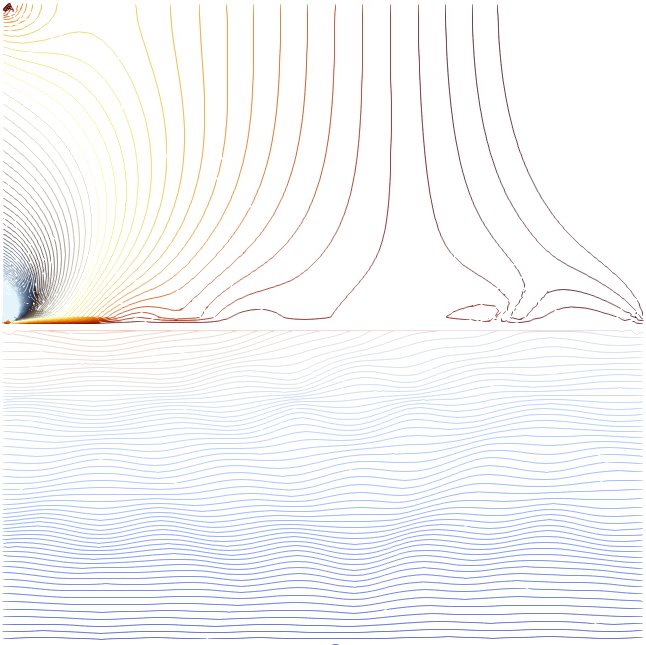}}
  \quad \subfloat[$t=15$ ]{\includegraphics[width=0.34\textwidth]{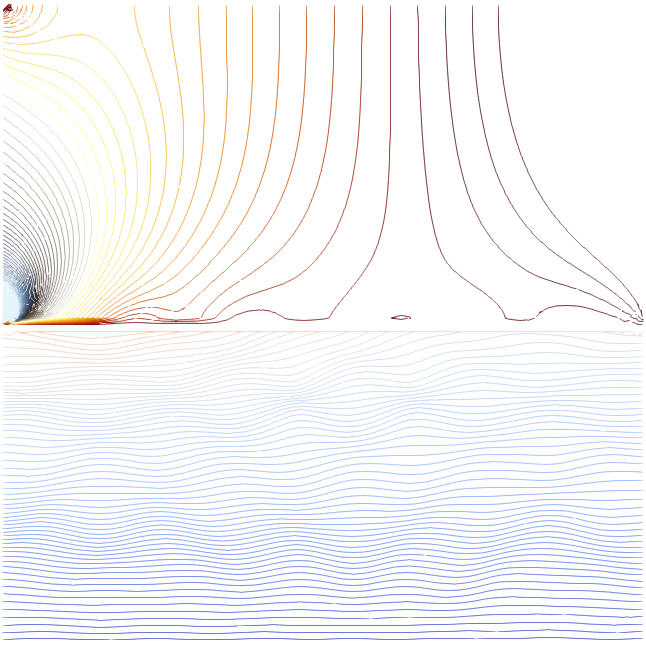}}
  \subfloat{\includegraphics[width=0.06\textwidth]{pressure-scale.png}}
  \caption{Pressure contours at times $t= \Delta t, 3, 6, 9, 12, 15$ for the example in \Cref{sec:example3}. The color represents the concentration values.}
  \label{fig:pres}
\end{figure}

\begin{figure}
  \centering \subfloat[
  $t=\Delta t$. ]{\includegraphics[width=0.34\textwidth]{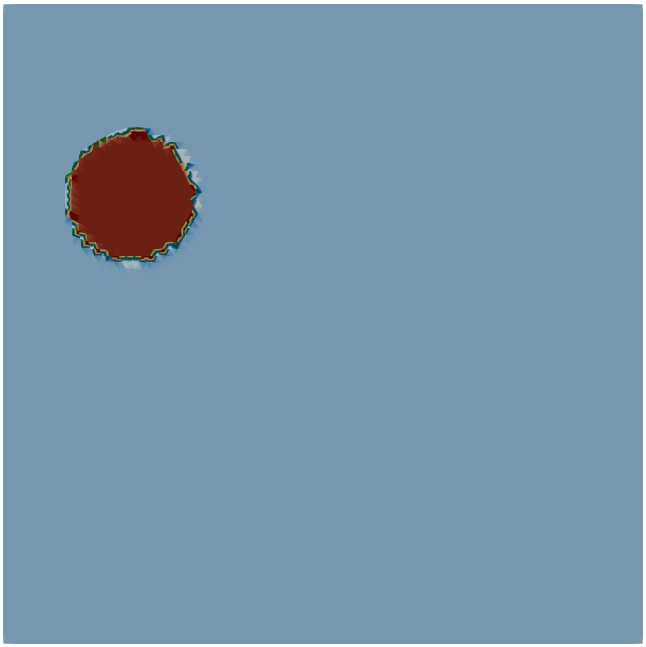}}
  \quad
  \subfloat[$t=3$. ]{\includegraphics[width=0.34\textwidth]{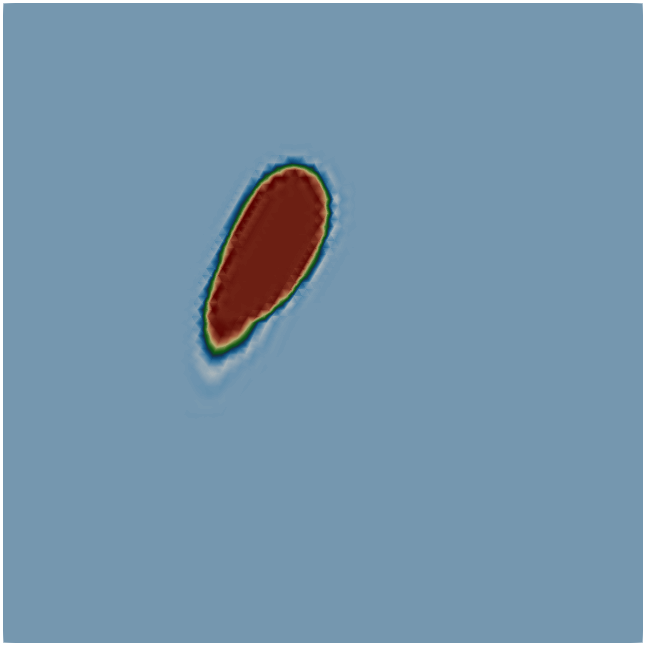}}
  \subfloat{\includegraphics[width=0.06\textwidth]{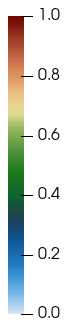}}
  \\
  \subfloat[
  $t=6.0$. ]{\includegraphics[width=0.34\textwidth]{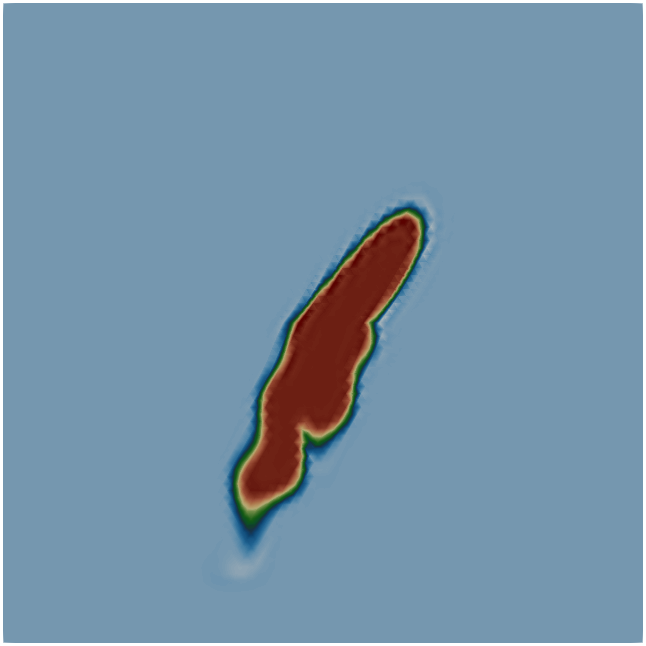}}
  \quad \subfloat[
  $t=9$. ]{\includegraphics[width=0.34\textwidth]{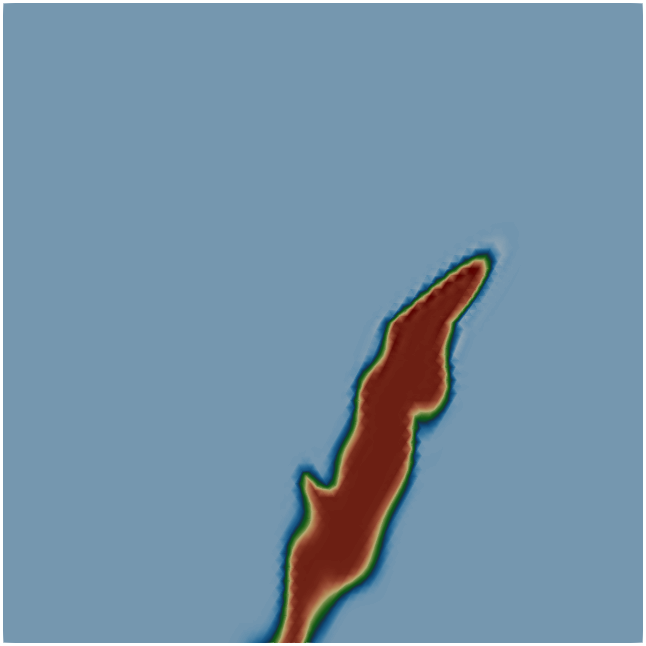}}
  \subfloat{\includegraphics[width=0.06\textwidth]{concentration-scale.png}}
  \\
  \subfloat[
  $t=12$. ]{\includegraphics[width=0.34\textwidth]{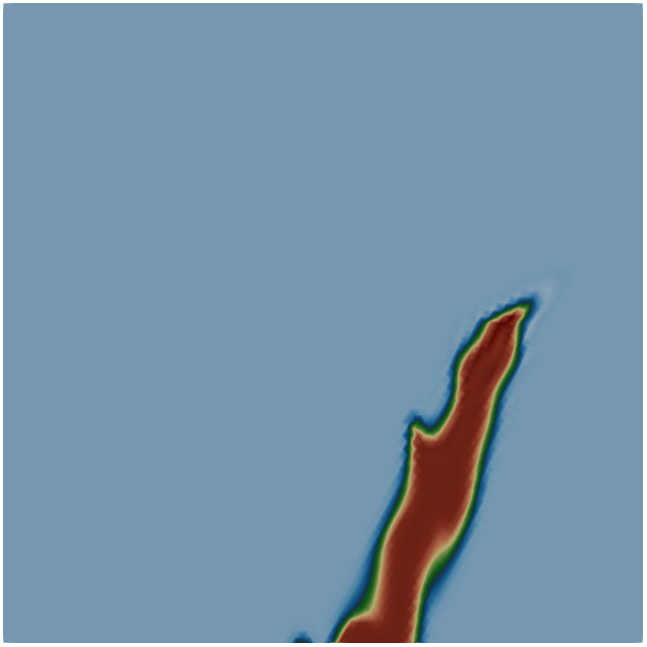}}
  \quad \subfloat[
  $t=15$. ]{\includegraphics[width=0.34\textwidth]{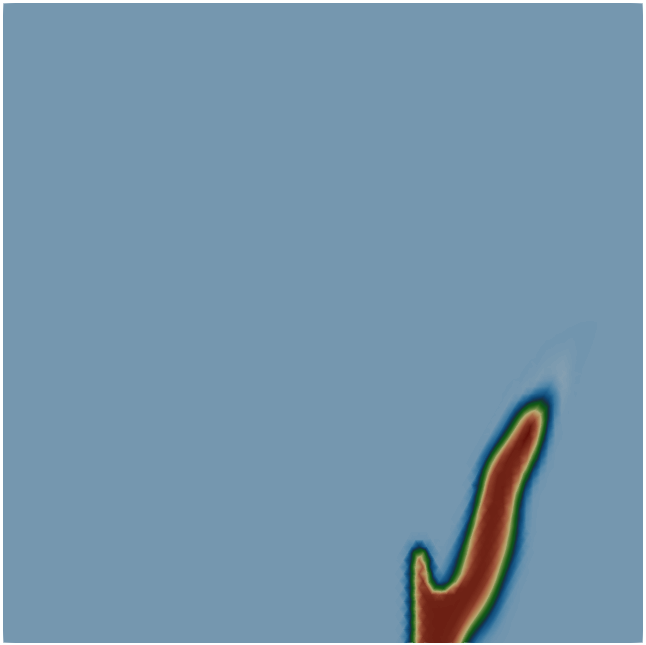}}
  \subfloat{\includegraphics[width=0.06\textwidth]{concentration-scale.png}}
  \caption{The plume of contaminant at times $t= \Delta t, 3, 6, 9, 12, 15$ for the example in \Cref{sec:example3}. The color represents the concentration values.}
  \label{fig:plume}
\end{figure}

We compute the solution using BDF3 time stepping. \Cref{fig:initfinalvel} shows the computed velocity field after one time step and at the final time.  In the Darcy region $\Omega^d$, the flow field avoids areas with low permeability as expected. \Cref{fig:pres} shows the pressure contours at various times which demonstrates the effect of the concentration on the pressure around the plume of contaminants, especially in the Stokes region. \Cref{fig:plume} presents the plume of contaminant spreading through the surface water region and infiltrating the porous medium. We plot the solution 6 different instances in time. The contaminant plume stays compact while in the surface water region. Once it reaches the subsurface region it spreads out following a path dictated by the heterogeneous permeability structure of the porous medium.

\section{Conclusions}
\label{sec:conclusions}
In this paper, we introduced and analyzed a fully discrete sequential method for the fully coupled Stokes/Darcy--transport problem. The spatial discretization uses the HDG method which is higher-order accurate, strongly mass conservative, and compatible.  We remark that the analysis also easily extends to the EDG-HDG method considered in \cite{Cesmelioglu:2020,Cesmelioglu:2021}. The sequential method discussed in the article linearizes the problem by time-lagging the concentration and decoupling the Stokes/Darcy and transport subproblems. We proved well-posedness and obtained a priori estimates in the energy norm. Finally, we presented numerical results demonstrating mass conservation and robustness with respect to varying permeability and optimal convergence in the $L^2$-norm for one-way coupling.

\bibliographystyle{abbrv}
\bibliography{references}

\end{document}